\pdfoutput=1
\documentclass[11pt]{amsart}
\usepackage{txfonts}
\usepackage{amsmath,amsthm,amscd,amssymb,mathtools}
\usepackage{graphicx}
\usepackage{tikz}
\usepackage[colorlinks=true,linkcolor=blue,citecolor=blue,urlcolor=blue,backref=page]{hyperref}
\renewcommand*{\backrefalt}[4]{%
  \ifcase #1 \or $\uparrow$\,#2\else $\uparrow$\,#2\fi}

\usepackage{enumitem}
\usepackage{orcidlink}

\numberwithin{equation}{section}

\newtheorem{theorem}{Theorem}[section]
\newtheorem{corollary}[theorem]{Corollary}
\newtheorem{lemma}[theorem]{Lemma}
\newtheorem{proposition}[theorem]{Proposition}
\newtheorem{theoremx}{Theorem}

\theoremstyle{definition}

\newtheorem{remark}[theorem]{Remark}
\newtheorem{question}[theorem]{Question}

\newcommand{\Z}{\mathbb{Z}}
\newcommand{\R}{\mathbb{R}}
\newcommand{\C}{\mathbb{C}}
\newcommand{\Q}{\mathbb{Q}}
\newcommand{\A}{\mathcal{A}}
\newcommand{\CP}{\mathbb{CP}}
\newcommand{\la}{\langle}
\newcommand{\ra}{\rangle}

\newcommand{\BS}{\operatorname{BS}}
\newcommand{\Linf}{L_\infty}

\newcommand{\sh}{\operatorname{sh}}
\newcommand{\cS}{\mathcal{S}}
\newcommand{\V}{\mathcal{V}}
\newcommand{\RR}{\mathcal{R}}
\newcommand{\lk}{\operatorname{lk}}
\newcommand{\nc}[1]{\langle\!\langle #1\rangle\!\rangle}

\DeclareMathOperator{\Zc}{Z}
\DeclareMathOperator{\Hom}{Hom}
\DeclareMathOperator{\PSL}{PSL}
\DeclareMathOperator{\SL}{SL}
\DeclareMathOperator{\id}{id}
\DeclareMathOperator{\pr}{pr}
\DeclareMathOperator{\tr}{tr}
\renewcommand{\Re}{\operatorname{Re}}
\renewcommand{\Im}{\operatorname{Im}}

\title[Twisted pencils and non-Hopfian Hurwitz curves]%
{Twisted pencils of quartics, braid monodromy factorizations, and non-Hopfian Hurwitz curves}

\author[Alexander~I.~Suciu]{Alexander~I.~Suciu\,\orcidlink{0000-0002-5060-7754}}
\address{Department of Mathematics, Northeastern University, Boston, MA 02115, USA}
\email{\href{mailto:a.suciu@northeastern.edu}{a.suciu@northeastern.edu}}
\urladdr{\href{https://suciu.sites.northeastern.edu}{suciu.sites.northeastern.edu}}

\date{November 7, 2026}



\begin{document}
\begin{abstract}
We study which groups arise as fundamental groups of complements of affine
Hurwitz curves.  In every such group the product of the meridians is central
and maps to the vector of component degrees in homology; this rules out, for
instance, the Baumslag--Solitar groups $\BS(p,q)$ with
$\lvert p\rvert\neq\lvert q\rvert$, and the groups of links with at least two
components and vanishing linking numbers.  O.~V.~Kulikova circumvented 
this obstruction for $\BS(2,3)$ by passing to the non-Hopfian group $\BS(2,3)\times\Z$, 
which Vik.~S.~Kulikov's doubling algorithm realizes in degree $768$.  
We give an explicit braid monodromy factorization of degree $12$ with 
this group.  The corresponding curve is not isotopic to an algebraic curve, 
yet it is obtained from a union of three quartics in a pencil by a single twist 
along a circle.  More generally, for every closed braid $\hat\beta$ we realize
$\pi_1(S^3\setminus\hat\beta)\times\Z$ by an explicit Hurwitz curve, again 
a twist of a union of quartics in a pencil.  In the other direction, a Hurwitz
curve with group $\BS(k,k+1)\times\Z$ has degree at least $4$, as every 
Hurwitz curve consisting of a conic and a line is H-isotopic to an algebraic 
one, after separation of its critical points.  Finally, we show that certain
semiholomorphic germs, among them those realizing the generalized 
Borromean rings, are braided surfaces but not Hurwitz curve germs.
\end{abstract}

\maketitle
\setcounter{tocdepth}{1}
\tableofcontents

\section{Introduction}
\label{sec:intro}

\subsection{Hurwitz curves and their groups}
\label{subsec:intro-hurwitz}
Hurwitz curves, introduced by Moishezon \cite{Moishezon-conm94} under the
name of semi-algebraic curves, imitate plane algebraic curves with respect
to a fixed pencil of lines: they are surfaces in $\CP^2$ that are smooth and
positively transverse to the lines of the pencil away from finitely many
points, near which they are complex analytic.  With negative nodes allowed,
they are the natural branch curves of the quasiholomorphic coverings of $\CP^2$
by symplectic $4$-manifolds constructed by Auroux and Katzarkov
\cite[Thm.~1]{Auroux-Katzarkov-inv00}, and their braid monodromy factorizations
form the combinatorial skeleton of that theory
\cite[Thm.~2]{Auroux-Katzarkov-inv00}, \cite{Kharlamov-Kulikov-izv03,Kulikov-usp07}.  
By the Zariski--van Kampen theorem the group $\pi_1(\C^2\setminus H)$ of an affine 
Hurwitz curve is read off from its braid monodromy, and Kulikov \cite[Thm.~6.1]{Kulikov-izv04a} 
identified the class of all such groups: they are exactly the \emph{Hurwitz C-groups}, the
C-groups admitting a presentation in which the product $l_m=x_1\cdots x_m$
of the C-generators is central.

In the algebraic setting Zariski \cite{Zariski-71} asked whether
$\pi_1(\CP^2\setminus C)$ is always residually finite; the question is open.
For Hurwitz curves the answer is negative: O.~V.~Kulikova \cite{Kulikova-izv05}
showed that the non-Hopfian group $\BS(2,3)\times\Z$ is a Hurwitz C-group, so
that, by Kulikov's realization theorem, it is $\pi_1(\C^2\setminus H)$ for a Hurwitz
curve $H$ with two irreducible components.  The realization is algorithmic
but not explicit.  The purpose of this paper is to explain precisely which
structural features force the passage to $\BS(2,3)\times\Z$, to make
Kulikova's example completely explicit, in degree $12$, to give a
construction of the same kind that realizes $G\times\Z$ for every link
group $G$ in a degree depending only on the number of strings, and to bound
the degree of such realizations from below.

\subsection{The centrality obstruction}
\label{subsec:intro-central}
The first observation is elementary, but it constrains everything that
follows.

\begin{proposition}
\label{prop:intro-central}
Let $\overline H\subset\CP^2$ be a Hurwitz curve of degree $m$ with
irreducible components of degrees $d_1,\dots,d_r$, let $\Linf$ be a line of
the pencil in general position with respect to $\overline H$, and let
$G=\pi_1(\C^2\setminus H)$, where $H=\overline H\setminus\Linf$.  Let
$x_1,\dots,x_m\in G$ be the meridians of $H$ given by a geometric basis of
loops in a generic line of the pencil, and let $l_m=x_1\cdots x_m$ be the
class of a loop in that line encircling all $m$ points of $H$.  Then $l_m$ is
central in $G$, and its image in $H_1(G;\Z)\cong\Z^r$ is
$(d_1,\dots,d_r)\neq0$.
\end{proposition}

The centrality of $l_m$ is noted in \cite[\S6]{Kulikov-izv04a} and
\cite[\S0]{Kulikova-izv05}; the point here is its image in homology.  
Thus the center of a Hurwitz curve group has non-zero image in
rational homology, whatever the number of components.  This excludes
every group with trivial center, for instance $\BS(p,q)$ with
$\lvert p\rvert\neq\lvert q\rvert$ or the fundamental group of a
hyperbolic $3$-manifold of finite volume (Corollary~\ref{cor:trivcenter}), and every 
link group with at least two components and vanishing pairwise linking numbers
(Proposition~\ref{prop:linking}); for knot groups it leaves exactly the
torus knot groups, as Kulikov observed \cite[Thm.~6.3]{Kulikov-izv04a}, and for
right-angled Artin groups exactly those of cones (Remark~\ref{rem:raag}).
Proposition~\ref{prop:linking} applies to the generalized
Borromean rings $L_n$ of the companion paper \cite{Suciu-borro}, which are
links of isolated singularities of semiholomorphic polynomial maps; the
argument uses Milnor's invariants, not the hyperbolicity of their
complements.  In the other direction, Kulikova's product trick
\cite[\S2]{Kulikova-izv05} shows that $G\times\Z$ is a Hurwitz C-group for
every C-group $G$ (Proposition~\ref{prop:product}).

\subsection{An explicit non-Hopfian Hurwitz curve}
\label{subsec:intro-explicit}
Write $\sigma_1,\dots,\sigma_{m-1}$ for the standard generators of the braid group
$B_m$, $\Delta_m^2$ for the full twist, $\cS_m$ for the semigroup of
factorizations in $B_m$ (finite sequences of braids, up to Hurwitz moves),
$G(s)$ for the group of a factorization $s$, the quotient of the free group
$F_m=\la x_1,\dots,x_m\ra$ by the relations $x_j=u(x_j)$ for all factors $u$
of $s$, $d(\cdot,\cdot,\cdot,\cdot)$ for Kulikov's doubling of
factorizations \cite[Thm.~3.2]{Kulikov-izv04a}, and $\lambda(b)$ for
simultaneous conjugation of all factors by $b$ (see \S\ref{sec:prelim} for
the precise conventions).

\begin{theoremx}
\label{thm:explicit}
Let $s=d(\Delta_3^2,\Delta_3^2,\Delta_3^2,\Delta_3^2)$, a factorization of
$\Delta_6^2$ with ten factors, and put
\[
g_1=\sigma_1^{-1}\sigma_2\sigma_3^{-1}\sigma_2,\qquad
g_2=\sigma_2\,\sigma_3^{-1}\sigma_4^{-2}\sigma_5^{-2}\sigma_4^{-1}\,
\sigma_3\,\sigma_2,\qquad b_k=g_k\,\sigma_1\,g_k^{-1}\in B_6,
\]
and $T_k=d\bigl(\lambda(b_k)(s),s,s,s\bigr)$.  Then, for $k=1,2$:
\begin{enumerate}[label=(\roman*),itemsep=1pt]
\item \label{a1}
$T_k$ is a factorization of $\Delta_{12}^2$ into $52$ factors: $16$
conjugates of $\Delta_3^2$ and $36$ conjugates of $\sigma_1$.
\item \label{a2}
The Zariski--van Kampen group of $T_k$ is
$G(T_k)\cong\BS(k,k+1)\times\Z$, the isomorphism sending the meridian $x_j$ to
$t$, $ta$, or $(t^2a)^{-1}z$ according as $j\equiv1,2,0\pmod 3$, where
$\BS(k,k+1)=\la a,t\mid t^{-1}a^kt=a^{k+1}\ra$ and $\Z=\la z\ra$.
\item \label{a3}
$T_k$ is the braid monodromy factorization of a Hurwitz curve
$\overline H_k\subset\CP^2$ of degree $12$, whose singular points are $16$
ordinary triple points; $\overline H_k$ has two irreducible components, of
degree $8$ and geometric genus $5$, and of degree $4$ and genus $3$, 
the latter smooth; each triple point is a node of the first component through 
which the second passes; and
\[
\pi_1(\C^2\setminus H_k)\cong\BS(k,k+1)\times\Z,\qquad
\pi_1(\CP^2\setminus \overline H_k)\cong\BS(k,k+1)\times\Z/4.
\]
\end{enumerate}
For $k=2$ both groups are non-Hopfian, hence not residually finite.
\end{theoremx}

For comparison, Kulikov's algorithm \cite[Thm.~3.6]{Kulikov-izv04a}
produces for the same group a curve of degree $3\cdot2^8=768$
(\S\ref{subsec:768}); the saving comes from writing the band $b_k$ directly
in $B_6$ and from a one-doubling version of Kulikov's band lemma
(Lemma~\ref{lem:band}).  By Theorem~\ref{thm:alex} the curves
$\overline H_k$ are not isotopic to algebraic curves, but
Theorem~\ref{thm:model} describes each of them as an explicit algebraic curve
of degree $12$ with one quarter twisted along a circle.

We do not know whether degree~$12$ is the smallest possible.  It is the
smallest among twisted doublings of iterated doublings of full twists
(Proposition~\ref{prop:min12}), and every Hurwitz curve with group
$\BS(k,k+1)\times\Z$ has degree at least~$4$ (Corollary~\ref{cor:deg3}).  For
the latter we show that every Hurwitz curve of degree~$3$ consisting of a
conic and a line is H-isotopic, after separation of its critical points, to
an algebraic curve (Proposition~\ref{prop:deg3}).

\subsection{Twisted doubling and link groups}
\label{subsec:intro-twisted}
Kulikov's band lemma conjugates one of the four inputs of a doubling by a
band.  Conjugating by an arbitrary braid is just as legitimate: the result
is again a positive factorization, and its group is computable
(Lemma~\ref{lem:twisted}).  This gives explicit curves whose degree does not
depend on the length of the relators.

\begin{theoremx}
\label{thm:twisted}
Let $\beta\in B_m$, let $s=d(\Delta_{m+1}^2)\in\cS_{2m+2}$, regard
$\beta$ as a braid on the first $m$ of $2m+2$ strings, and put
$U_\beta=d\bigl(\lambda(\beta)(s),s,s,s\bigr)\in\cS_{4m+4}$.
\begin{enumerate}[label=(\roman*),itemsep=1pt]
\item \label{b1}
$U_\beta$ is a factorization of $\Delta_{4m+4}^2$ into $16$ conjugates of
$\Delta_{m+1}^2$ and $12(m+1)$ conjugates of $\sigma_1$.
\item \label{b2}
$G(U_\beta)\cong\pi_1(S^3\setminus\hat\beta)\times\Z$, where $\hat\beta$
is the closed braid; the meridian $x_j$ goes to the meridian $x_i$ of the
$i$-th string of $\hat\beta$ if $j\equiv i\pmod{m+1}$ with $1\le i\le m$, and to
$(x_1\cdots x_m)^{-1}z$ if $j\equiv0\pmod{m+1}$.
\item \label{b3}
The corresponding Hurwitz curve $\overline H_\beta\subset\CP^2$ has
degree $4m+4$; its singular points are $16$ ordinary $(m+1)$-fold points; its
components are one component of degree $4\ell$ for each component of
$\hat\beta$ with $\ell$ strings, and one further smooth quartic; and
$\pi_1(\CP^2\setminus\overline H_\beta)\cong
\pi_1(S^3\setminus\hat\beta)\times\Z/4$.
\end{enumerate}
In particular, for the braid $\gamma_n\in B_3$ of \eqref{eq:gamma},
whose closure is $L_n$, $U_n\coloneqq U_{\gamma_n}$ is a factorization
of $\Delta_{16}^2$ into $64$ factors, and
\[
\pi_1(\C^2\setminus H_{\gamma_n})\cong\pi_1(S^3\setminus L_n)\times\Z,\qquad
\pi_1(\CP^2\setminus\overline H_{\gamma_n})\cong\pi_1(S^3\setminus L_n)\times\Z/4,
\]
where $\overline H_{\gamma_n}$ is the union of four smooth quartics meeting
pairwise transversally in the same $16$ points.  The curves
$\overline H_{\gamma_n}$ are not isotopic to algebraic curves.
\end{theoremx}

The same construction realizes, in degree $4m$, every C-group of the form
\[
\la x_1,\dots,x_m\mid l_m\text{ central},\ x_j=\beta_i(x_j)\ra
\]
with at most three braids $\beta_i\in B_m$ (Corollary~\ref{cor:braidtype}).  By contrast,
the degree produced by Kulikov's algorithm grows exponentially with the
lengths of the relators (\S\ref{subsec:768}).

In the other direction, if $L$ has at least two components and vanishing
linking numbers, the
components of any Hurwitz curve with group $\pi_1(S^3\setminus L)\times\Z$
meet only at points common to all of them, with local intersection numbers
proportional to the products of their degrees
(Theorem~\ref{thm:concurrency}).  Consequently, if $L$ is also non-split,
the degree of such a curve exceeds the number of components of $L$ by at
least two, and with the
C-structure of Theorem~\ref{thm:twisted} the curves $\overline H_{\gamma_n}$
have the smallest degree among Kulikov doublings (Proposition~\ref{prop:doubling-links}).

\subsection{Alexander invariants and non-algebraicity}
\label{subsec:intro-alex}
Abelian invariants of the groups $\widetilde G_k=\BS(k,k+1)\times\Z$ show
that the curves of Theorem~\ref{thm:explicit} are not isotopic to algebraic
curves.

\begin{theoremx}
\label{thm:alex}
Let $\widetilde G_k=\BS(k,k+1)\times\Z$ with the C-structure of
Theorem~\ref{thm:explicit}, $k\ge1$, and $q_k=3k^2+3k+1$.
\begin{enumerate}[label=(\roman*),itemsep=1pt]
\item \label{c1}
The Alexander invariant of $\widetilde G_k$ with respect to the
total-degree map is $A_0\cong\Z\oplus\Z/q_k$, with $t$ acting trivially on
$\Z$ and with order $3$ on $\Z/q_k$; in particular
$\Delta_{\widetilde G_k}(t)\doteq t-1$.
\item \label{c2}
$\V_1(\widetilde G_k)=\{(1,1)\}\cup\{(\zeta_k,\zeta_k^{-2})\}$,
$\zeta_k=k/(k+1)$, in the coordinates given by the two component meridians.
\item \label{c3}
Neither $\widetilde G_k$ nor $\BS(k,k+1)\times\Z/4$ is the fundamental
group of a smooth quasi-projective variety.  Hence the curves $\overline H_k$
of Theorem~\ref{thm:explicit} are not isotopic to algebraic curves in $\CP^2$.
\end{enumerate}
\end{theoremx}

The non-algebraicity in \ref{c3} is independent of Zariski's residual
finiteness problem: it holds for $k=1$, where $\BS(1,2)$ is residually
finite.  The obstruction is invisible to the one-variable Alexander
polynomial and to the Betti numbers of all finite cyclic branched covers,
but the torsion $\Z/q_k$ in \ref{c1} shows up integrally in the $n$-fold cyclic
covers with $3\mid n$ (Proposition~\ref{prop:covers}); it is $H_1$ of the
$3$-fold cyclic branched covering of $S^4$ along a knotted $2$-sphere with
group $\BS(k,k+1)$ (Remark~\ref{rem:branched}).  In general, the
characteristic variety of an affine Hurwitz curve group lies in the
hypersurface $\prod_c\lambda_c^{d_c}=1$ cut out by the degrees of the
components (Proposition~\ref{prop:V1-central}).

The groups $\pi_1(S^3\setminus L_n)\times\Z$ realized by
Theorem~\ref{thm:twisted} behave differently: they are not
quasi-projective, although their characteristic varieties are unions of
torsion-translated subtori, so that the torsion conditions of Arapura
\cite{Arapura-jag97} and of Artal Bartolo, Cogolludo and Matei
\cite{Artal-Cogolludo-Matei-gt13} are satisfied, and the obstruction comes
from the relative position of the components
(Proposition~\ref{prop:multivar-fails}).  For the group of the horizontal
arrangement $\A(31425)$ of \cite{Matei-Suciu-top00,Suciu-planes-r4}, also realized
by Theorem~\ref{thm:twisted} after a product with $\Z$, the obstruction comes
instead from the tangent cone of $\V_1$ at the identity, an irreducible cubic,
which is not a finite union of linear subspaces (Remark~\ref{rem:exotic24}).

\subsection{An algebraic model}
\label{subsec:intro-model}
The curves of Theorems~\ref{thm:explicit} and~\ref{thm:twisted} are one twist
away from explicit algebraic curves.  For $n\ge1$ put
\[
Q_n(x)=(x^2-n)(x^2-n-2),
\]
and for $n\ge2$ let $D_n\subset\CP^2$ be the closure of the affine curve
$\prod_{i=1}^n\bigl(Q_n(v)-(i+1)Q_1(u)\bigr)=0$, projected by $\pr(u,v)=u$.  The
quartics $Q_n(v)=\lambda Q_1(u)$ form a pencil, whose sixteen base points
$\{\pm1,\pm\sqrt3\}\times\{\pm\sqrt n,\pm\sqrt{n+2}\}$ are real.  Let $\Omega$
be the component of $\{u\in\C\mid\lvert u^2-2\rvert\le r_n\}$ in the
half-plane $\Re u>0$, where $r_n=\sqrt{(n+3)/(n+2)}$, and let
$u_b=\sqrt{2+r_n}\in\partial\Omega$.  Over $\Omega$ the curve $D_n$ splits into
its parts $D_n^\pm$ in the half-spaces $\pm\Re v>0$, and the $2n$ points of
$D_n^+$ over $u_b$ are real.  For $\gamma\in B_{2n}$, acting on these points
in their real order, let $D_n^\gamma$ be the Hurwitz curve obtained from
$D_n$ by cutting $D_n^+$ along $\pr^{-1}(\partial\Omega)$ and regluing it by a
fiber-preserving diffeomorphism that realizes $\gamma$; it is well defined
up to H-isotopy, and has the same singular points and tangencies as $D_n$
(\S\ref{subsec:pullback}).  Put $s_n=d(\Delta_n^2)$ and
$\omega=\Delta_{n,0}\Delta_{2n}^{-1}\in B_{2n}$.

\begin{figure}[ht]
\centering
\begin{tikzpicture}[x=0.82cm,y=0.82cm,
  q2/.style={line width=0.6pt},
  q3/.style={line width=0.6pt,dashed},
  q4/.style={line width=0.7pt,densely dotted}]
\fill[black!12] (0.951,0) rectangle (1.759,2.9);
\begin{scope}
\clip (-2.4,-2.9) rectangle (2.4,2.9);
\draw[q2] plot coordinates {(2.280,-2.900) (1.938,-2.505) (1.764,-2.284) (1.710,-2.200) (1.673,-2.128) (1.652,-2.063) (1.645,-2.001) (1.654,-1.925) (1.682,-1.841) (1.947,-1.294) (2.028,-1.088) (2.092,-0.888) (2.147,-0.671) (2.186,-0.450) (2.209,-0.228) (2.217,-0.004) (2.210,0.217) (2.188,0.435) (2.151,0.650) (2.099,0.866) (2.034,1.071) (1.953,1.280) (1.677,1.852) (1.649,1.950) (1.648,2.037) (1.661,2.095) (1.686,2.157) (1.783,2.309) (1.953,2.523) (2.280,2.900)};
\draw[q2] plot coordinates {(-2.280,2.900) (-1.938,2.505) (-1.764,2.284) (-1.710,2.200) (-1.673,2.128) (-1.652,2.063) (-1.645,2.001) (-1.654,1.925) (-1.682,1.841) (-1.947,1.294) (-2.028,1.088) (-2.092,0.888) (-2.147,0.671) (-2.186,0.450) (-2.209,0.228) (-2.217,0.004) (-2.210,-0.218) (-2.188,-0.435) (-2.151,-0.650) (-2.099,-0.866) (-2.034,-1.071) (-1.953,-1.280) (-1.677,-1.852) (-1.649,-1.950) (-1.648,-2.037) (-1.661,-2.095) (-1.686,-2.157) (-1.783,-2.309) (-1.953,-2.523) (-2.280,-2.900)};
\draw[q2] plot coordinates {(-0.042,-2.577) (-0.219,-2.564) (-0.390,-2.533) (-0.560,-2.483) (-0.718,-2.418) (-0.869,-2.335) (-0.995,-2.240) (-1.083,-2.146) (-1.112,-2.099) (-1.129,-2.055) (-1.136,-1.979) (-1.114,-1.900) (-1.061,-1.809) (-0.971,-1.700) (-0.843,-1.573) (-0.695,-1.450) (-0.549,-1.348) (-0.405,-1.266) (-0.303,-1.222) (-0.201,-1.190) (-0.102,-1.170) (-0.006,-1.164) (0.090,-1.169) (0.189,-1.187) (0.291,-1.218) (0.393,-1.261) (0.522,-1.331) (0.654,-1.419) (0.789,-1.526) (0.916,-1.642) (1.016,-1.751) (1.085,-1.845) (1.125,-1.929) (1.137,-2.001) (1.131,-2.048) (1.112,-2.099) (1.081,-2.150) (1.034,-2.204) (0.974,-2.258) (0.902,-2.313) (0.733,-2.411) (0.549,-2.487) (0.357,-2.540) (0.159,-2.571) (-0.042,-2.577)} -- cycle;
\draw[q2] plot coordinates {(-0.072,1.167) (-0.180,1.185) (-0.291,1.218) (-0.405,1.266) (-0.522,1.331) (-0.647,1.414) (-0.772,1.512) (-0.893,1.620) (-0.997,1.729) (-1.061,1.809) (-1.103,1.878) (-1.128,1.939) (-1.137,1.997) (-1.129,2.055) (-1.104,2.113) (-1.063,2.171) (-1.003,2.233) (-0.907,2.309) (-0.789,2.382) (-0.662,2.443) (-0.528,2.494) (-0.390,2.533) (-0.249,2.560) (-0.105,2.575) (0.039,2.577) (0.186,2.568) (0.330,2.546) (0.471,2.511) (0.609,2.465) (0.740,2.407) (0.863,2.338) (0.970,2.262) (1.051,2.186) (1.100,2.121) (1.129,2.055) (1.137,1.990) (1.124,1.925) (1.089,1.852) (1.031,1.769) (0.948,1.675) (0.843,1.573) (0.723,1.472) (0.600,1.381) (0.480,1.306) (0.363,1.247) (0.252,1.204) (0.141,1.177) (0.033,1.164) (-0.072,1.167)} -- cycle;
\draw[q3] plot coordinates {(2.162,-2.900) (1.922,-2.563) (1.778,-2.331) (1.734,-2.240) (1.703,-2.157) (1.684,-2.077) (1.678,-2.001) (1.686,-1.907) (1.711,-1.798) (1.920,-1.202) (1.977,-1.006) (2.022,-0.819) (2.060,-0.616) (2.087,-0.413) (2.104,-0.210) (2.109,-0.004) (2.104,0.199) (2.089,0.399) (2.064,0.595) (2.027,0.794) (1.980,0.995) (1.923,1.192) (1.707,1.812) (1.682,1.939) (1.681,2.048) (1.693,2.121) (1.715,2.193) (1.794,2.360) (1.934,2.581) (2.162,2.900)};
\draw[q3] plot coordinates {(-2.162,2.900) (-1.922,2.563) (-1.778,2.331) (-1.734,2.240) (-1.703,2.157) (-1.684,2.077) (-1.678,2.001) (-1.686,1.907) (-1.711,1.798) (-1.920,1.202) (-1.977,1.006) (-2.022,0.819) (-2.060,0.616) (-2.087,0.413) (-2.104,0.210) (-2.109,0.004) (-2.104,-0.199) (-2.089,-0.399) (-2.064,-0.595) (-2.027,-0.794) (-1.980,-0.995) (-1.923,-1.192) (-1.707,-1.812) (-1.682,-1.939) (-1.681,-2.048) (-1.693,-2.121) (-1.715,-2.193) (-1.794,-2.360) (-1.934,-2.581) (-2.162,-2.900)};
\draw[q3] plot coordinates {(-0.051,-2.675) (-0.240,-2.656) (-0.426,-2.610) (-0.603,-2.541) (-0.765,-2.450) (-0.908,-2.338) (-1.010,-2.222) (-1.046,-2.164) (-1.070,-2.106) (-1.084,-2.048) (-1.088,-1.986) (-1.078,-1.918) (-1.054,-1.841) (-1.017,-1.762) (-0.963,-1.671) (-0.802,-1.457) (-0.603,-1.242) (-0.432,-1.093) (-0.273,-0.990) (-0.198,-0.955) (-0.126,-0.932) (-0.054,-0.918) (0.018,-0.916) (0.090,-0.924) (0.162,-0.942) (0.237,-0.972) (0.315,-1.013) (0.477,-1.129) (0.654,-1.293) (0.844,-1.508) (0.988,-1.711) (1.035,-1.798) (1.066,-1.874) (1.084,-1.947) (1.088,-2.012) (1.080,-2.074) (1.061,-2.131) (1.029,-2.193) (0.986,-2.255) (0.930,-2.316) (0.862,-2.378) (0.705,-2.487) (0.528,-2.574) (0.342,-2.634) (0.147,-2.669) (-0.051,-2.675)} -- cycle;
\draw[q3] plot coordinates {(-0.042,0.917) (-0.117,0.929) (-0.192,0.953) (-0.270,0.988) (-0.348,1.034) (-0.519,1.165) (-0.703,1.345) (-0.883,1.559) (-0.955,1.660) (-1.011,1.751) (-1.051,1.834) (-1.076,1.907) (-1.087,1.976) (-1.085,2.041) (-1.070,2.106) (-1.044,2.168) (-1.002,2.233) (-0.948,2.298) (-0.879,2.363) (-0.801,2.425) (-0.718,2.479) (-0.625,2.530) (-0.438,2.607) (-0.237,2.656) (-0.036,2.676) (0.168,2.666) (0.348,2.633) (0.522,2.576) (0.687,2.498) (0.835,2.400) (0.954,2.291) (1.038,2.179) (1.064,2.124) (1.080,2.070) (1.088,2.015) (1.085,1.958) (1.069,1.885) (1.040,1.809) (0.994,1.722) (0.933,1.628) (0.759,1.407) (0.561,1.203) (0.396,1.067) (0.243,0.975) (0.171,0.945) (0.099,0.925) (0.027,0.916) (-0.042,0.917)} -- cycle;
\draw[q4] plot coordinates {(2.090,-2.900) (1.898,-2.588) (1.778,-2.356) (1.741,-2.262) (1.714,-2.171) (1.698,-2.088) (1.693,-2.001) (1.700,-1.892) (1.724,-1.765) (1.896,-1.150) (1.976,-0.779) (2.006,-0.584) (2.027,-0.392) (2.040,-0.197) (2.044,-0.004) (2.040,0.188) (2.028,0.377) (1.980,0.754) (1.899,1.138) (1.718,1.791) (1.696,1.932) (1.695,2.059) (1.706,2.135) (1.725,2.215) (1.795,2.392) (1.913,2.614) (2.090,2.900)};
\draw[q4] plot coordinates {(-2.090,2.900) (-1.898,2.588) (-1.778,2.356) (-1.741,2.262) (-1.714,2.171) (-1.698,2.088) (-1.693,2.001) (-1.700,1.892) (-1.724,1.765) (-1.896,1.150) (-1.976,0.779) (-2.006,0.584) (-2.027,0.391) (-2.040,0.197) (-2.044,0.004) (-2.040,-0.188) (-2.028,-0.377) (-1.980,-0.754) (-1.899,-1.138) (-1.718,-1.791) (-1.696,-1.932) (-1.695,-2.059) (-1.706,-2.135) (-1.725,-2.215) (-1.795,-2.393) (-1.913,-2.614) (-2.090,-2.900)};
\draw[q4] plot coordinates {(-0.081,-2.755) (-0.267,-2.729) (-0.447,-2.675) (-0.617,-2.596) (-0.769,-2.494) (-0.899,-2.374) (-0.953,-2.309) (-0.996,-2.244) (-1.028,-2.179) (-1.050,-2.113) (-1.062,-2.048) (-1.065,-1.983) (-1.055,-1.903) (-1.033,-1.820) (-0.997,-1.726) (-0.948,-1.624) (-0.794,-1.370) (-0.591,-1.091) (-0.411,-0.879) (-0.258,-0.737) (-0.189,-0.688) (-0.123,-0.654) (-0.060,-0.634) (0.000,-0.628) (0.063,-0.635) (0.129,-0.657) (0.198,-0.694) (0.270,-0.746) (0.348,-0.816) (0.435,-0.905) (0.624,-1.133) (0.825,-1.417) (0.968,-1.664) (1.014,-1.765) (1.044,-1.856) (1.061,-1.939) (1.064,-2.019) (1.055,-2.092) (1.035,-2.160) (1.004,-2.229) (0.961,-2.298) (0.905,-2.367) (0.838,-2.436) (0.760,-2.501) (0.678,-2.559) (0.590,-2.610) (0.501,-2.653) (0.312,-2.718) (0.117,-2.752) (-0.081,-2.755)} -- cycle;
\draw[q4] plot coordinates {(-0.039,0.631) (-0.102,0.646) (-0.168,0.676) (-0.237,0.721) (-0.312,0.782) (-0.480,0.956) (-0.672,1.197) (-0.853,1.461) (-0.979,1.686) (-1.019,1.780) (-1.046,1.863) (-1.061,1.943) (-1.064,2.019) (-1.056,2.088) (-1.037,2.157) (-1.006,2.226) (-0.963,2.295) (-0.909,2.363) (-0.846,2.429) (-0.774,2.490) (-0.693,2.549) (-0.522,2.644) (-0.336,2.712) (-0.141,2.750) (0.054,2.757) (0.243,2.734) (0.427,2.682) (0.603,2.603) (0.760,2.501) (0.896,2.378) (0.950,2.313) (0.993,2.248) (1.027,2.182) (1.050,2.113) (1.062,2.048) (1.064,1.979) (1.054,1.896) (1.029,1.805) (0.989,1.707) (0.932,1.595) (0.759,1.320) (0.543,1.031) (0.369,0.836) (0.291,0.764) (0.219,0.708) (0.153,0.668) (0.087,0.641) (0.024,0.629) (-0.039,0.631)} -- cycle;

\end{scope}
\draw[->,black!60] (-2.5,0) -- (2.6,0) node[right,black] {$u$};
\draw[->,black!60] (0,-3.0) -- (0,3.1) node[above,black] {$v$};
\foreach \a in {-1.732,-1,1,1.732}{\foreach \b in {-2.236,-1.732,1.732,2.236}{\fill (\a,\b) circle (1.5pt);}}
\node[below] at (0,-3.2) {(a)};
\end{tikzpicture}\hspace{2em}%
\begin{tikzpicture}[x=1.9cm,y=1.9cm]
\fill[black!12] plot coordinates {(1.759,0.000) (1.759,0.011) (1.759,0.022) (1.758,0.033) (1.757,0.043) (1.755,0.054) (1.754,0.065) (1.752,0.076) (1.749,0.086) (1.747,0.097) (1.744,0.107) (1.741,0.118) (1.737,0.128) (1.733,0.139) (1.729,0.149) (1.725,0.159) (1.720,0.169) (1.715,0.179) (1.709,0.188) (1.704,0.198) (1.698,0.207) (1.691,0.217) (1.685,0.226) (1.678,0.235) (1.671,0.244) (1.664,0.252) (1.656,0.261) (1.648,0.269) (1.640,0.277) (1.631,0.285) (1.623,0.292) (1.614,0.300) (1.604,0.307) (1.595,0.314) (1.585,0.320) (1.575,0.327) (1.565,0.333) (1.555,0.339) (1.544,0.344) (1.533,0.349) (1.522,0.354) (1.510,0.359) (1.499,0.363) (1.487,0.367) (1.475,0.371) (1.463,0.374) (1.451,0.377) (1.438,0.380) (1.425,0.382) (1.412,0.384) (1.399,0.385) (1.386,0.386) (1.373,0.387) (1.359,0.387) (1.346,0.387) (1.332,0.386) (1.318,0.385) (1.304,0.384) (1.290,0.381) (1.276,0.379) (1.262,0.376) (1.248,0.372) (1.234,0.368) (1.220,0.363) (1.206,0.358) (1.192,0.352) (1.177,0.346) (1.163,0.339) (1.150,0.331) (1.136,0.323) (1.122,0.314) (1.109,0.304) (1.096,0.294) (1.083,0.283) (1.070,0.271) (1.058,0.259) (1.046,0.246) (1.034,0.232) (1.023,0.218) (1.013,0.203) (1.003,0.187) (0.994,0.170) (0.985,0.153) (0.977,0.136) (0.971,0.117) (0.965,0.099) (0.960,0.079) (0.956,0.060) (0.953,0.040) (0.952,0.020) (0.951,0.000) (0.952,-0.020) (0.953,-0.040) (0.956,-0.060) (0.960,-0.079) (0.965,-0.099) (0.971,-0.117) (0.977,-0.136) (0.985,-0.153) (0.994,-0.170) (1.003,-0.187) (1.013,-0.203) (1.023,-0.218) (1.034,-0.232) (1.046,-0.246) (1.058,-0.259) (1.070,-0.271) (1.083,-0.283) (1.096,-0.294) (1.109,-0.304) (1.122,-0.314) (1.136,-0.323) (1.150,-0.331) (1.163,-0.339) (1.177,-0.346) (1.192,-0.352) (1.206,-0.358) (1.220,-0.363) (1.234,-0.368) (1.248,-0.372) (1.262,-0.376) (1.276,-0.379) (1.290,-0.381) (1.304,-0.384) (1.318,-0.385) (1.332,-0.386) (1.346,-0.387) (1.359,-0.387) (1.373,-0.387) (1.386,-0.386) (1.399,-0.385) (1.412,-0.384) (1.425,-0.382) (1.438,-0.380) (1.451,-0.377) (1.463,-0.374) (1.475,-0.371) (1.487,-0.367) (1.499,-0.363) (1.510,-0.359) (1.522,-0.354) (1.533,-0.349) (1.544,-0.344) (1.555,-0.339) (1.565,-0.333) (1.575,-0.327) (1.585,-0.320) (1.595,-0.314) (1.604,-0.307) (1.614,-0.300) (1.623,-0.292) (1.631,-0.285) (1.640,-0.277) (1.648,-0.269) (1.656,-0.261) (1.664,-0.252) (1.671,-0.244) (1.678,-0.235) (1.685,-0.226) (1.691,-0.217) (1.698,-0.207) (1.704,-0.198) (1.709,-0.188) (1.715,-0.179) (1.720,-0.169) (1.725,-0.159) (1.729,-0.149) (1.733,-0.139) (1.737,-0.128) (1.741,-0.118) (1.744,-0.107) (1.747,-0.097) (1.749,-0.086) (1.752,-0.076) (1.754,-0.065) (1.755,-0.054) (1.757,-0.043) (1.758,-0.033) (1.759,-0.022) (1.759,-0.011)} -- cycle;
\draw[line width=0.5pt] plot coordinates {(1.759,0.000) (1.759,0.011) (1.759,0.022) (1.758,0.033) (1.757,0.043) (1.755,0.054) (1.754,0.065) (1.752,0.076) (1.749,0.086) (1.747,0.097) (1.744,0.107) (1.741,0.118) (1.737,0.128) (1.733,0.139) (1.729,0.149) (1.725,0.159) (1.720,0.169) (1.715,0.179) (1.709,0.188) (1.704,0.198) (1.698,0.207) (1.691,0.217) (1.685,0.226) (1.678,0.235) (1.671,0.244) (1.664,0.252) (1.656,0.261) (1.648,0.269) (1.640,0.277) (1.631,0.285) (1.623,0.292) (1.614,0.300) (1.604,0.307) (1.595,0.314) (1.585,0.320) (1.575,0.327) (1.565,0.333) (1.555,0.339) (1.544,0.344) (1.533,0.349) (1.522,0.354) (1.510,0.359) (1.499,0.363) (1.487,0.367) (1.475,0.371) (1.463,0.374) (1.451,0.377) (1.438,0.380) (1.425,0.382) (1.412,0.384) (1.399,0.385) (1.386,0.386) (1.373,0.387) (1.359,0.387) (1.346,0.387) (1.332,0.386) (1.318,0.385) (1.304,0.384) (1.290,0.381) (1.276,0.379) (1.262,0.376) (1.248,0.372) (1.234,0.368) (1.220,0.363) (1.206,0.358) (1.192,0.352) (1.177,0.346) (1.163,0.339) (1.150,0.331) (1.136,0.323) (1.122,0.314) (1.109,0.304) (1.096,0.294) (1.083,0.283) (1.070,0.271) (1.058,0.259) (1.046,0.246) (1.034,0.232) (1.023,0.218) (1.013,0.203) (1.003,0.187) (0.994,0.170) (0.985,0.153) (0.977,0.136) (0.971,0.117) (0.965,0.099) (0.960,0.079) (0.956,0.060) (0.953,0.040) (0.952,0.020) (0.951,0.000) (0.952,-0.020) (0.953,-0.040) (0.956,-0.060) (0.960,-0.079) (0.965,-0.099) (0.971,-0.117) (0.977,-0.136) (0.985,-0.153) (0.994,-0.170) (1.003,-0.187) (1.013,-0.203) (1.023,-0.218) (1.034,-0.232) (1.046,-0.246) (1.058,-0.259) (1.070,-0.271) (1.083,-0.283) (1.096,-0.294) (1.109,-0.304) (1.122,-0.314) (1.136,-0.323) (1.150,-0.331) (1.163,-0.339) (1.177,-0.346) (1.192,-0.352) (1.206,-0.358) (1.220,-0.363) (1.234,-0.368) (1.248,-0.372) (1.262,-0.376) (1.276,-0.379) (1.290,-0.381) (1.304,-0.384) (1.318,-0.385) (1.332,-0.386) (1.346,-0.387) (1.359,-0.387) (1.373,-0.387) (1.386,-0.386) (1.399,-0.385) (1.412,-0.384) (1.425,-0.382) (1.438,-0.380) (1.451,-0.377) (1.463,-0.374) (1.475,-0.371) (1.487,-0.367) (1.499,-0.363) (1.510,-0.359) (1.522,-0.354) (1.533,-0.349) (1.544,-0.344) (1.555,-0.339) (1.565,-0.333) (1.575,-0.327) (1.585,-0.320) (1.595,-0.314) (1.604,-0.307) (1.614,-0.300) (1.623,-0.292) (1.631,-0.285) (1.640,-0.277) (1.648,-0.269) (1.656,-0.261) (1.664,-0.252) (1.671,-0.244) (1.678,-0.235) (1.685,-0.226) (1.691,-0.217) (1.698,-0.207) (1.704,-0.198) (1.709,-0.188) (1.715,-0.179) (1.720,-0.169) (1.725,-0.159) (1.729,-0.149) (1.733,-0.139) (1.737,-0.128) (1.741,-0.118) (1.744,-0.107) (1.747,-0.097) (1.749,-0.086) (1.752,-0.076) (1.754,-0.065) (1.755,-0.054) (1.757,-0.043) (1.758,-0.033) (1.759,-0.022) (1.759,-0.011)} -- cycle;
\draw[->,black!60] (-0.3,0) -- (2.5,0) node[right,black] {$\Re u$};
\draw[->,black!60] (0,-1.15) -- (0,1.2) node[above,black] {$\Im u$};

\fill (1.000,0) circle (1.6pt);
\fill (1.732,0) circle (1.6pt);
\draw[line width=0.4pt] (1.065,0) circle (1.1pt);
\draw[line width=0.4pt] (1.088,0) circle (1.1pt);
\draw[line width=0.4pt] (1.137,0) circle (1.1pt);
\draw[line width=0.4pt] (1.645,0) circle (1.1pt);
\draw[line width=0.4pt] (1.678,0) circle (1.1pt);
\draw[line width=0.4pt] (1.693,0) circle (1.1pt);
\draw[line width=0.5pt] (2.044,-0.04) -- ++(0,0.08);
\draw[line width=0.5pt] (2.109,-0.04) -- ++(0,0.08);
\draw[line width=0.5pt] (2.217,-0.04) -- ++(0,0.08);
\draw[line width=0.5pt] (-0.04,0.424) -- ++(0.08,0);
\draw[line width=0.5pt] (-0.04,-0.424) -- ++(0.08,0);
\draw[line width=0.5pt] (-0.04,0.670) -- ++(0.08,0);
\draw[line width=0.5pt] (-0.04,-0.670) -- ++(0.08,0);
\draw[line width=0.5pt] (-0.04,0.957) -- ++(0.08,0);
\draw[line width=0.5pt] (-0.04,-0.957) -- ++(0.08,0);
\draw[line width=0.5pt,fill=white] (1.759,0) ++(-1.5pt,-1.5pt) rectangle ++(3pt,3pt);
\node[above right] at (1.759,0.02) {$u_b$};
\node at (1.36,0.17) {$\Omega$};
\node[below] at (1.1,-1.3) {(b)};
\end{tikzpicture}
\caption{(a) The real points of $D_3$: the quartics
$(v^2-3)(v^2-5)=\lambda(u^2-1)(u^2-3)$ for $\lambda=2$ (solid), $3$ (dashed),
and $4$ (dotted), and their sixteen base points.  The shaded strip is the
real part of $\pr^{-1}(\Omega)\cap\{\Re v>0\}$, which contains the real part
of $D_3^+$.  (b) The half $\Re u\ge0$ of the $u$-line, with the region
$\Omega$ and the base point $u_b$; the picture is symmetric under
$u\mapsto-u$.  Dots mark the triple-point fibers $u=1,\sqrt3$; circles, the
two clusters of three critical values where $D_3$ is tangent to the
pencil at $v=\pm2$; ticks, the critical values where it is tangent at
$v=0$, three on the real axis and three pairs on the imaginary axis.}
\label{fig:D3}
\end{figure}
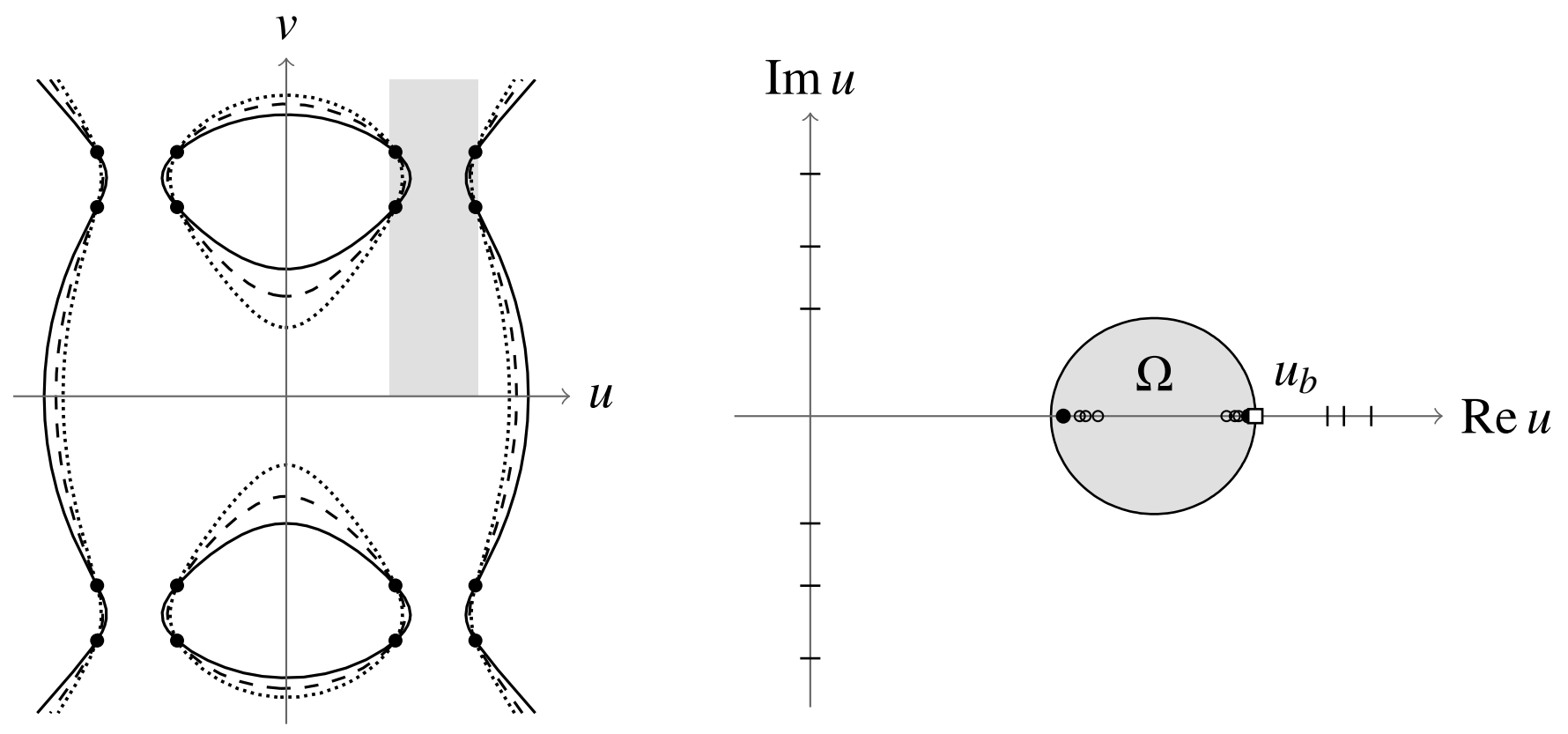

\begin{theoremx}
\label{thm:model}
\begin{enumerate}[label=(\roman*),itemsep=1pt]
\item \label{mod1}
$D_n$ is the union of $n$ smooth quartics of the pencil, which meet pairwise
transversally at its base points; these are ordinary $n$-fold points, and the
braid monodromy factorization of $D_n$ has the type of $d(s_n)$.
\item \label{mod2}
For $\gamma\in B_{2n}$, the braid monodromy factorization of $D_n^\gamma$
has the type of
\[
d\bigl(\lambda(\omega^{-1}\gamma\omega)s_n,s_n,s_n,s_n\bigr).
\]
\item \label{mod3}
Consequently, the curve $\overline H_k$ of Theorem~\ref{thm:explicit} and
$D_3^{\,\omega b_k\omega^{-1}}$ become H-isotopic once the critical points
lying in common fibers are moved apart by small deformations
(\S\ref{subsec:hurwitz}); in particular, they are isotopic through Hurwitz
curves with the same singular points and tangencies.  The same holds for the
curve $\overline H_\beta$ of Theorem~\ref{thm:twisted} and
$D_{m+1}^{\,\omega\beta\omega^{-1}}$.
\end{enumerate}
\end{theoremx}

Thus, up to such an isotopy, $\overline H_k$ is obtained from the union $D_3$ 
of the three quartics $(v^2-3)(v^2-5)=\lambda(u^2-1)(u^2-3)$, $\lambda=2,3,4$, 
whose group is $F_2\times\Z$, by one twist over $\Omega$; the twist preserves the 
sixteen triple points and the $36$ tangencies, merges two of the quartics, and turns
$F_2\times\Z$ into $\BS(k,k+1)\times\Z$.  This twist is a braiding of $D_3$
along an annulus over $\partial\Omega$, in the sense of Auroux, Donaldson and
Katzarkov \cite[\S3]{Auroux-Donaldson-Katzarkov-mathann03}.  In their
reinterpretation of Moishezon's infinite families of pairwise non-isotopic
cuspidal curves \cite{Moishezon-conm94}, the preimage of a suitable algebraic
curve under a generic cubic map $\CP^2\to\CP^2$ is, topologically, a union of
cubics of a pencil modified near the base points, and the curves of the family
arise by braiding two of these cubics repeatedly along an annulus
\cite[\S4]{Auroux-Donaldson-Katzarkov-mathann03}.  The proof rests on a geometric
reading of Kulikov's formula: $d(s)$ is the braid monodromy of the preimage,
under $(u,v)\mapsto(u^2-\alpha_0,v^2-\beta_0)$, of a curve with braid monodromy
$s$, its four inputs corresponding to the four quarters cut out by the signs
of $\Re u$ and $\Re v$ (Lemma~\ref{lem:geomdouble} and
Figure~\ref{fig:pullback}).  Applied to $n$
concurrent lines, it shows that $s_n=d(\Delta_n^2)$ is the braid monodromy of
$n$ conics of a pencil with four base points, and applied once more, that
$d(s_n)$ is that of $D_n$.  Figure~\ref{fig:D3} shows the real points of
$D_3$ and the region $\Omega$; the twisting bands for $\overline H_1$ and
$\overline H_2$ are written down in Remark~\ref{rem:twistbands}.

\subsection{Local versus global}
\label{subsec:intro-local}
The motivating local picture is that of \cite{Suciu-borro}: the links
$L_n$ are links of isolated singularities of semiholomorphic maps
$f\colon(\R^4,0)\to(\R^2,0)$, holomorphic in one variable.  Their zero sets
are braided surfaces, in a slightly generalized sense of Rudolph
\cite{Rudolph-cmh83}, but they
are \emph{not} Hurwitz curve germs: the boundary braid has exponent sum
zero, whereas every non-trivial Hurwitz germ has boundary braid of positive
exponent sum (\S\ref{sec:local}).  The horizontal
$2$-arrangements $\A(\sigma)$ of \cite{Matei-Suciu-top00,Suciu-planes-r4},
a separate family of braided-surface germs, are Hurwitz germs only for
$\sigma=\id$ (\S\ref{subsec:arrangements}).  On the
group-theoretic side, every C-group is a local Hurwitz group
\cite[Thm.~7.2]{Kulikov-usp07}, so there is no local obstruction at all;
the global obstruction is exactly Proposition~\ref{prop:intro-central}.

\subsection*{Organization}
Section~\ref{sec:prelim} recalls Kulikov's conventions for braids, Hurwitz
curves, C-groups and doublings.  Section~\ref{sec:central} proves
Proposition~\ref{prop:intro-central} and collects its consequences, and
Section~\ref{sec:kulikova} reviews
Kulikova's example and the degree produced by Kulikov's algorithm, and shows
that curves with isomorphic groups can be told apart by their meridians,
although not by the fundamental groups of their $2$-fold branched coverings.
Theorems~\ref{thm:explicit}, \ref{thm:twisted}, \ref{thm:alex}
and~\ref{thm:model} are proved in Sections~\ref{sec:explicit},
\ref{sec:twisted}, \ref{sec:alexander} and~\ref{sec:model}.
Section~\ref{sec:mindeg} collects what we know
about the minimal degree of a Hurwitz curve with group $\BS(k,k+1)\times\Z$
or $\pi_1(S^3\setminus L)\times\Z$,
Section~\ref{sec:local} treats the local side, and
Section~\ref{sec:questions} lists open problems.

\section{Braids, Hurwitz curves, and C-groups}
\label{sec:prelim}

We follow the conventions of Kulikov \cite{Kulikov-izv04a} throughout
(in the survey \cite{Kulikov-usp07} braids act on the right), since the 
explicit factorizations of \S\ref{sec:explicit} are meant
to be checked against his formulas.  We recall the action of braids on the
free group (\S\ref{subsec:braids}), factorizations and their groups $G(s)$
(\S\ref{subsec:factorizations}), Hurwitz curves and the Zariski--van Kampen
theorem (\S\ref{subsec:hurwitz}), C-groups (\S\ref{subsec:cgroups}), and
Kulikov's doubling formula, our basic tool for producing new factorizations
of the full twist (\S\ref{subsec:doubling}).

\subsection{Braids acting on the free group}
\label{subsec:braids}
Let $F_m=\la x_1,\dots,x_m\ra$ be the fundamental group of a disk with $m$
punctures on a horizontal diameter, freely generated by the standard loops
$x_1,\dots,x_m$ around the punctures, from left to right.  The braid group
$B_m$, with standard generators $\sigma_1,\dots,\sigma_{m-1}$, acts on the
left on $F_m$ by
\begin{equation}
\label{eq:action}
\sigma_i(x_i)=x_ix_{i+1}x_i^{-1},\qquad \sigma_i(x_{i+1})=x_i,\qquad \sigma_i(x_j)=x_j\
(j\neq i,i+1),
\end{equation}
with $(\beta\gamma)(x)=\beta(\gamma(x))$.  This action is faithful.  We write
$l_m=x_1\cdots x_m$; it is fixed by every braid.  The full twist
$\Delta_m^2=(\sigma_1\cdots \sigma_{m-1})^m$ generates the center
of $B_m$ and acts by
\begin{equation}
\label{eq:fulltwist}
\Delta_m^2(x_i)=l_m\,x_i\,l_m^{-1},\qquad i=1,\dots,m .
\end{equation}
For $k+i\le m$ let $\Delta_{k,i}$ be the Garside element of the subgroup
generated by $\sigma_{i+1},\dots,\sigma_{i+k-1}$, and write $\Delta_k=\Delta_{k,0}$.

An element $b\in B_m$ conjugate to $\sigma_1$, $b=g\,\sigma_1\,g^{-1}$, is a
\emph{band}; its \emph{endpoint loops} are $y_1=g(x_1)$ and $y_2=g(x_2)$.

\begin{lemma}[{\cite[Lem.~1.10]{Kulikov-izv04a}}]
\label{lem:bandrel}
Let $b=g\,\sigma_1g^{-1}$ be a band with endpoint loops $y_1,y_2$.  Then the
normal closure of $\{x^{-1}b(x)\mid x\in F_m\}$ is $\nc{y_1^{-1}y_2}$.  In
particular $b(w)\in w\,\nc{y_1^{-1}y_2}$ for all $w\in F_m$.
\end{lemma}

\begin{proof}
The elements $g(x_1),\dots,g(x_m)$ form a free basis, and $b$ acts on it as
$\sigma_1$ acts on $x_1,\dots,x_m$: $b(y_1)=y_1y_2y_1^{-1}$, $b(y_2)=y_1$, 
and $b$ fixes $g(x_j)$ for $j\ge3$.  Since $(uv)^{-1}b(uv)=v^{-1}(u^{-1}b(u))v\cdot
v^{-1}b(v)$, the normal closure in question is the normal
closure of the $y^{-1}b(y)$ for $y$ in this basis, i.e., of $y_2y_1^{-1}$ and
$y_2^{-1}y_1$.
\end{proof}

\subsection{Factorization semigroups}
\label{subsec:factorizations}
A \emph{factorization} in $B_m$ is a finite sequence
$s=u_1\cdot u_2\cdots u_n$ of elements of $B_m$, considered up to Hurwitz
moves $u_i\cdot u_{i+1}\mapsto u_{i+1}\cdot u_{i+1}^{-1}u_iu_{i+1}$ (and their
inverses); these form a semigroup $\cS_m$ under juxtaposition, with product
homomorphism $\alpha(s)=u_1u_2\cdots u_n\in B_m$.  Hurwitz moves, and
Hurwitz equivalence, make sense for finite sequences of elements of any
group, and preserve their product.  All factors below will be
conjugates of $\sigma_1$ or of full twists $\Delta_k^2$.  For $\beta\in B_m$ put
$\lambda(\beta)(s)=\beta u_1\beta^{-1}\cdots\beta u_n\beta^{-1}$, so that
$\alpha(\lambda(\beta)s)=\beta\,\alpha(s)\beta^{-1}$.  The shift
$\sh=\sh_m\colon B_m\to B_{2m}$, $\sigma_j\mapsto \sigma_{m+j}$, induces
$\sh\colon\cS_m\to\cS_{2m}$, and the inclusion $B_m\subset B_{2m}$ on the
first $m$ strings induces $\cS_m\subset\cS_{2m}$.

The \emph{group} of a factorization is the C-group
\begin{equation}
\label{eq:Gs}
G(s)=\la x_1,\dots,x_m\mid x_j=u_i(x_j),\ 1\le i\le n,\ 1\le j\le m\ra
=F_m/N(s),
\end{equation}
where $N(s)$ is the normal closure of the relators $x_j^{-1}u_i(x_j)$; it
depends only on the subgroup generated by the factors \cite[(22)]{Kulikov-izv04a}.
Two facts are used repeatedly: $N(s\cdot s')=N(s)N(s')$, and
\begin{equation}
\label{eq:Nconj}
N(\lambda(\beta)s)=\beta\bigl(N(s)\bigr),
\end{equation}
since $x^{-1}\beta u\beta^{-1}(x)=\beta\bigl(w^{-1}u(w)\bigr)$ with
$w=\beta^{-1}(x)$.  In particular $G(\lambda(\beta)s)\cong G(s)$.

\subsection{Hurwitz curves}
\label{subsec:hurwitz}
Let $\pr\colon\CP^2\setminus\{p\}\to\CP^1$ be the linear projection
from a point $p$.  Following \cite[Def.~5.1]{Kulikov-izv04a} and
\cite[\S0]{Kulikova-izv05}, a \emph{Hurwitz curve} of degree $m$ is the image
$\overline H=f(S)$ of a smooth map $f\colon S\to\CP^2\setminus\{p\}$ of a closed
oriented surface, such that for some finite set $Z\subset\overline H$:
$f$ is an embedding off $f^{-1}(Z)$, and $\overline H$ meets the lines of the
pencil transversally and positively off $Z$; near each point of $Z$,
$\overline H$ is a complex-analytic curve with its complex orientation; and
$\pr|_{\overline H}$ has degree $m$.  The \emph{irreducible
components} of $\overline H$ are the images of the components of $S$.

Fix a line $\Linf$ of the pencil in general position with respect to
$\overline H$ and put $H=\overline H\setminus\Linf\subset\C^2$.  We identify
$B_m$ with the mapping class group of $(\C,K)$, for a set $K$ of $m$ real
numbers, that is, the group of isotopy classes of diffeomorphisms of $\C$
preserving $K$ and supported in a large disk, which is the mapping class
group of that disk with $m$ marked points, relative to its boundary
\cite[\S9.1.3]{Farb-Margalit-12}.  We do so in such a way that the product of
braids is the composite of diffeomorphisms and braids act on
$\pi_1(\C\setminus K,*)\cong F_m$, with $*$ of large modulus and
$x_1,\dots,x_m$ the standard loops of \S\ref{subsec:braids} for the real
order of $K$, by~\eqref{eq:action} (cf.\ \cite[Thm.~1.10]{Birman-74});
this is the identification implicit in Theorem~\ref{thm:zvk} below.  Lifting
a geometric basis of loops around the critical values of $\pr|_{\overline H}$,
based at a point of the boundary of a disk that contains them and with product
the boundary loop of that disk, defines the \emph{braid monodromy
factorization} $\operatorname{bmf}(\overline H)\in\cS_m$.  For algebraic curves this is
Moishezon's braid monodromy; we refer to \cite{Cohen-Suciu-cmh97} for an
account based on Hansen's polynomial covers, which applies verbatim to
Hurwitz curves, since off the critical fibers $H$ is the zero set of a
simple Weierstrass polynomial with continuous coefficients
\cite[\S1]{Kulikov-usp07}.  Another geometric basis changes
$\operatorname{bmf}(\overline H)$ by Hurwitz moves \cite[\S4.3]{Kulikov-usp07},
and so does any sequence of loops that is Hurwitz equivalent to a geometric
basis; we use such sequences in place of geometric bases.

The factor along a loop around a critical value $u_1$ is the braid
monodromy of the germ of $\overline H$ over $u_1$, transported to the base
fiber.  Let $p_1,\dots,p_k$ be the critical points in the fiber over $u_1$,
and choose a closed disk $W\ni u_1$ containing no other critical value and
disjoint closed disks $E_1,\dots,E_k$ in the $v$-line, with $p_j\in W\times
E_j$, such that $\overline H$ does not meet $W\times\partial E_j$ and is
complex analytic on $\Sigma_j=\overline H\cap(W\times E_j)$.  Then the factor
is the product of $k$ braids with disjoint supports, the motions of the points
of the $\Sigma_j$ over $\partial W$, the other points being held fixed,
transported to the base fiber \cite[\S4.2]{Kulikov-usp07}; we call them its
\emph{local factors}.  They commute, and we always replace the factor by the
sequence of its local factors, in any order.  This changes neither the
product nor the group of the factorization: in a geometric basis of $F_m$
whose loops around the points of each support avoid the other supports, the
factor and its $j$-th local factor agree on the loops around the points of
the $j$-th support, and the local factors fix all other loops, so the
argument in the proof of Lemma~\ref{lem:bandrel} applies.  Thus the factors
of $\operatorname{bmf}(\overline H)$ correspond to the critical points of
$\pr|_{\overline H}$, and are conjugate to the braid monodromies of the
analytic germs of $\overline H$ at these points
--- a simple tangency with a line of the pencil ($v^2=u$) contributes a
conjugate of $\sigma_1$, an ordinary triple point ($v^3=u^3$, equivalently
$v(v^2-u^2)=0$) a conjugate of $\Delta_3^2$ --- and
\begin{equation}
\label{eq:alpha}
\alpha\bigl(\operatorname{bmf}(\overline H)\bigr)=\Delta_m^2
\end{equation}
\cite[Lem.~5.2]{Kulikov-izv04a}.  Conversely, every factorization of
$\Delta_m^2$ whose factors are conjugates of braid monodromies of analytic
germs is the braid monodromy factorization of a Hurwitz curve
\cite[Thm.~5.3]{Kulikov-izv04a}, \cite{Moishezon-conm94};
Figure~\ref{fig:double} draws such a factorization.  The \emph{type} of a
factorization $s$ is its orbit $\{\lambda(\beta)s\}$ under simultaneous
conjugation, which corresponds to braid-equivalence of braid monodromy
homomorphisms in the sense of \cite[Def.~3.6]{Cohen-Suciu-cmh97}.

The natural equivalence relation for Hurwitz curves is Kulikov's
\emph{H-isotopy}, the H standing for Hurwitz \cite[\S5.2]{Kulikov-usp07}:
two Hurwitz curves are H-isotopic if one is carried to the other by an
isotopy of $\CP^2$ that preserves the center of the pencil and passes through
Hurwitz curves, two critical points lying in a common fiber at one time if
and only if they do so at all times.  It plays the role that equisingular
deformation plays for algebraic curves, and it preserves the singularities,
the tangencies and the groups of the complements.  The type of
$\operatorname{bmf}(\overline H)$ does not depend on the choices made in
defining it, and it is an invariant of H-isotopy.  A much coarser relation is
regular homotopy, which allows in addition the creation and cancellation of
pairs of nodes of opposite signs: any two irreducible cuspidal Hurwitz curves
of the same degree, with the same numbers of cusps and nodes, are regular
homotopic.  This is \cite[Thm.~0.1]{Auroux-Kulikov-Shevchishin-izv04}, applied
in the Hirzebruch surface $\mathbb{F}_1\cong\CP^2\#\overline{\CP}{}^2$, the
blow-up of $\CP^2$ at the center of the pencil, in which Hurwitz curves of 
degree $m$ avoid the exceptional curve and all represent the same homology 
class (compare \cite[\S4.2]{Auroux-Kulikov-Shevchishin-izv04}).

An H-isotopy cannot move critical points out of a common fiber, but a small
deformation can.  With $u_1$, $W$, $E_j$ and $\Sigma_j$ as above and $k\ge2$,
choose diffeomorphisms $\varphi_j$ of $W$, equal to the identity near
$\partial W$ and to translations $u\mapsto u+\varepsilon_j$ near $u_1$, with
$\varepsilon_1,\dots,\varepsilon_k$ small and distinct, and replace each
$\Sigma_j$ by its image under $(u,v)\mapsto(\varphi_j(u),v)$.  This map
carries lines of the pencil to lines of the pencil, preserving orientations,
and is holomorphic near $p_j$.  Doing this over every critical value whose
fiber contains several critical points, we obtain a Hurwitz curve
$\overline H'$ with the same singular points and tangencies, each fiber of
$\pr$ now containing at most one of them; we call $\overline H'$ a
\emph{separation} of $\overline H$.  Shrinking the $\varepsilon_j$ to $0$
gives an ambient isotopy from $\overline H'$ to $\overline H$ through
Hurwitz curves, which is an H-isotopy away from its endpoint.  Over $W$, the
curve $\overline H'$ is the disjoint union of the images of the $\Sigma_j$
and of some sections, so in a geometric basis for $\overline H'$ the loop
around $u_1$ is replaced by loops around the points $u_1+\varepsilon_j$, along
which the monodromies are the local factors.  Hence $\overline H$ and
$\overline H'$ have the same braid monodromy factorization.

Conversely, two Hurwitz curves whose critical points are of inseparable type
and lie in distinct fibers are H-isotopic as soon as their braid monodromy
factorizations have the same type \cite[Thm.~5.2]{Kulikov-usp07}.  Ordinary
$k$-fold points and simple tangencies are of inseparable type, since their
braid monodromies $\Delta_k^2$ and $\sigma_1$ are the standard forms of the
singularities $v^k=u^k$ and $v^2=u$, which are inseparable by
\cite[Lem.~5.5]{Kulikov-usp07}.  Hence, if the critical points of two Hurwitz
curves are ordinary multiple points and simple tangencies, and their braid
monodromy factorizations have the same type, then any separations of them
are H-isotopic; in particular, the two curves are isotopic through Hurwitz
curves with the same singular points and tangencies.

\begin{theorem}[Zariski--van Kampen {\cite[\S6]{Kulikov-izv04a}, \cite[\S7]{Kulikov-usp07}}]
\label{thm:zvk}
If $s=\operatorname{bmf}(\overline H)$, then $\pi_1(\C^2\setminus H)\cong G(s)$,
the generators $x_j$ being meridians of $H$ in a generic line of the pencil,
and $\pi_1(\CP^2\setminus\overline H)\cong G(s)/\nc{l_m}$.
\end{theorem}

The same holds for pseudoholomorphic curves with respect to a pencil of
pseudoholomorphic lines, and for Hurwitz curves with negative nodes
\cite[\S0]{Kulikova-izv05}, \cite[\S4.2, Prop.~4.3, \S7]{Kulikov-usp07};
Proposition~\ref{prop:intro-central}
and its consequences apply to all of these, since they use only
Theorem~\ref{thm:zvk} and \eqref{eq:alpha}.

\subsection{C-groups}
\label{subsec:cgroups}
A \emph{C-group} \cite{Kulikov-izv95} is a group with a finite presentation
$\la x_1,\dots,x_m\mid x_i=w_{i,j,k}^{-1}x_jw_{i,j,k}\ra$ in which each relation
conjugates a generator to a generator; the conjugates of the $x_i$ are the
\emph{C-generators}.  Every C-group is isomorphic to $G(\beta_1\cdots\beta_n)$
for suitable bands $\beta_i$ \cite[Thm.~1.14]{Kulikov-izv04a}.  Following
\cite{Kulikova-izv05}, a C-group $G$ is a \emph{Hurwitz C-group of degree $m$} if
it has C-generators $x_1,\dots,x_m$ that generate $G$ and whose product
$l_m$ is central.  By \eqref{eq:alpha} and \eqref{eq:fulltwist}, every group
$G(s)$ with $\alpha(s)=\Delta_m^2$ is a Hurwitz C-group of degree $m$;
conversely every Hurwitz C-group of degree $m$ is C-isomorphic to
$\pi_1(\C^2\setminus H)$ for a Hurwitz curve whose degree is $m$ times a
power of $2$, all of whose singular points are of type $v^m=u^m$
\cite[Thms.~3.6, 6.1, 6.2]{Kulikov-izv04a}, \cite[Thm.~7.3]{Kulikov-usp07}.  Let
$\nu\colon G\to\Z$ send every C-generator to $1$ and put $N=\ker\nu$.  The
\emph{Alexander invariant} of $G$ with respect to $\nu$ is $A_0(G)=N/N'$, the
first homology of the infinite cyclic cover defined by $\nu$, a module over
$\Lambda=\Z[t^{\pm1}]$ via conjugation by a C-generator; Kulikov calls it the
Alexander module of $G$.  If $G$ is irreducible, then $N=G'$ and $A_0(G)$ is
the Alexander invariant $G'/G''$.  We write $\Delta_G(t)$ for the
characteristic polynomial of $t$ on $A_0(G)\otimes\C$.

\subsection{Kulikov's doubling}
\label{subsec:doubling}
\begin{figure}[ht]
\centering
\begin{tikzpicture}[line width=0.9pt, line cap=round,
    halo/.style={draw=white, line width=3.6pt, line cap=butt}, gr/.style={black!35}]
\draw[line width=0.3pt,black!50,densely dotted]
  (5.070,-0.25) -- ++(0,2.60) (6.540,-0.25) -- ++(0,2.60) (8.910,-0.25) -- 
  ++(0,2.60) (9.480,-0.25) -- ++(0,2.60) (10.950,-0.25) -- ++(0,2.60);
\draw[gr]
  (4.500,0.00) -- (6.540,0.00) (6.990,0.00) -- (8.340,0.00) (8.910,0.00) -- (10.950,0.00)
  (11.400,0.00) -- (12.750,0.00) (4.500,0.42) -- (5.070,0.42) (5.520,0.42) -- (5.970,0.42)
  (7.440,0.42) -- (7.890,0.42) (8.790,0.42) -- (9.480,0.42) (9.930,0.42) -- (10.380,0.42)
  (11.850,0.42) -- (12.300,0.42) (13.200,0.42) -- (13.320,0.42) (6.420,0.84) -- (6.990,0.84)
  (8.340,0.84) -- (8.910,0.84) (10.830,0.84) -- (11.400,0.84) (12.750,0.84) -- (13.320,0.84)
  (6.420,1.26) -- (6.990,1.26) (8.340,1.26) -- (8.910,1.26) (10.830,1.26) -- (11.400,1.26)
  (12.750,1.26) -- (13.320,1.26) (4.500,1.68) -- (5.070,1.68) (5.520,1.68) -- (5.970,1.68)
  (7.440,1.68) -- (7.890,1.68) (8.790,1.68) -- (9.480,1.68) (9.930,1.68) -- (10.380,1.68)
  (11.850,1.68) -- (12.300,1.68) (13.200,1.68) -- (13.320,1.68) (4.500,2.10) -- (6.540,2.10)
  (6.990,2.10) -- (8.340,2.10) (8.910,2.10) -- (10.950,2.10) (11.400,2.10) -- (12.750,2.10);
\draw
  (0.000,0.00) -- (4.500,0.00) (8.790,0.00) -- (8.910,0.00) (13.200,0.00) -- (14.370,0.00)
  (0.000,0.42) -- (4.500,0.42) (6.420,0.42) -- (6.540,0.42) (10.830,0.42) -- (10.950,0.42)
  (13.320,0.42) -- (14.370,0.42) (0.000,0.84) -- (4.500,0.84) (4.950,0.84) -- (5.070,0.84)
  (9.360,0.84) -- (9.480,0.84) (13.320,0.84) -- (14.370,0.84) (0.000,1.26) -- (4.500,1.26)
  (4.950,1.26) -- (5.070,1.26) (9.360,1.26) -- (9.480,1.26) (13.320,1.26) -- (14.370,1.26)
  (0.000,1.68) -- (4.500,1.68) (6.420,1.68) -- (6.540,1.68) (10.830,1.68) -- (10.950,1.68)
  (13.320,1.68) -- (14.370,1.68) (0.000,2.10) -- (4.500,2.10) (8.790,2.10) -- (8.910,2.10)
  (13.200,2.10) -- (14.370,2.10);
\draw[gr]
  (5.070,0.840) .. controls (5.295,0.840) and (5.295,0.420) .. (5.520,0.420)
  (5.970,0.420) .. controls (6.195,0.420) and (6.195,0.840) .. (6.420,0.840)
  (6.540,0.420) .. controls (6.765,0.420) and (6.765,0.000) .. (6.990,0.000)
  (6.990,0.840) .. controls (7.215,0.840) and (7.215,0.420) .. (7.440,0.420)
  (7.890,0.420) .. controls (8.115,0.420) and (8.115,0.840) .. (8.340,0.840)
  (8.340,0.000) .. controls (8.565,0.000) and (8.565,0.420) .. (8.790,0.420)
  (9.480,0.840) .. controls (9.705,0.840) and (9.705,0.420) .. (9.930,0.420)
  (10.380,0.420) .. controls (10.605,0.420) and (10.605,0.840) .. (10.830,0.840)
  (10.950,0.420) .. controls (11.175,0.420) and (11.175,0.000) .. (11.400,0.000)
  (11.400,0.840) .. controls (11.625,0.840) and (11.625,0.420) .. (11.850,0.420)
  (12.300,0.420) .. controls (12.525,0.420) and (12.525,0.840) .. (12.750,0.840)
  (12.750,0.000) .. controls (12.975,0.000) and (12.975,0.420) .. (13.200,0.420);
\draw
  (4.500,1.260) .. controls (4.725,1.260) and (4.725,0.840) .. (4.950,0.840)
  (5.070,1.680) .. controls (5.295,1.680) and (5.295,1.260) .. (5.520,1.260)
  (5.520,1.260) .. controls (5.745,1.260) and (5.745,0.840) .. (5.970,0.840)
  (5.970,1.260) .. controls (6.195,1.260) and (6.195,1.680) .. (6.420,1.680)
  (6.540,2.100) .. controls (6.765,2.100) and (6.765,1.680) .. (6.990,1.680)
  (6.990,1.680) .. controls (7.215,1.680) and (7.215,1.260) .. (7.440,1.260)
  (7.440,1.260) .. controls (7.665,1.260) and (7.665,0.840) .. (7.890,0.840)
  (7.890,1.260) .. controls (8.115,1.260) and (8.115,1.680) .. (8.340,1.680)
  (8.340,1.680) .. controls (8.565,1.680) and (8.565,2.100) .. (8.790,2.100)
  (8.910,1.260) .. controls (9.135,1.260) and (9.135,0.840) .. (9.360,0.840)
  (9.480,1.680) .. controls (9.705,1.680) and (9.705,1.260) .. (9.930,1.260)
  (9.930,1.260) .. controls (10.155,1.260) and (10.155,0.840) .. (10.380,0.840)
  (10.380,1.260) .. controls (10.605,1.260) and (10.605,1.680) .. (10.830,1.680)
  (10.950,2.100) .. controls (11.175,2.100) and (11.175,1.680) .. (11.400,1.680)
  (11.400,1.680) .. controls (11.625,1.680) and (11.625,1.260) .. (11.850,1.260)
  (11.850,1.260) .. controls (12.075,1.260) and (12.075,0.840) .. (12.300,0.840)
  (12.300,1.260) .. controls (12.525,1.260) and (12.525,1.680) .. (12.750,1.680)
  (12.750,1.680) .. controls (12.975,1.680) and (12.975,2.100) .. (13.200,2.100);
\draw[gr,preaction={halo}]
  (5.070,1.260) .. controls (5.295,1.260) and (5.295,1.680) .. (5.520,1.680)
  (5.970,1.680) .. controls (6.195,1.680) and (6.195,1.260) .. (6.420,1.260)
  (6.540,1.680) .. controls (6.765,1.680) and (6.765,2.100) .. (6.990,2.100)
  (6.990,1.260) .. controls (7.215,1.260) and (7.215,1.680) .. (7.440,1.680)
  (7.890,1.680) .. controls (8.115,1.680) and (8.115,1.260) .. (8.340,1.260)
  (8.340,2.100) .. controls (8.565,2.100) and (8.565,1.680) .. (8.790,1.680)
  (9.480,1.260) .. controls (9.705,1.260) and (9.705,1.680) .. (9.930,1.680)
  (10.380,1.680) .. controls (10.605,1.680) and (10.605,1.260) .. (10.830,1.260)
  (10.950,1.680) .. controls (11.175,1.680) and (11.175,2.100) .. (11.400,2.100)
  (11.400,1.260) .. controls (11.625,1.260) and (11.625,1.680) .. (11.850,1.680)
  (12.300,1.680) .. controls (12.525,1.680) and (12.525,1.260) .. (12.750,1.260)
  (12.750,2.100) .. controls (12.975,2.100) and (12.975,1.680) .. (13.200,1.680);
\draw[preaction={halo}]
  (4.500,0.840) .. controls (4.725,0.840) and (4.725,1.260) .. (4.950,1.260)
  (5.070,0.420) .. controls (5.295,0.420) and (5.295,0.840) .. (5.520,0.840)
  (5.520,0.840) .. controls (5.745,0.840) and (5.745,1.260) .. (5.970,1.260)
  (5.970,0.840) .. controls (6.195,0.840) and (6.195,0.420) .. (6.420,0.420)
  (6.540,0.000) .. controls (6.765,0.000) and (6.765,0.420) .. (6.990,0.420)
  (6.990,0.420) .. controls (7.215,0.420) and (7.215,0.840) .. (7.440,0.840)
  (7.440,0.840) .. controls (7.665,0.840) and (7.665,1.260) .. (7.890,1.260)
  (7.890,0.840) .. controls (8.115,0.840) and (8.115,0.420) .. (8.340,0.420)
  (8.340,0.420) .. controls (8.565,0.420) and (8.565,0.000) .. (8.790,0.000)
  (8.910,0.840) .. controls (9.135,0.840) and (9.135,1.260) .. (9.360,1.260)
  (9.480,0.420) .. controls (9.705,0.420) and (9.705,0.840) .. (9.930,0.840)
  (9.930,0.840) .. controls (10.155,0.840) and (10.155,1.260) .. (10.380,1.260)
  (10.380,0.840) .. controls (10.605,0.840) and (10.605,0.420) .. (10.830,0.420)
  (10.950,0.000) .. controls (11.175,0.000) and (11.175,0.420) .. (11.400,0.420)
  (11.400,0.420) .. controls (11.625,0.420) and (11.625,0.840) .. (11.850,0.840)
  (11.850,0.840) .. controls (12.075,0.840) and (12.075,1.260) .. (12.300,1.260)
  (12.300,0.840) .. controls (12.525,0.840) and (12.525,0.420) .. (12.750,0.420)
  (12.750,0.420) .. controls (12.975,0.420) and (12.975,0.000) .. (13.200,0.000);
\draw[line width=0.5pt,fill=white] (0.130,-0.170) rectangle (0.830,1.010)
  node[midway] {$s_1$};
\draw[line width=0.5pt,fill=white] (0.960,-0.170) rectangle (1.660,1.010)
  node[midway] {$s_2$};
\draw[line width=0.5pt,fill=white] (1.790,1.090) rectangle (2.490,2.270)
  node[midway] {$s_3$};
\draw[line width=0.5pt,fill=white] (2.620,1.090) rectangle (3.320,2.270)
  node[midway] {$s_4$};
\draw[line width=0.5pt,fill=white] (3.450,1.090) rectangle (4.370,2.270)
  node[midway] {$\Delta_{3,3}^{-1}$};
\draw[line width=0.5pt,fill=white] (13.320,1.090) rectangle (14.240,2.270)
  node[midway] {$\Delta_{3,3}$};
\foreach \x/\i in {4.725/3,5.745/2,7.665/1,9.135/3,10.155/2,12.075/1}
  \node at (\x,-0.48) {\scriptsize $\bar c_{3,\i}$};
\foreach \k in {1,...,6} {
  \node[left] at (-0.03,{(\k-1)*0.42}) {\footnotesize $\k$};
  \node[right] at (14.40,{(\k-1)*0.42}) {\footnotesize $\k$};}
\end{tikzpicture}
\caption{The factorization $d(s_1,s_2,s_3,s_4)$ of~\eqref{eq:double} for
$m=3$, drawn as a braid on six strands, with its factors in order from left to
right and the positions of the strands numbered from bottom to top.  The boxes
stand for the factorizations $s_1,s_2$ and $\sh(s_3),\sh(s_4)$ of the full twist
on the strands in positions $1,2,3$ and $4,5,6$; each remaining factor is a
band, a conjugate of $\sigma_1$.  As
$c_{3,i}=\Delta_{3,3}^{-1}\bar c_{3,i}\Delta_{3,3}$, the conjugating
half-twists cancel between consecutive bands, so $\tilde r_3\cdot\tilde r_3$
appears as the bands $\bar c_{3,3},\bar c_{3,2},\bar c_{3,1}$, twice, between
$\Delta_{3,3}^{-1}$ and $\Delta_{3,3}$; within each band, the two strands it
joins, in positions $i$ and $7-i$, are black.  The product of all the factors
is $\Delta_6^2$.  In the geometric reading of \S\ref{subsec:pullback}
(Figure~\ref{fig:pullback}), the boxes correspond to the four quarters of the
curve and the bands to its tangencies at $v=0$; with $\Delta_3^2$ in each box,
this is the factorization $s$ of~\eqref{eq:s}.}
\label{fig:double}
\end{figure}
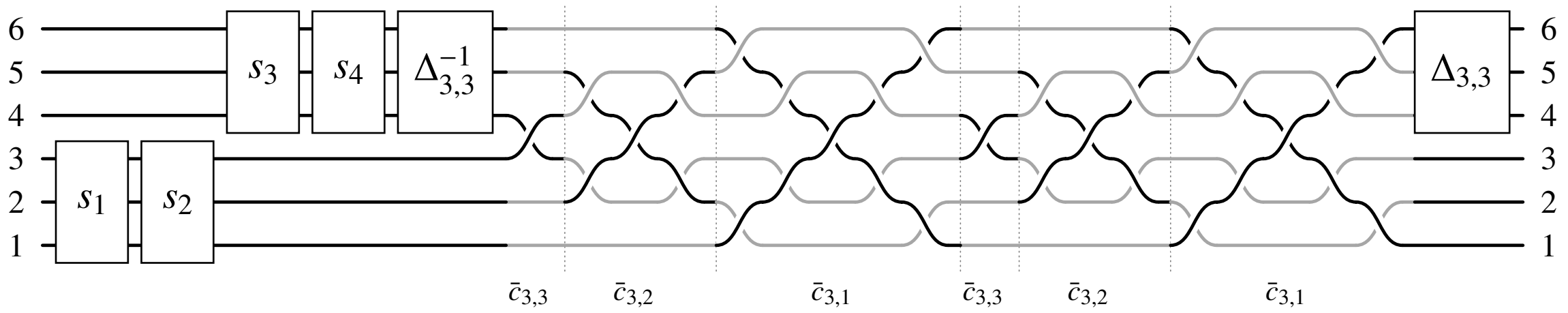

In $B_{2m}$ put $\sigma_{m,j}=\sigma_j\sigma_{2m-j}$,
\[
\bar c_{m,i}=(\sigma_{m,i}\sigma_{m,i+1}\cdots \sigma_{m,m-1})\,\sigma_m\,
(\sigma_{m,i}\sigma_{m,i+1}\cdots \sigma_{m,m-1})^{-1}\ (1\le i<m),
\qquad \bar c_{m,m}=\sigma_m,
\]
$c_{m,i}=\Delta_{m,m}^{-1}\bar c_{m,i}\Delta_{m,m}$, and
$\tilde r_m=c_{m,m}\cdot c_{m,m-1}\cdots c_{m,1}\in\cS_{2m}$.  Each $c_{m,i}$ is
a band with
\begin{equation}
\label{eq:cmi}
c_{m,i}(x_{m+i})=x_i,\qquad N(c_{m,i})=\nc{x_{m+i}^{-1}x_i}
\end{equation}
\cite[Assertion~1.12]{Kulikov-izv04a}.  For $s_1,\dots,s_4\in\cS_m$ with
$\alpha(s_i)=\Delta_m^2$, Kulikov's full twist formula
\cite[Thm.~3.2]{Kulikov-izv04a} states that
\begin{equation}
\label{eq:double}
d(s_1,s_2,s_3,s_4)=s_1\cdot s_2\cdot\sh(s_3)\cdot\sh(s_4)\cdot\tilde r_m\cdot
\tilde r_m
\end{equation}
satisfies $\alpha\bigl(d(s_1,s_2,s_3,s_4)\bigr)=\Delta_{2m}^2$ (see
Figure~\ref{fig:double} for $m=3$).  Write
$d(s)=d(s,s,s,s)$.  By \eqref{eq:cmi}, eliminating $x_{m+1},\dots,x_{2m}$ gives
\begin{equation}
\label{eq:Gdouble}
G\bigl(d(s_1,s_2,s_3,s_4)\bigr)\cong F_m\big/N(s_1)N(s_2)N(s_3)N(s_4),
\end{equation}
and in particular $G(d(s))\cong G(s)$ \cite[Assertions~1.13,
3.3]{Kulikov-izv04a}.

\section{The centrality obstruction}
\label{sec:central}

In this section we prove Proposition~\ref{prop:intro-central} and draw its
consequences: groups with trivial center (Corollary~\ref{cor:trivcenter}) and
link groups with vanishing linking numbers
(Proposition~\ref{prop:linking}) are not Hurwitz curve groups, whereas
$G\times\Z$ is a Hurwitz C-group for every C-group $G$
(Proposition~\ref{prop:product}).  In \S\ref{subsec:other} we recall, for
later reference, the consequences of centrality for $\ker\nu$ and the
Alexander invariant, and compare them on right-angled Artin groups
(Remark~\ref{rem:raag}).

\subsection{Proof of Proposition~\ref{prop:intro-central}}
\label{subsec:proofcentral}

\begin{lemma}
\label{lem:central}
If $\alpha(s)=\Delta_m^2$, then $l_m$ is central in $G(s)$.  In particular
$l_m$ is central in $\pi_1(\C^2\setminus H)$ for every affine Hurwitz curve
$H$ of degree $m$.
\end{lemma}

\begin{proof}
Each factor acts trivially on $G(s)$, hence so does their product
$\Delta_m^2$; by \eqref{eq:fulltwist}, $l_mx_il_m^{-1}=x_i$ in $G(s)$ for all
$i$.  The second statement follows from Theorem~\ref{thm:zvk} and
\eqref{eq:alpha}.
\end{proof}

This is noted in \cite[\S6]{Kulikov-izv04a} (see (27) there) and in
\cite[\S0]{Kulikova-izv05}.

\begin{proposition}
\label{prop:H1}
Let $\overline H$ have irreducible components $\overline H_1,\dots,\overline H_r$
of degrees $d_1,\dots,d_r$ and let $G=\pi_1(\C^2\setminus H)$.  Then
$H_1(G;\Z)\cong\Z^r$, with the meridians of $\overline H_c$ mapping to the
$c$-th basis vector, and the central element $l_m$ maps to
$(d_1,\dots,d_r)$.  In particular $l_m$ has infinite order and the image of
$\Zc(G)$ in $H_1(G;\Q)$ is non-zero.
\end{proposition}

\begin{proof}
Abelianizing \eqref{eq:Gs}, the relation $x_j=u_i(x_j)$ becomes
$x_j=x_{\pi_i(j)}$, where $\pi_i$ is the permutation of $u_i$.  Hence
$H_1(G;\Z)$ is free abelian on the orbits of the permutation group generated
by the $\pi_i$, which are the sets of strands belonging to the irreducible
components.  A generic line of the pencil meets $\overline H_c$ in $d_c$
points, so $l_m=x_1\cdots x_m\mapsto(d_1,\dots,d_r)$.
\end{proof}

Lemma~\ref{lem:central} and Proposition~\ref{prop:H1} together prove
Proposition~\ref{prop:intro-central}.

\subsection{Groups excluded by the obstruction}
\label{subsec:excluded}

\begin{corollary}
\label{cor:trivcenter}
A group $G$ such that $\Zc(G)\to H_1(G;\Q)$ is zero is not isomorphic to
$\pi_1(\C^2\setminus H)$ for any affine Hurwitz curve $H$.  This applies to
every group with trivial center: to $\BS(p,q)$ with
$\lvert p\rvert\neq\lvert q\rvert$, and to the fundamental group of any
hyperbolic $3$-manifold of finite volume.
\end{corollary}

Note that the obstruction does not depend on the C-structure: it rules out
$G$ as an abstract group.  That $\Zc(\BS(p,q))=1$ for
$\lvert p\rvert\neq\lvert q\rvert$ is standard: the axis or fixed subtree of a
central element of $\BS(p,q)=\la a,t\mid t^{-1}a^pt=a^q\ra$ in the
Bass--Serre tree is invariant under the vertex-transitive action, hence is
the whole tree; as the tree has valence $\lvert p\rvert+\lvert q\rvert\ge3$,
it is a fixed subtree, so the element lies in the normal core of $\la a\ra$,
which is trivial by Britton's lemma.

\begin{proposition}
\label{prop:linking}
Let $L\subset S^3$ be a link with $\mu\ge2$ components and vanishing pairwise
linking numbers, and $G=\pi_1(S^3\setminus L)$.  Then $\Zc(G)\subseteq G'$.
Consequently $G$ is not the group of any affine Hurwitz curve.
\end{proposition}

\begin{proof}
For a group $K$, write $\gamma_1K=K$ and $\gamma_{k+1}K=[K,\gamma_kK]$ for
its lower central series.  Let $F$ be the free group on meridians
$x_1,\dots,x_\mu$ of the components $L_1,\dots,L_\mu$ of $L$.  By Milnor
\cite{Milnor-proc57},
\[
G/\gamma_3G\cong\la x_1,\dots,x_\mu\mid[x_i,\lambda_i]\ (1\le i\le\mu),\
\gamma_3F\ra,
\]
where $\lambda_i\in F$ represents a longitude of $L_i$, and its image in
$F/\gamma_2F\cong\Z^\mu$ is given by the linking numbers $\lk(L_i,L_j)$.  As
these vanish, $\lambda_i\in\gamma_2F$, so $[x_i,\lambda_i]\in\gamma_3F$ and
$G/\gamma_3G\cong F/\gamma_3F$, the free nilpotent group of class $2$ on
$\mu\ge2$ generators.  Its center is $\gamma_2F/\gamma_3F$.  Hence every
central element of $G$ maps into $\gamma_2G/\gamma_3G$, and so lies in
$\gamma_2G=G'$.  Now apply Corollary~\ref{cor:trivcenter}.
\end{proof}

For knots the obstruction leaves exactly the torus knot groups: a knot group
with its Wirtinger C-structure is a Hurwitz C-group if and only if the knot
is a torus knot \cite[Thm.~6.3]{Kulikov-izv04a}.

\subsection{Products with \texorpdfstring{$\Z$}{Z}}
\label{subsec:products}

\begin{proposition}[Kulikova's product trick]
\label{prop:product}
Let $G$ be a C-group with C-generators $x_1,\dots,x_m$ generating $G$, and
let $z$ generate a free cyclic factor.  Then $G\times\la z\ra$ is a Hurwitz
C-group of degree $m+1$, with C-generators $x_1,\dots,x_m$ and
$x_{m+1}=l_m^{-1}z$.  Its irreducible components are those of $G$ together
with one new component, of degree $1$ in the degree-$(m+1)$ presentation.
\end{proposition}

\begin{proof}
From a C-presentation $\la x_1,\dots,x_m\mid R\ra$ of $G$ we get
$G\times\la z\ra=\la x_1,\dots,x_m,z\mid R,\ [z,x_i]\ (i\le m)\ra$.
Substituting $z=l_{m+1}=x_1\cdots x_{m+1}$ (Tietze) gives the presentation
$\la x_1,\dots,x_{m+1}\mid R,\ x_i=l_{m+1}^{-1}x_il_{m+1}\ (i\le m+1)\ra$, the
relation for $i=m+1$ being a consequence of the others.  All relations are
C-relations and $l_{m+1}$ is central.  The class of $x_{m+1}$ in $H_1$ is
$-\sum_{i\le m}[x_i]+[z]$, which is not the class of any $x_i$, $i\le m$.
\end{proof}

This is the argument of \cite[\S2]{Kulikova-izv05}, where it is applied to
$\BS(2,3)$; it works for every C-group.  Combined with \cite[Thm.~7.3]{Kulikov-usp07},
it shows that every C-group is a direct factor, with complement $\Z$, of the
group of some Hurwitz curve.

\subsection{Consequences of centrality}
\label{subsec:other}
Centrality of $l_m$ has consequences for $N$ and $A_0(G)$ that, unlike
Corollary~\ref{cor:trivcenter}, depend on the C-structure, but can be tested
without determining the center.  For the reader's orientation we list those
used below, together with the corresponding realization results.
\begin{enumerate}[label=(\arabic*),itemsep=2pt]
\item \label{k1}
\emph{Finite presentability.}  For every Hurwitz C-group, $N=\ker\nu$ is
finitely presented \cite[Thm.~0.2]{Kulikova-izv05}; for irreducible $G$ this says
that $G'$ is finitely presented \cite[Thm.~0.1]{Kulikova-izv05}.  This follows
from centrality: as $\nu(l_m)=m$, the subgroup $\nu^{-1}(m\Z)$, of index~$m$
in $G$, is $N\times\la l_m\ra$, so it is finitely presented, and so is its
quotient~$N$.
\item \label{k2}
\emph{The Alexander invariant} \cite[Thm.~10.1, Cor.~10.1]{Kulikov-usp07},
\cite{Greuel-Kulikov-izv05}.  For a Hurwitz curve of degree $m$ with $r$ irreducible
components, $A_0$ is a finitely generated abelian group, a quotient of
$\Z\oplus\bigl(\Lambda/(t^m-1)\bigr)^{m-2}$, the action of $t$ on
$A_0\otimes\C$ is semisimple, $\Delta(t)$ divides $(t-1)(t^m-1)^{m-2}$, and
the root $t=1$ has multiplicity exactly $r-1$.  Moreover, $t^m$ acts as the identity 
on $A_0$ \cite[Prop.~3.4]{Kulikov-usp07}.  These properties, too, follow from
centrality: $t^m$ is induced by conjugation by $x_1^m\in l_mN$, hence is the
identity; the classes of the elements $x_ix_1^{-1}$, $1<i<m$, and the
$t$-invariant class of $l_mx_1^{-m}\in N$ generate the $\Lambda$-module $A_0$;
and, $t$ being semisimple, the multiplicity of the root $t=1$ is the
dimension $b_1(G)-1$ of the coinvariants of $t$ on $A_0\otimes\C$, 
by the Wang sequence.  For irreducible curves, $A_0$ is
$(t-1)$-invertible, and conversely every Noetherian, $(t-1)$-invertible,
$t$-unipotent $\Lambda$-module occurs, with unipotence index dividing the
degree \cite[Thm.~0.3]{Kulikov-izv08}; the possible Alexander polynomials
are determined in \cite{Kulikov-izv06}.
\end{enumerate}
For the irreducible C-group $G_k=\BS(k,k+1)$ of \eqref{eq:Gk} below all three
constraints fail: its center is trivial; $G_k'$ is not finitely generated;
and its Alexander invariant $\Lambda/\bigl((k+1)t-k\bigr)\cong\Z[\tfrac{1}{k(k+1)}]$
\cite[Ex.~3.20]{Kulikov-izv08} is
not finitely generated over $\Z$ ($\Delta_{G_k}$ has the non-cyclotomic root
$k/(k+1)$).

\begin{remark}
\label{rem:raag}
Let $\Gamma$ be a finite simplicial graph with vertex set $V$, and
$n=\lvert V\rvert$.  The right-angled Artin group $A_\Gamma$, generated by $V$
with a relation $v=u^{-1}vu$ for each edge $\{u,v\}$, is a C-group with $n$
irreducible components, one for each vertex, and $N=\ker\nu$ is the
Bestvina--Brady group $N_\Gamma$ \cite{Bestvina-Brady-inv97,Papadima-Suciu-jlms07}.
If $u,v\in V$ are not adjacent, the retraction of $A_\Gamma$ onto
$\la u,v\ra\cong F_2$ that kills the other vertices maps central elements
to~$1$, so their images in $H_1(A_\Gamma;\Z)=\Z^V$ have vanishing $u$- and
$v$-coordinates.  Hence, by Corollary~\ref{cor:trivcenter}, and by
Proposition~\ref{prop:product} together with the realization theorem recalled
in \S\ref{subsec:cgroups}, $A_\Gamma$ is the group of an affine Hurwitz curve
if and only if $\Gamma$ is a cone, that is, has a vertex $w$ adjacent to all
others, so that $A_\Gamma=A_{\Gamma'}\times\la w\ra$, with $\Gamma'$ the subgraph
induced on $V\setminus\{w\}$; and, by Proposition~\ref{prop:H1}, the
C-structure above is a Hurwitz one only if $\Gamma$ is complete,
$A_\Gamma=\Z^n$ being then the group of $n$ lines in general position.  For a
cone, $\ker\nu$ in the C-structure of Proposition~\ref{prop:product} projects
isomorphically onto the kernel of the length homomorphism $A_{\Gamma'}\to\Z$
reduced modulo~$n$, a subgroup of index~$n$ in
$A_{\Gamma'}\cong N_\Gamma$ \cite[Ex.~2.5]{Papadima-Suciu-jlms07}.  For instance,
if $\Gamma'$ is discrete, $A_\Gamma=F_{n-1}\times\Z$ is the group of $n$
concurrent lines, whose meridians give this C-structure, and $\ker\nu$ is the
fundamental group, free of rank $(n-1)^2$, of their Milnor fiber.

The consequences \ref{k1} and \ref{k2} are strictly weaker than centrality,
and \ref{k2} does not imply \ref{k1}.  Indeed, $N_\Gamma$ is finitely presented
if and only if the flag complex of $\Gamma$ is simply connected
\cite{Bestvina-Brady-inv97}, and if $\Gamma$ is connected, then
$A_0(A_\Gamma)=H_1(N_\Gamma)\cong\Z^{n-1}$, with $t$ acting trivially
\cite[Props.~4.2 and~5.3]{Papadima-Suciu-jlms07}, so that the conditions
in~\ref{k2} hold, with $r=n$ and any $m\ge n$.  Thus, for the $4$-cycle,
$A_\Gamma=F_2\times F_2$ satisfies~\ref{k2} but not~\ref{k1}, while for the path
on four vertices $A_\Gamma$ satisfies both, $N_\Gamma$ being then free of
rank~$3$ \cite[Ex.~2.4]{Papadima-Suciu-jlms07}, but is not the group of a
Hurwitz curve.
\end{remark}

\section{Kulikova's example and Kulikov's algorithm}
\label{sec:kulikova}

We recall how Kulikova realizes $\BS(k,k+1)\times\Z$ as a Hurwitz C-group of
degree $3$, and explain why Kulikov's realization algorithm, applied to this
C-group for $k=2$, produces a Hurwitz curve of degree $768$
(\S\ref{subsec:768}).  In \S\ref{subsec:samegroup} we use the same algorithm
to show that the C-structure can distinguish Hurwitz curves whose groups are
isomorphic.

For $k\ge1$ put
\begin{equation}
\label{eq:Gk}
G_k=\la x_1,x_2\mid (x_1^{-1}x_2)^k\,x_1\,(x_1^{-1}x_2)^{-k}=x_2\ra
\cong\BS(k,k+1),
\end{equation}
via $t=x_1$, $a=x_1^{-1}x_2$, under which the relation becomes
$t^{-1}a^kt=a^{k+1}$ \cite[\S7.2]{Kulikov-usp07}.  The groups $G_k$ are irreducible
C-groups, non-Hopfian for $k\ge2$ \cite{Baumslag-Solitar-bams62}, hence not residually
finite; $G_1=\BS(1,2)$ is residually finite, being a finitely generated linear
group.  Since \eqref{eq:Gk} has two generators and one relation, $G_k$ has
C-deficiency one, so it is the group of a knotted $2$-sphere in $S^4$
\cite[\S3.4, Ex.~3.20]{Kulikov-izv08}.
By Proposition~\ref{prop:product},
\begin{equation}
\label{eq:tildeGk}
\widetilde G_k=G_k\times\la z\ra
=\la x_1,x_2,x_3\mid (x_1^{-1}x_2)^kx_1(x_1^{-1}x_2)^{-k}=x_2,\
[x_1x_2x_3,x_i]=1\ (i=1,2,3)\ra
\end{equation}
is a Hurwitz C-group of degree $3$, with $z=x_1x_2x_3$ and two irreducible
components: $\{x_1,x_2\}$ and $\{x_3\}$ \cite[\S2]{Kulikova-izv05},
\cite[\S7.2]{Kulikov-usp07}.  Kulikova \cite[Cor.~0.7, 0.8]{Kulikova-izv05} (see also
\cite[Thms.~7.5, 7.6]{Kulikov-usp07}) concludes that there is a two-component
Hurwitz curve, all of whose singular points are ordinary triple points, with
non-Hopfian $\pi_1(\C^2\setminus H)$ and $\pi_1(\CP^2\setminus\overline H)$,
and a pseudoholomorphic curve of the same kind with non-Hopfian
$\pi_1(\CP^2\setminus\overline H)$.  The projective groups are
$G_{3,j}=\widetilde G_2/\nc{l_3^{\,j}}\cong\BS(2,3)\times\Z/j$ (the cyclic
factor has order $j$, not $3j$ as printed in \cite[\S2]{Kulikova-izv05}; this does not
affect the conclusion).

\subsection{The degree produced by the algorithm}
\label{subsec:768}
The realization goes through \cite[Thm.~3.6]{Kulikov-izv04a}.  Its proof starts
from $\Delta_m^2$ and installs the non-central C-relations one at a time; a
relation $x_i=w^{-1}x_jw$ whose conjugating word has letter length $l$ is
installed by passing from the current factorization $s_n\in\cS_{M_n}$ to
$\tilde s_n=d^{\,l+2}(s_n)$ and then to
$s_{n+1}=d\bigl(\lambda(b)(d(\tilde s_n)),d(\tilde s_n),d(\tilde s_n),d(\tilde
s_n)\bigr)$ for a suitable band $b$; thus $M_{n+1}=2^{l+4}M_n$ exactly.  For
\eqref{eq:tildeGk} with $k=2$ there is one non-central relation, with
conjugator $(x_1^{-1}x_2)^2$ of length $l=4$, so the curve produced has
degree exactly $3\cdot2^8=768$.  The source of the cost is visible in the proof:
the word $w$ must be ``moved apart'' \cite[Assertion~1.7]{Kulikov-izv04a} so that its letters
occupy distinct strings, and the band lemma \cite[Lem.~3.4]{Kulikov-izv04a} uses
an extra doubling to create a partner string for the endpoint of the band.
Section~\ref{sec:explicit} removes both costs for \eqref{eq:tildeGk}.

We record a property of the algorithm used in \S\ref{subsec:samegroup}.  Let
$\mathcal A_m$ be the smallest set of factorizations that contains
$\Delta_m^2\in\cS_m$ and, with each $s\in\cS_M$, contains $d(s)$ and every
$s'=d\bigl(\lambda(b)s,s,s,s\bigr)$ in which $b\in B_M$ is a band whose
endpoint loops have the same image in $G(s')$.  The factorizations produced by
the algorithm, and their iterated doubles, lie in $\mathcal A_m$: the band
$b=g\sigma_1g^{-1}$ of \cite[Lem.~3.4(i)]{Kulikov-izv04a} has endpoint loops
$g(x_1)=y$ and $g(x_2)=x_i$, and $x_i=y$ holds in the resulting group, being
the relation that the step installs.

\begin{lemma}
\label{lem:meridprod}
Let $s\in\cS_M$ be in $\mathcal A_m$, $m\ge2$, and let $\rho\colon F_M\to G(s)$ be the
projection.  Then $\rho(x_j)=\rho(x_{j'})$ whenever $j\equiv j'\pmod m$, and
every factor of $s$ conjugate to $\Delta_m^2$ has the form
$h\,\Delta_{m,i}^2h^{-1}$, with $h\in B_M$ and $0\le i\le M-m$, such that
$\rho\bigl(h(x_{i+1}\cdots x_{i+m})\bigr)$ is conjugate to $\rho(l_m)$.
\end{lemma}

\begin{proof}
We argue by induction along the construction of $s$; for $s=\Delta_m^2$ take
$h=1$ and $i=0$.  Let $s\in\cS_M$ satisfy the assertions, let
$s'=d(s_1,s,s,s)$ with $s_1=s$ or $s_1=\lambda(b)s$ as above, and let
$\rho'\colon F_{2M}\to G(s')$ be the projection.  By \eqref{eq:Gdouble},
$\rho'(x_{M+j})=\rho'(x_j)$ for $j\le M$, and the restriction of $\rho'$ to
$F_M$ factors through $\rho$, as $N(s)\subseteq N(s_1)N(s)$.  This gives the
first assertion for $s'$, and shows that $\rho'(\sh(w))=\rho'(w)$ for $w\in
F_M$.  By \eqref{eq:double}, the factors of $s'$ conjugate to $\Delta_m^2$ are
those of $s_1$ and of $s$ and the shifts of those of $s$, since the bands of
$\tilde r_M\cdot\tilde r_M$ have exponent sum $1<m(m-1)$.  A factor
$h\,\Delta_{m,i}^2h^{-1}$ of $s$ and its shift
$\sh(h)\,\Delta_{m,M+i}^2\sh(h)^{-1}$ give $w=h(x_{i+1}\cdots x_{i+m})$ and
$\sh(w)$, whose images under $\rho'$ are conjugate to $\rho'(l_m)$.  If
$s_1=\lambda(b)s$, the corresponding factor of $s_1$ is
$(bh)\,\Delta_{m,i}^2(bh)^{-1}$, which gives $b(w)$, and $\rho'(b(w))=\rho'(w)$ by
Lemma~\ref{lem:bandrel}, since the endpoint loops of $b$ have the same image
under $\rho'$.
\end{proof}

\subsection{Curves with isomorphic groups}
\label{subsec:samegroup}
The C-structure is part of the data of a Hurwitz curve group, and it can
distinguish curves whose groups are isomorphic; compare the line
arrangements of \cite[\S7.4]{Cohen-Suciu-cmh97}, which have isomorphic groups
but braid monodromies of different types.  We illustrate this with the
meridians of the trefoil group found in \cite{Suciu-mpcps85}.  Let
$\pi=\la s,t\mid sts=tst\ra\cong B_3$ be the trefoil knot group, and put
$u=t^{-1}s$.  For $k\ge1$, Suciu's ribbon $2$-knot $R_k\subset S^4$ has group
\begin{equation}
\label{eq:Rk}
\pi\cong\la x_1,x_2,x_3\mid x_1=V^{-1}x_3V,\ x_2=W^{-1}x_1W\ra,\qquad
V=(x_1x_2^{-1})^{k-1}x_3^{-1},\ W=x_3x_1x_2^{-1}x_3^{-1},
\end{equation}
a C-presentation in which the meridian $x_3$ corresponds to $t_k=u^kt$
\cite[\S4]{Suciu-mpcps85}, \cite[Prop.~2.1]{Kanenobu-Sumi-jktr20}; $R_1$ is the
spun trefoil.  Under the isomorphism $\pi/\Zc(\pi)\cong\PSL(2,\Z)$,
no automorphism takes the image of $t_k$ to a conjugate of the image
of $t_l^{\pm1}$ when $k\neq l$ \cite[Lem.~4.1]{Suciu-mpcps85}.

\begin{proposition}
\label{prop:samegroup}
Let $r\ge1$ and $N=2r+11$.  There are Hurwitz curves
$\overline Q_1,\dots,\overline Q_r\subset\CP^2$ of the same degree
$4\cdot2^N$, each with two irreducible components, of degrees $3\cdot2^N$
and $2^N$, each with $4^N$ ordinary quadruple points as its only
singularities, and each with $4^{N+1}-2^{N+2}$ simple tangencies with the
lines of the pencil, such that
\[
\pi_1(\C^2\setminus Q_k)\cong\pi\times\Z,\qquad
\pi_1(\CP^2\setminus\overline Q_k)\cong\pi\times\Z/2^N
\]
for all $k$, but no homeomorphism of $\CP^2$ carries $\overline Q_k$ onto
$\overline Q_l$ for $k\neq l$.
\end{proposition}

\begin{proof}
By Proposition~\ref{prop:product}, \eqref{eq:Rk} gives $\pi\times\la z\ra$ the
structure of a Hurwitz C-group of degree $4$, with C-generators $x_1,x_2,x_3$
and $x_4=(x_1x_2x_3)^{-1}z$.  The algorithm of \S\ref{subsec:768}
\cite[Thm.~3.6]{Kulikov-izv04a} realizes it by a factorization obtained from
$\Delta_4^2$ through $N_k$ doublings, some of them twisted by bands.  The
central relations $x_i=l_4^{-1}x_il_4$ are provided by $\Delta_4^2$ itself
\cite[Assertion~3.5]{Kulikov-izv04a}, and the two relations of
\eqref{eq:Rk} are installed one at a time, at a cost of $l+4$ doublings each,
where $l$ is the letter length of the conjugating word.  As $V$ and $W$ have
letter lengths $2k-1$ and $4$, we get $N_k=(2k+3)+8=2k+11$.  By
\eqref{eq:double}, the factors of $d(s_1,s_2,s_3,s_4)$ are those of the $s_i$,
shifted or not, together with the $2m$ bands of $\tilde r_m\cdot\tilde r_m$,
and $\lambda(b)$ does not change the types of the factors.  Hence the result
has $4^{N_k}$ factors conjugate to $\Delta_4^2$, that is, ordinary quadruple
points, and all other factors are bands.  A further doubling $d(s)$ does not
change the C-group, by \eqref{eq:Gdouble}; doubling the factorization for
$R_k$ a further $2(r-k)$ times, we may assume that $N_k=N=2r+11$ for
$k=1,\dots,r$.  The product of the factorization is $\Delta_M^2$, with
$M=4\cdot2^N$, whose exponent sum is $M(M-1)$; as $\Delta_4^2$ has
exponent sum $12$ and each band exponent sum $1$, the number of bands is
$M(M-1)-12\cdot4^N=4^{N+1}-2^{N+2}$.  The strings of $x_1,x_2,x_3$ and those of
$x_4$ form the two components, of degrees $3\cdot2^N$ and $2^N$.  In
$G(d(s))$ we have $l_{2m}=l_m^2$, so $l_{4\cdot2^N}=z^{2^N}$, and
Theorem~\ref{thm:zvk} gives the projective group $\pi\times\Z/2^N$.

Suppose that a homeomorphism $f$ of $\CP^2$ carries $\overline Q_k$ onto
$\overline Q_l$.  As the two components have different degrees, $f$ maps the
first component of $\overline Q_k$ onto that of $\overline Q_l$, so the
induced isomorphism of projective groups takes the conjugacy class of a
meridian of the first to that of a meridian of the second, or of its
inverse.  Identifying both groups with $\pi\times\Z/2^N$, these meridians
are the conjugates of $t_k$ and of $t_l$, respectively.  Dividing by the
centers, which are characteristic, we obtain an automorphism of
$\pi/\Zc(\pi)\cong\PSL(2,\Z)$ taking the image of $t_k$ to a conjugate
of the image of $t_l^{\pm1}$, so $k=l$.
\end{proof}

The $2$-knots $R_k$ are also distinguished by the fundamental groups of
their $2$-fold branched coverings \cite{Kanenobu-Sumi-jktr20}.  For the curves
$\overline Q_k$ the corresponding groups carry no information.

\begin{proposition}
\label{prop:2fold}
Let $\overline Q_1,\dots,\overline Q_r$ be the curves constructed in the proof
of Proposition~\textup{\ref{prop:samegroup}}.  Since their degree is even,
there is a $2$-fold covering $Y_k\to\CP^2$ branched along $\overline Q_k$.
Its singular points lie over the quadruple points of $\overline Q_k$; let
$\overline X_k$ be its minimal resolution, and $E_k\subset\overline X_k$ the
union of the exceptional curves.  Then, for every $k$,
\[
\pi_1(\overline X_k)\cong\Z/3,\qquad
\pi_1(\overline X_k\setminus E_k)\cong\Z/6.
\]
\end{proposition}

\begin{proof}
Write $G=\pi_1(\CP^2\setminus\overline Q_k)=\pi\times\la z\ra$, with $z$ of
order $2^N$.  The covering $Y_k$ corresponds to the homomorphism $G\to\Z/2$
sending all meridians to $1$.  Its kernel is $K=\pi^{(2)}\times\la z\ra$,
where $\pi^{(2)}$ is the subgroup of index two of $\pi$, since
$z=x_1x_2x_3x_4$ has total degree $4$.  Away from the quadruple points,
$\overline Q_k$ is smooth, so $Y_k$ minus its singular points is a smooth
$4$-manifold, with fundamental group $\Gamma=K/\nc{\mu^2}$, the normal
closure being taken in $G$ over all meridians $\mu$.  The meridians are the
conjugates of $x_1,\dots,x_4$, and $x_1,x_2,x_3$ are conjugate in $\pi$, so
$\Gamma=K/\nc{x_1^2,x_4^2}$.

By \cite[Thm.~6.2]{Kulikov-izv04a}, near each quadruple point
$\overline Q_k$ is given, in local analytic coordinates, by $v^4=u^4$, a
union of four lines; so $Y_k$ has there the simple elliptic singularity
$w^2=u^4-v^4$.  Its minimal resolution
replaces it by an elliptic curve $E$, the $2$-fold covering of the
exceptional line of the blow-up branched at the four tangent directions,
with $E\cdot E=-2$.  The link is the pullback to $E$ of the Hopf circle
bundle, so its fiber is a lift of the Hopf fiber, that is, of the loop
around the four branches in a nearby line of the pencil, the product of the
four local meridians.  Now $\overline X_k\setminus E_k$ is $Y_k$ minus its
singular points, so $\pi_1(\overline X_k\setminus E_k)=\Gamma$.  Since
$\pi_1(E)$ is the quotient of the fundamental group of the link by the
fiber, the van Kampen theorem shows that $\pi_1(\overline X_k)$ is the
quotient of $\Gamma$ by the normal closure of these products.

Each such product is conjugate to $z$.  Indeed, the braid monodromy
factorization of $\overline Q_k$ lies in the class $\mathcal A_4$ of
\S\ref{subsec:768}, so by Lemma~\ref{lem:meridprod} each of its factors
conjugate to $\Delta_4^2$ has the form $h\,\Delta_{4,i}^2h^{-1}$, with
$h(x_{i+1}\cdots x_{i+4})$ conjugate in $G$ to $x_1x_2x_3x_4=z$.  Under the
identification of \S\ref{subsec:hurwitz}, $\Delta_{4,i}^2$ is the Dehn twist
about a circle enclosing the points $i+1,\dots,i+4$, so this factor is the
Dehn twist about the image $c$ of that circle under $h$.  It is also the
local monodromy of the corresponding quadruple point, transported to the base
fiber, that is, the Dehn twist in the same direction about the transport $c'$
of a small circle enclosing the four branches in a nearby fiber.  Both
circles enclose $4$ of the $4\cdot2^N$ points of the base fiber, so they are
essential, and a Dehn twist about an essential circle determines the circle
up to isotopy \cite[Fact~3.6]{Farb-Margalit-12}; so $c$ and $c'$ are isotopic, 
and the product of the four local meridians, a loop once around $c'$, is 
conjugate to $h(x_{i+1}\cdots x_{i+4})$, hence to $z$.
As $z$ is central, $\pi_1(\overline X_k)=\Gamma/\la z\ra$.

Killing $z$ turns $x_4^2$ into $(x_1x_2x_3)^2$, so
$\pi_1(\overline X_k)=\pi^{(2)}/\nc{x_1^2,(x_1x_2x_3)^2}$.  Here
$\Pi=\pi^{(2)}/\nc{x_1^2}$ is the group of the $2$-fold branched covering of
$R_k$.  By the Reidemeister--Schreier method, $\pi^{(2)}$ is generated by
the elements $x_ix_1^{-1}$ and $x_1x_i$, $i\le3$; as
$x_1x_i\equiv(x_ix_1^{-1})^{-1}$ modulo the squares of the meridians, $\Pi$ is
generated by $p=x_2x_1^{-1}$ and $q=x_3x_1^{-1}$.
Conjugation by $x_1$ induces the covering involution, which inverts $p$ and
$q$.  Hence
\[
(x_1x_2x_3)^2=x_1\,p\,x_1\,q\,x_1\cdot x_1\,p\,x_1\,q\,x_1
=p^{-1}q\,x_1\,p^{-1}q\,x_1=p^{-1}qpq^{-1}
\]
in $\Pi$.  Killing this commutator of the two generators makes $\Pi$ abelian,
so $\pi_1(\overline X_k)\cong H_1(\Pi)$.  Writing
$\pi^{(2)}=\pi'\rtimes\la x_1^2\ra$, we get $H_1(\Pi)=A/(t^2-1)A$, where
$A=\pi'/\pi''\cong\Z[t^{\pm1}]/(t^2-t+1)$ is the Alexander module of $\pi$,
which does not depend on the choice of meridian.  Hence $H_1(\Pi)\cong\Z/3$,
in accordance with the presentation of $\Pi$ in
\cite[Prop.~2.2]{Kanenobu-Sumi-jktr20}.

Finally, $z$ is central in $\Gamma$ and $\Gamma/\la z\ra\cong\Z/3$ is cyclic,
so $\Gamma$ is abelian.  In $\Gamma$, the relation $x_4^2=1$ reads
$(x_1x_2x_3)^2=z^2$, and the left side is a commutator, so $z^2=1$.  Thus
$\Gamma$ is a quotient of $\Z/3\times\Z/2$.  It maps onto $\Z/3$, and onto
$\Z/2$ by $K\to\la z\ra\to\Z/2$, which kills $x_1^2$, $x_4^2$ and their
conjugates.  Hence $\Gamma\cong\Z/6$.
\end{proof}

\begin{remark}
\label{rem:zariski}
The trefoil group also appears in Zariski's classical example, but the
relationship with the curves $\overline Q_k$ is tenuous.  Let
$\overline Z\subset\CP^2$ be the sextic $Q^3+C^2=0$, where the conic $Q=0$
and the cubic $C=0$ meet transversally, so that its only singularities are
six cusps on the conic.  Its affine group is $\pi\cong B_3$ with the
classical C-structure, in which the meridians are the conjugates of $s$
\cite{Zariski-ajm29}, \cite[Ex.~4.6]{Kulikov-izv08}.  The curves
$\overline Q_k$ are reducible, of degree $4\cdot2^N$, and their affine group
is $\pi\times\Z$, with the C-structure on $\pi$ coming from $R_k$; this is
the classical one only for $k=1$, when $t_1=t^{-1}st$.  The $2$-fold
branched coverings also differ.  For $\overline Z$, the minimal resolution
is a K3 surface, hence simply connected, while the complement of its
exceptional curves has $H_1\cong\Z/3$ \cite[Thm.~0.5, Ex.~4.6]{Kulikov-izv08}.
For $\overline Q_k$ the corresponding groups are $\Z/3$ and $\Z/6$, by
Proposition~\ref{prop:2fold}.  In both cases this $\Z/3$ is $A/(t^2-1)A$,
where $A=\Z[t^{\pm1}]/(t^2-t+1)$ is the Alexander module of the trefoil
group, but only for $\overline Q_k$ does it survive the resolution.
\end{remark}

\section{An explicit factorization of degree twelve}
\label{sec:explicit}

This section proves Theorem~\ref{thm:explicit}.  The key input is a
refinement of Kulikov's band lemma in which a single doubling installs a
C-relation, provided one endpoint loop of the band is a standard generator
with a partner string fixed by the band (\S\ref{subsec:bandlemma}).  
The factorizations $T_k$ are written down in
\S\ref{subsec:factorizations-explicit}, and their group, singularities and
components are determined in \S\ref{subsec:proofB}.

\subsection{A one-doubling band lemma}
\label{subsec:bandlemma}
Throughout this subsection $M\ge2$ is an integer, and we use the notation of
\S\ref{sec:prelim}: $\cS_M$ is the semigroup of factorizations in $B_M$,
with product map $\alpha\colon\cS_M\to B_M$; $\lambda(b)$ is simultaneous
conjugation by $b\in B_M$; $G(s)=F_M/N(s)$ is the group~\eqref{eq:Gs} of
$s\in\cS_M$; $d$ is Kulikov's doubling~\eqref{eq:double}, which takes four
factorizations of $\Delta_M^2$ to one of $\Delta_{2M}^2$; and a band
$b=g\,\sigma_1g^{-1}$ has endpoint loops $g(x_1)$ and $g(x_2)$, as in
Lemma~\ref{lem:bandrel}.

\begin{lemma}
\label{lem:band}
Let $s\in\cS_M$ with $\alpha(s)=\Delta_M^2$, and let $b=g\,\sigma_1g^{-1}\in B_M$ 
be a band.  Suppose there are indices $i\neq p$ such that
\begin{enumerate}[label=(\alph*),itemsep=1pt]
\item \label{h1}
One endpoint loop of $b$ is the standard generator $x_i$; call the other
one $y$.
\item \label{h2}
$b(x_p)=x_p$ in $F_M$.
\item \label{h3}
$x_p=x_i$ in $G(s)$.
\end{enumerate}
Then $T=d\bigl(\lambda(b)(s),s,s,s\bigr)\in\cS_{2M}$ satisfies
$\alpha(T)=\Delta_{2M}^2$, and
\[
G(T)\cong G(s)\big/\nc{x_i^{-1}y}
\]
by an isomorphism sending $x_j$ and $x_{M+j}$ to the class of $x_j$.
\end{lemma}

\begin{proof}
Since $\alpha(\lambda(b)s)=b\Delta_M^2b^{-1}=\Delta_M^2$, the first
assertion is \eqref{eq:double}.  By \eqref{eq:Gdouble} and \eqref{eq:Nconj},
$G(T)\cong F_M/N(s)\,b(N(s))$.  Let $K=\nc{x_i^{-1}y}$.  By
Lemma~\ref{lem:bandrel}, $b(w)\in wK$ for every $w\in F_M$, so
$b(N(s))\subseteq N(s)K$.  Conversely, $x_p^{-1}x_i\in N(s)$ by \ref{h3}, hence
$x_p^{-1}b(x_i)=b(x_p^{-1}x_i)\in b(N(s))$ by \ref{h2}.  If $x_i=g(x_2)$, then
$b(x_i)=g(\sigma_1(x_2))=g(x_1)=y$; if $x_i=g(x_1)$, then
$b(x_i)=g(x_1x_2x_1^{-1})=x_iyx_i^{-1}$.  Multiplying on the left by $x_i^{-1}x_p\in N(s)$
gives $x_i^{-1}y\in N(s)\,b(N(s))$ in the first case and $yx_i^{-1}\in N(s)\,b(N(s))$
in the second; as $N(s)\,b(N(s))$ is normal, $K\subseteq N(s)\,b(N(s))$ in
both cases.  Hence $N(s)\,b(N(s))=N(s)K$.
\end{proof}

\begin{remark}
\label{rem:kulikov34}
Kulikov's band lemma \cite[Lem.~3.4]{Kulikov-izv04a} takes
$\bar s=d\bigl(\lambda(b)(d(s)),d(s),d(s),d(s)\bigr)\in\cS_{4M}$ with
$g(x_2)=x_i$: there, hypotheses \ref{h2} and \ref{h3} hold automatically for
$p=M+i$, because $b\in B_M$ fixes $x_{M+i}$ and $c_{M,i}$ identifies
$x_{M+i}$ with $x_i$.  Lemma~\ref{lem:band} saves this doubling whenever a
partner string $p$ already exists in $s$.  Applied as stated to
$s=d(\Delta_3^2)$, Kulikov's lemma with the bands below gives degree $24$.
Applied to $s=\Delta_3^2$ it would give degree $12$, but no band then
produces a non-abelian group: if $b=g\sigma_1g^{-1}\in B_3$, then
$(y_1,y_2,y_3)=(g(x_1),g(x_2),g(x_3))$ is a free basis with
$y_1y_2y_3=l_3$, $N(\Delta_3^2)$ is the normal closure of the commutators $[y_j,l_3]$, and
$G(\Delta_3^2)/\nc{y_1^{-1}y_2}\cong\la y_1,y_3\mid y_1^2y_3\text{ central}\ra
\cong\Z^2$.  Lemma~\ref{lem:band} does not apply to $s=\Delta_M^2$ either:
there is no partner string, and in any case $\lambda(b)(\Delta_M^2)=\Delta_M^2$.
\end{remark}

\subsection{The factorizations}
\label{subsec:factorizations-explicit}
\begin{table}[ht]
\centering
\renewcommand{\arraystretch}{1.25}
\begin{tabular}{c|l|l}
& factor & relations \\ \hline
$1,2$ & $\Delta_3^2=(\sigma_1\sigma_2)^3$ & $[x_i,x_1x_2x_3]=1$, $i=1,2,3$\\
$3,4$ & $\sh(\Delta_3^2)=(\sigma_4\sigma_5)^3$ & $[x_i,x_4x_5x_6]=1$, $i=4,5,6$\\
$5,8$ & $c_{3,3}=\Delta_{3,3}^{-1}\,\sigma_3\,\Delta_{3,3}$ & $x_6=x_3$\\
$6,9$ & $c_{3,2}=\Delta_{3,3}^{-1}(\sigma_2\sigma_4)\,\sigma_3\,(\sigma_2\sigma_4)^{-1}
\Delta_{3,3}$ & $x_5=x_2$\\
$7,10$ & $c_{3,1}=\Delta_{3,3}^{-1}(\sigma_1\sigma_5\sigma_2\sigma_4)\,
\sigma_3\,(\sigma_1\sigma_5\sigma_2\sigma_4)^{-1}
\Delta_{3,3}$ & $x_4=x_1$
\end{tabular}
\medskip
\caption{The factorization $s=d(\Delta_3^2)$ of $\Delta_6^2$.}
\label{tab:s}
\end{table}

Let
\begin{equation}
\label{eq:s}
s=d(\Delta_3^2)=\Delta_3^2\cdot\Delta_3^2\cdot\sh(\Delta_3^2)\cdot\sh(\Delta_3^2)
\cdot c_{3,3}\cdot c_{3,2}\cdot c_{3,1}\cdot c_{3,3}\cdot c_{3,2}\cdot c_{3,1}
\in\cS_6 ,
\end{equation}
whose ten factors are listed in Table~\ref{tab:s}, together with the
relations they contribute to $G(s)$; here $\Delta_{3,3}=\sigma_4\sigma_5\sigma_4$.
The ten factors multiply to the full twist, $\alpha(s)=\Delta_6^2$;
Figure~\ref{fig:double}, with $\Delta_3^2$ in each box, draws $s$ as a braid.
Geometrically, $s$ has the type of the braid monodromy factorization of the
curve $C_3$ of \S\ref{subsec:pencil}, formed by the three conics
$V^2-1=\lambda(U^2-1)$, $\lambda=2,3,4$, of a pencil with four base points;
$C_3$ is the preimage of three concurrent lines under $(U,V)\mapsto(U^2,V^2)$.
Its four ordinary triple points $(\pm1,\pm1)$ give the four conjugates of
$\Delta_3^2$, and its six simple tangencies with the vertical lines, two on
each conic, at the points with $V=0$, give the six bands.
By \eqref{eq:Gdouble} and \eqref{eq:fulltwist},
\begin{equation}
\label{eq:Gs3}
G(s)\cong G(\Delta_3^2)=\la x_1,x_2,x_3\mid [x_i,x_1x_2x_3]\ (i=1,2,3)\ra
\cong F(x_1,x_2)\times\la z\ra,\qquad z=x_1x_2x_3,
\end{equation}
with $x_{j+3}=x_j$.  With $g_k$, $b_k$ as in Theorem~\ref{thm:explicit}, a direct
computation with \eqref{eq:action} gives the endpoint loops:
\begin{align*}
k=2:\quad & g_2(x_1)=x_1,\\
& g_2(x_2)=\eta_2\coloneqq x_2x_3x_2^{-1}x_4x_6^{-1}x_5^{-1}x_4^{-1}x_2^{-1}
x_4x_5x_4^{-1}x_2x_4x_5x_6x_4^{-1}x_2x_3^{-1}x_2^{-1},\\
k=1:\quad & g_1(x_1)=x_2,\qquad
g_1(x_2)=\eta_1\coloneqq x_2^{-1}x_1x_2x_3x_2^{-1}x_1^{-1}x_2x_4x_2^{-1}
x_1x_2x_3^{-1}x_2^{-1}x_1^{-1}x_2,
\end{align*}
and $b_2(x_4)=x_4$, $b_1(x_5)=x_5$.  So the hypotheses of
Lemma~\ref{lem:band} hold with $(i,p)=(1,4)$ for $k=2$ and $(i,p)=(2,5)$ for
$k=1$, \ref{h3} being part of \eqref{eq:Gs3}.

\subsection{Proof of Theorem~\ref{thm:explicit}}
\label{subsec:proofB}
\ref{a1} The factor count is immediate from \eqref{eq:double}: $T_k$ consists of
four copies of the ten factors of $s$ (conjugated by $b_k$, or shifted) and
the twelve bands of $\tilde r_6\cdot\tilde r_6$; the triple-point factors are
the $4\cdot4=16$ conjugates of $\Delta_3^2$.  The identity
$\alpha(T_k)=\Delta_{12}^2$ is \eqref{eq:double}; as a consistency check,
the exponent sums add up to $16\cdot6+36=132=12\cdot11$.

\ref{a2} By Lemma~\ref{lem:band} and \eqref{eq:Gs3}, $G(T_k)$ is the quotient of
$F(x_1,x_2)\times\la z\ra$ by the relation $x_i=\eta_k$, where $\eta_k$ is
rewritten using $x_{j+3}=x_j$ and $x_3=(x_1x_2)^{-1}z$.  The $z$-exponent of
$\eta_k$ is zero ($\eta_k$ is conjugate to $x_5$, resp.\ $x_4$), so its
image in $G(s)$ lies in $F(x_1,x_2)=F(t,a)$, where $t=x_1$ and
$a=x_1^{-1}x_2$; explicitly,
\[
k=2:\ \ \eta_2=t^{-1}a^{-2}t\,a^3\,t,\qquad
k=1:\ \ \eta_1=(ta)^{-1}\bigl(a^{-1}t^{-1}a\,t\,a^{-1}\bigr)(ta)^2 .
\]
Thus $x_1=\eta_2$ reads $t^{-1}a^2t=a^3$, and $x_2=\eta_1$ (i.e., $ta=\eta_1$) reads
$t^{-1}at=a^2$.  Hence
$G(T_k)\cong\la a,t\mid t^{-1}a^kt=a^{k+1}\ra\times\la z\ra$, with $x_1\mapsto t$,
$x_2\mapsto ta$, $x_3\mapsto(t^2a)^{-1}z$.

\ref{a3} Every factor of $T_k$ is a conjugate of $\sigma_1$ or of $\Delta_3^2$, the braid
monodromies of a simple tangency and of an ordinary triple point, so by
\cite[Thm.~5.3]{Kulikov-izv04a} (with $N=1$, i.e., in $\CP^2$) there is a Hurwitz
curve $\overline H_k$ with $\operatorname{bmf}(\overline H_k)=T_k$, and
Theorem~\ref{thm:zvk} gives the two fundamental groups, since
$l_{12}=(x_1x_2x_3)^4\mapsto z^4$.  The permutations of the factors have two
orbits, $\{j:j\not\equiv0\}$ and $\{j:j\equiv0\}\pmod3$, of sizes $8$ and $4$;
these are the two components $\overline H'$ and $\overline H''$.  The
permutation of each band is a transposition within one orbit, and the bands
transposing strings $\equiv0\pmod3$ are the two copies of $c_{3,3}$
(conjugated by $b_k$, or shifted) in each of the four copies of $s$, and the
bands $c_{6,3}$ and $c_{6,6}$ in each copy of $\tilde r_6$; so $8+4=12$ of the
$36$ tangencies lie on $\overline H''$, and $24$ on $\overline H'$.  For each
of the $16$ triple points, the strands $h(x_1),h(x_2),h(x_3)$ of the
corresponding factor $h\Delta_3^2h^{-1}$ consist of two strands of
$\overline H'$ and one of $\overline H''$.  Riemann--Hurwitz for the
normalizations gives $\chi(\widetilde H')=2\cdot8-24=-8$ and
$\chi(\widetilde H'')=2\cdot4-12=-4$, i.e., genera $5$ and $3$.  As a check,
the adjunction formula gives $\binom{7}{2}-16=5$ (the $16$ nodes of
$\overline H'$) and $\binom{3}{2}=3$, and the intersection number
$\overline H'\cdot\overline H''=8\cdot4=32$ is accounted for by the two
transverse intersections at each triple point.

Finally, $\BS(2,3)$ is non-Hopfian \cite{Baumslag-Solitar-bams62}, and so are
$\BS(2,3)\times\Z$ and $\BS(2,3)\times\Z/4$ (extend a non-injective
surjective endomorphism by the identity on the second factor); finitely
generated residually finite groups are Hopfian \cite{Lyndon-Schupp-01}.\qed

\begin{remark}
\label{rem:search}
The words $g_1$ and $g_2$, of lengths $4$ and $9$, are the shortest found by a
computer search for bands satisfying the hypotheses of Lemma~\ref{lem:band}
and installing the $\BS(k,k+1)$ relation; any other such band would serve.
\end{remark}

\begin{remark}
\label{rem:mindegree}
Degree $12$ is the smallest degree of a twisted doubling of an iterated
doubling of a full twist with group $\BS(k,k+1)\times\Z$
(Proposition~\ref{prop:min12}), but we do not know the minimal degree of a
Hurwitz curve with group $\BS(2,3)\times\Z$; see Question~\ref{q:mindegree}.  
Kulikova's curve is described in \cite[Cor.~0.8]{Kulikova-izv05} 
and \cite[Thm.~7.5]{Kulikov-usp07} as having only simplest
triple points; the curves $\overline H_k$ have exactly this singularity type,
and $\overline H_2$ has the same C-group $\widetilde G_2$.  The curves differ,
however.  The realization of $\widetilde G_2$ by Kulikov's algorithm goes
through eight doublings of $\Delta_3^2$, so it has degree $768$
(\S\ref{subsec:768}), $4^8$ ordinary triple points,
$3\cdot2^8(2^8-1)$ simple tangencies, and, as in the proof of
Proposition~\ref{prop:samegroup}, projective group $\BS(2,3)\times\Z/2^8$.
The curve $\overline H_2$ has degree $12$, $16$ triple points, $36$
tangencies, and projective group $\BS(2,3)\times\Z/4$.
\end{remark}

\section{Twisted doubling}
\label{sec:twisted}

Kulikov's band lemma conjugates one input of a doubling by a band.  Here we
conjugate up to three inputs by arbitrary braids and compute the resulting
groups, which are especially simple when the strings moved by the braids have
partners; the same argument removes the partner hypothesis from the band
lemma (\S\ref{subsec:twisted}).  Applied to the braids $\gamma_n$ of the
generalized Borromean rings (\S\ref{subsec:Ln}), this proves
Theorem~\ref{thm:twisted} (\S\ref{subsec:proofD}).

\subsection{Twisting by arbitrary braids}
\label{subsec:twisted}

\begin{lemma}
\label{lem:twisted}
Let $s\in\cS_M$ with $\alpha(s)=\Delta_M^2$, and let
$\beta_1,\beta_2,\beta_3\in B_M$.  Then
\[
T=d\bigl(\lambda(\beta_1)(s),\lambda(\beta_2)(s),\lambda(\beta_3)(s),s\bigr)
\]
is a factorization of $\Delta_{2M}^2$ whose factors are conjugates of the
factors of $s$ and of $\sigma_1$, and
\[
G(T)\cong F_M\big/N(s)\,\beta_1\bigl(N(s)\bigr)\,\beta_2\bigl(N(s)\bigr)\,
\beta_3\bigl(N(s)\bigr).
\]
\end{lemma}

\begin{proof}
$\alpha(\lambda(\beta_i)s)=\beta_i\Delta_M^2\beta_i^{-1}=\Delta_M^2$, so
\eqref{eq:double} applies; the group is \eqref{eq:Gdouble} combined with
\eqref{eq:Nconj}.
\end{proof}

For $\beta\in B_M$ let $N(\beta)=\nc{x_j^{-1}\beta(x_j)\colon 1\le j\le M}$.
Since $\beta$ induces the identity of $F_M/N(\beta)$, we have $\beta(w)\in
wN(\beta)$ for all $w$, hence $\beta(N(s))\subseteq N(s)N(\beta)$.  The reverse
inclusion needs partners.

\begin{lemma}
\label{lem:partners}
In Lemma~\ref{lem:twisted}, suppose that for every $\ell$ and every $j$ with
$\beta_\ell(x_j)\neq x_j$ there is an index $p$ with $\beta_\ell(x_p)=x_p$ and
$x_p=x_j$ in $G(s)$.  Then
$G(T)\cong G(s)\big/\nc{x_j^{-1}\beta_\ell(x_j)\colon j\le M,\ \ell\le3}$.
\end{lemma}

\begin{proof}
It remains to show $N(\beta_\ell)\subseteq N(s)\beta_\ell(N(s))$.  If
$\beta_\ell(x_j)\ne x_j$, then
$x_j^{-1}\beta_\ell(x_j)=(x_j^{-1}x_p)\cdot\beta_\ell(x_p^{-1}x_j)$, a product
of an element of $N(s)$ and an element of $\beta_\ell(N(s))$.
\end{proof}

\begin{corollary}
\label{cor:braidtype}
Let $s_0\in\cS_m$ with $\alpha(s_0)=\Delta_m^2$ and $\beta_1,\beta_2,\beta_3\in B_m$,
viewed as braids on the first $m$ of $2m$ strings.  Then
$d\bigl(\lambda(\beta_1)d(s_0),\lambda(\beta_2)d(s_0),\lambda(\beta_3)d(s_0),
d(s_0)\bigr)\in\cS_{4m}$ has group $G(s_0)/\nc{x_j^{-1}\beta_\ell(x_j)}$.  In
particular every C-group
$\la x_1,\dots,x_m\mid l_m\text{ central},\ x_j=\beta_\ell(x_j)\ra$ with at most
three braids $\beta_\ell$ is the group of an explicit Hurwitz curve of degree
$4m$, and one with at most $3k$ braids is the group of one of degree
$2^{k+1}m$.
\end{corollary}

\begin{proof}
In $d(s_0)$ the string $m+j$ is a partner of $x_j$ by \eqref{eq:cmi}, and it
is fixed by $B_m$.  Take $s_0=\Delta_m^2$ for the second statement.  For
more braids, iterate: the output of one step has the partners $x_{m+j}$ and
$x_{2m+j}$ of $x_j$, fixed by the next three braids.
\end{proof}

The same argument, applied to a single band, weakens the hypotheses of
Lemma~\ref{lem:band}: a partner string is not needed.

\begin{lemma}
\label{lem:fixband}
Let $s\in\cS_M$ with $\alpha(s)=\Delta_M^2$, let $b=g\,\sigma_1g^{-1}$ be a band with
endpoint loops $y_1=g(x_1)$, $y_2=g(x_2)$, and suppose that some $y'\in F_M$
with $b(y')=y'$ equals $y_1$ or $y_2$ in $G(s)$.  Then
$T=d(\lambda(b)(s),s,s,s)$ satisfies $G(T)\cong G(s)/\nc{y_1^{-1}y_2}$.  The
fixed subgroup of $b$ contains $y_1y_2$, $g(x_3),\dots,g(x_M)$, and $l_M$.
\end{lemma}

\begin{proof}
As in Lemma~\ref{lem:band}, $b(N(s))\subseteq N(s)\nc{y_1^{-1}y_2}$.  We have
$b(y_1)=y_1y_2y_1^{-1}$ and $b(y_2)=y_1$.  If $y'=y_2$ in $G(s)$, then
$y_2y_1^{-1}=(y_2y'^{-1})\cdot b(y'y_2^{-1})\in N(s)\,b(N(s))$.  If $y'=y_1$ in
$G(s)$, then $(y_1y'^{-1})\cdot b(y'y_1^{-1})=y_1(y_1y_2^{-1})y_1^{-1}$ lies in the
normal subgroup $N(s)\,b(N(s))$, hence so does $y_1y_2^{-1}$.  The last
sentence holds because $\sigma_1$ fixes $x_1x_2,x_3,\dots,x_M$ and
$b(l_M)=l_M$.
\end{proof}

Lemma~\ref{lem:band} is the case $y_1$ or $y_2=x_i$, $y'=x_p$.  If $G(s)$ is
a product $F\times\la l_M\ra$, the hypothesis of Lemma~\ref{lem:fixband}, for
$y'$ in the subgroup generated by $y_1y_2,g(x_3),\dots,g(x_M)$, amounts to a
membership problem in a finitely generated subgroup of the
free group $F$, which Stallings folding decides \cite{Stallings-invent83};
see also \cite{Kapovich-Myasnikov-jalg02}.

\subsection{The generalized Borromean rings}
\label{subsec:Ln}
The generalized Borromean rings of \cite{Suciu-borro} are the closures
$L_n=\hat\gamma_n$, $n\ge1$, of the pure braids
\begin{equation}
\label{eq:gamma}
\gamma_n=\sigma_1^{2n-1}\sigma_2^{-1}\sigma_1\sigma_2^{-(2n-1)}
\sigma_1\sigma_2^{-1}\in B_3
\end{equation}
(as in \cite[(1.1)]{Suciu-borro});
$\gamma_1=(\sigma_1\sigma_2^{-1})^3$ gives the Borromean rings.  Only the 
group of the closed braid is used below, so a difference in orientation conventions,
which would replace $L_n$ by its mirror image, is immaterial.  We write
\[
\Gamma_n=\pi_1(S^3\setminus L_n)=\la x_1,x_2,x_3\mid x_j=\gamma_n(x_j)\ra
\]
for the link group, denoted $G_n$ in \cite{Suciu-borro}.

\subsection{Proof of Theorem~\ref{thm:twisted}}
\label{subsec:proofD}
Let $s=d(\Delta_{m+1}^2)$.  By \eqref{eq:Gdouble},
$G(s)\cong G(\Delta_{m+1}^2)=\la x_1,\dots,x_{m+1}\mid l_{m+1}\text{ central}\ra$,
which is $F(x_1,\dots,x_m)\times\la z\ra$ with $z=l_{m+1}$ and
$x_{m+1}=(x_1\cdots x_m)^{-1}z$, and $x_{m+1+j}=x_j$ in $G(s)$.  The braid
$\beta$ fixes $x_{m+1},\dots,x_{2m+2}$.  Lemma~\ref{lem:partners} with
$p(j)=m+1+j$ gives
\[
G(U_\beta)\cong\bigl(F(x_1,\dots,x_m)\times\la z\ra\bigr)\big/\nc{x_j^{-1}\beta(x_j)}
\cong\la x_1,\dots,x_m\mid x_j=\beta(x_j)\ra\times\Z,
\]
and the first factor is $\pi_1(S^3\setminus\hat\beta)$ by the closed-braid
presentation.  This is \ref{b2}.  For \ref{b1}, $s$ has $4$ factors $\Delta_{m+1}^2$ and
$2(m+1)$ bands, and $U_\beta$ consists of four copies of $s$ (one conjugated by
$\beta$) and $\tilde r_{2m+2}$ twice.  For \ref{b3}, a conjugate of
$\Delta_{m+1}^2$ is the braid monodromy of an ordinary $(m+1)$-fold point and a
conjugate of $\sigma_1$ that of a simple tangency with a fiber.  The permutation of
$\beta c_{m+1,j}\beta^{-1}$ is the transposition of $\pi_\beta(j)$ and $m+1+j$,
all other factors preserve the classes of strings modulo $m+1$, so the
components of $\overline H_\beta$ correspond to the cycles of $\pi_\beta$ and to
the residue class $0$; each residue class contains $4$ strings.  The
component of the class $0$ has $12$ tangencies and one branch through each
$(m+1)$-fold point, so it has genus $3$ by Riemann--Hurwitz and is a smooth
quartic by adjunction.  Finally, $l_{4m+4}=l_{m+1}^4=z^4$ in $G(U_\beta)$, and
Theorem~\ref{thm:zvk} gives the projective group.

For $\beta=\gamma_n$, which is pure, the same count applies to all four
residue classes: $\overline H_{\gamma_n}$ is the union of four smooth
quartics, each with $12$ tangencies, passing once through each of the $16$
quadruple points; since distinct components of a Hurwitz curve can meet only
at points of $Z$, they meet nowhere else, in accordance with B\'ezout.
The last assertion of the theorem, that $\overline H_{\gamma_n}$ is not
isotopic to an algebraic curve, is Corollary~\ref{cor:Ln-nonalg}.  \qed

\section{Alexander invariants and non-algebraicity}
\label{sec:alexander}

In this section we compute the Alexander invariant (\S\ref{subsec:A0}) and the
first characteristic variety (\S\ref{subsec:cvar}) of the groups
$\widetilde G_k$, and prove Theorem~\ref{thm:alex}: the isolated non-torsion
point of $\V_1$ rules out quasi-projectivity.  Along the way we show that the
characteristic variety of every Hurwitz curve group lies in a hypersurface
determined by the degrees of the components
(Proposition~\ref{prop:V1-central}).  In \S\ref{subsec:linkgroups} we treat
the groups $\pi_1(S^3\setminus L_n)\times\Z$, which satisfy the torsion
conditions on $\V_1$ but are nevertheless not quasi-projective.

In \S\S\ref{subsec:A0}--\ref{subsec:cvar}, $\widetilde G_k=G_k\times\la z\ra$
carries the C-structure \eqref{eq:tildeGk}, so $\nu(t)=1$, $\nu(a)=0$,
$\nu(z)=3$, and $q_k=3k^2+3k+1$.

\subsection{The Alexander invariant}
\label{subsec:A0}
For a product with $\Z$, the Alexander invariant is governed by Kulikov's
derived modules $A_n(M)=M/(t^n-1)M$ \cite[\S2.2]{Kulikov-izv08}.

\begin{proposition}
\label{prop:A0-product}
Let $G$ be a C-group with C-generators $x_1,\dots,x_m$, and give
$G\times\la z\ra$ the Hurwitz C-structure of Proposition~\ref{prop:product}.
Then
\[
A_0(G\times\Z)\cong A_0(G)/(t^{m+1}-1)A_0(G)\,\oplus\,\Z,
\]
with $t$ acting trivially on $\Z$.  If $G$ is irreducible, the first summand
is the derived module $A_{m+1}\bigl(A_0(G)\bigr)$, which is isomorphic to
$H_1$ of the $(m+1)$-fold cyclic branched covering of $S^{n+2}$ along any
knotted $n$-manifold with group $G$ \cite[Thm.~0.2]{Kulikov-izv08}.
\end{proposition}

\begin{proof}
Let $N_G$ be the kernel of the total-degree map of $G$, so that
$A_0(G)=N_G/N_G'$.  The added C-generator $x_{m+1}=l_m^{-1}z$ has total
degree $1$, so $\nu(z)=m+1$, and $\ker\nu=N_G\rtimes\la w\ra$ with
$w=x_1^{-(m+1)}z$, which acts on $N_G$ as conjugation by $x_1^{-(m+1)}$.
Hence $H_1(\ker\nu)$ is the module of coinvariants $A_0(G)/(t^{m+1}-1)$,
plus the free summand $\Z w$, which conjugation by $x_1$ fixes.
\end{proof}

\begin{proposition}
\label{prop:A0}
$A_0(\widetilde G_k)\cong\Z\oplus\Lambda/\bigl(q_k,\,(k+1)t-k\bigr)\cong
\Z\oplus\Z/q_k$, where $t$ acts trivially on $\Z$ and as multiplication by
$k/(k+1)$, an automorphism of order $3$, on $\Z/q_k$.  In particular
$\Delta_{\widetilde G_k}(t)\doteq t-1$.
\end{proposition}

\begin{proof}
By Proposition~\ref{prop:A0-product} with $G=G_k$ and $m=2$,
$A_0(\widetilde G_k)\cong A_3\bigl(A_0(G_k)\bigr)\oplus\Z$.  Here
$A_0(G_k)=\Lambda/\bigl((k+1)t-k\bigr)$, which is $\Z[\frac1{k(k+1)}]$ with $t$
acting by $k/(k+1)$ \cite[Ex.~3.20]{Kulikov-izv08}, and its derived module
$A_3$ is cyclic of order $(k+1)^3-k^3=q_k$ \cite[Ex.~2.14]{Kulikov-izv08}:
indeed $(k/(k+1))^3-1=-q_k/(k+1)^3$ and $\gcd(q_k,k(k+1))=1$.  On $\Z/q_k$,
$t^3=1$ and $t-1=-1/(k+1)$ is a unit, so $t$ has order exactly $3$.
\end{proof}

This is consistent with the constraints listed in \ref{k2} of \S\ref{subsec:other},
for a curve of degree $m=3$ with $r=2$ irreducible components: $A_0$ is a
finitely generated abelian group, a quotient of $\Z\oplus\Lambda/(t^3-1)$,
$t$ acts semisimply over $\C$, and the root $t=1$ of $\Delta$ has multiplicity
$1=r-1$.  For $k=1,2$ the torsion is $\Z/7$ and $\Z/19$.

\begin{remark}
\label{rem:branched}
Since $G_k$ is the group of a knotted $2$-sphere in $S^4$
(\S\ref{sec:kulikova}), Proposition~\ref{prop:A0-product} identifies the
torsion $\Z/q_k$ of $A_0(\widetilde G_k)$ with $H_1$ of the $3$-fold cyclic
branched covering of $S^4$ along such a sphere; compare
\cite[\S4]{Kulikov-izv08}.  Likewise, the curves $\overline H_\beta$ of
Theorem~\ref{thm:twisted} carry the C-structure of
Proposition~\ref{prop:product}, so if $\hat\beta$ is a knot $K$, with
$\beta\in B_m$, then $A_0\cong H_1\bigl(\Sigma_{m+1}(K)\bigr)\oplus\Z$, where
$\Sigma_{m+1}(K)$ is the $(m+1)$-fold cyclic branched covering of $S^3$
along $K$.
\end{remark}

\begin{proposition}
\label{prop:covers}
Let $X'_n\to\C^2\setminus H_k$ be the $n$-fold cyclic covering associated
with $\nu$.  Then
\[
H_1(X'_n;\Z)\cong
\begin{cases}\Z^2\oplus\Z/q_k & \text{if } 3\mid n,\\ \Z^2 & \text{if }
3\nmid n.\end{cases}
\]
For every $n$, the resolved $n$-fold cyclic covering of $\CP^2$ branched
along $\overline H_k$ (and possibly along $\Linf$) has first Betti number
zero.
\end{proposition}

\begin{proof}
Milnor's exact sequence for the infinite cyclic cover (see
\cite[\S10.3]{Kulikov-usp07}) gives
$0\to A_0/(t^n-1)A_0\to H_1(X'_n;\Z)\to\Z\to0$, which splits.  By
Proposition~\ref{prop:A0}, $A_0/(t^n-1)A_0\cong\Z\oplus\Z/q_k$ if $3\mid n$, and
$\Z$ otherwise: then $t^n-1$ is $t-1=-1/(k+1)$ or $t^2-1=(t-1)(t+1)$, where
$t+1=(2k+1)/(k+1)$, and both are units on $\Z/q_k$, as
$4q_k=3(2k+1)^2+1$.  The last statement
is \cite[Thm.~10.2]{Kulikov-usp07} ($b_1=r_{n,\neq1}$, the number of $n$-th roots of
unity $\neq1$ among the roots of $\Delta$), as $\Delta\doteq t-1$.
\end{proof}

\subsection{Characteristic varieties}
\label{subsec:cvar}
For a finitely generated group $G$ let $\V_1(G)=\{\rho\in\Hom(G,\C^*)\mid
H^1(G;\C_\rho)\neq0\}$.  The centrality of $l_m$ constrains these loci for
every Hurwitz curve.

\begin{proposition}
\label{prop:V1-central}
Let $\overline H$ be a Hurwitz curve with irreducible components of degrees
$d_1,\dots,d_r$, let $G=\pi_1(\C^2\setminus H)$, and for
$\rho\in\Hom(G,\C^*)$ let $\lambda_c$ be the value of $\rho$ on the
meridians of the $c$-th component.  If $H^q(G;\C_\rho)\neq0$ for some $q$, then
$\lambda_1^{d_1}\cdots\lambda_r^{d_r}=1$.  In particular,
$\V_1(G)\subseteq\{\lambda\in(\C^*)^r\mid\lambda_1^{d_1}\cdots\lambda_r^{d_r}=1\}$.
\end{proposition}

\begin{proof}
By Proposition~\ref{prop:intro-central}, $l_m$ is central of infinite order,
and $\rho(l_m)=\lambda_1^{d_1}\cdots\lambda_r^{d_r}$.  If $\rho(l_m)\neq1$,
then $H^q(\la l_m\ra;\C_\rho)=0$ for all $q$, and the Lyndon--Hochschild--Serre
spectral sequence of the central extension
$1\to\la l_m\ra\to G\to G/\la l_m\ra\to1$ gives $H^*(G;\C_\rho)=0$.
\end{proof}

For an irreducible curve of degree $m$ this recovers the fact that the roots
of its Alexander polynomial are $m$-th roots of unity
\cite[Cor.~10.1]{Kulikov-usp07}; the argument is a multivariable version of
\cite[Lem.~3.7]{Kulikov-izv08}.  For $r\ge2$ the hypersurface
$\prod_c\lambda_c^{d_c}=1$ is a union of $\gcd(d_1,\dots,d_r)$
torsion-translated subtori of codimension one; for instance, no character
with $\lvert\lambda_c\rvert>1$ for all $c$ lies in $\V_1(G)$, and the
characteristic varieties of curves of bounded degree lie in a fixed finite
union of such subtori.

For the groups $\widetilde G_k$ we use the component meridians as
coordinates on $\Hom(\widetilde G_k,\C^*)\cong(\C^*)^2$: $\rho(x_1)=\rho(x_2)=\lambda_1$
(so $\rho(t)=\lambda_1$, $\rho(a)=1$) and $\rho(x_3)=\lambda_2$, so that
$\rho(z)=\lambda_1^2\lambda_2$.

\begin{proposition}
\label{prop:cvar}
$\V_1(\widetilde G_k)=\{(1,1)\}\cup\{(\zeta_k,\zeta_k^{-2})\}$ with
$\zeta_k=k/(k+1)$.
\end{proposition}

\begin{proof}
A crossed homomorphism $\delta$ on $G_k$ with values in $\C_\rho$, where
$\rho(t)=\lambda_1$ and $\rho(a)=1$, satisfies
\[
\delta\bigl(t^{-1}a^kta^{-(k+1)}\bigr)=\bigl(k\lambda_1^{-1}-(k+1)\bigr)\,\delta(a);
\]
hence for $\lambda_1\neq1$, $H^1(G_k;\C_\rho)\neq0$ if and only if
$\lambda_1=\zeta_k$, and $\V_1(G_k)=\{1,\zeta_k\}$.  By the K\"unneth formula for
$K(G_k,1)\times S^1$, $H^1(G_k\times\Z;\C_\rho\otimes\C_\mu)$ vanishes for
$\mu\neq1$ and equals $H^1(G_k;\C_\rho)\oplus H^0(G_k;\C_\rho)$ for $\mu=1$.
So $\V_1(\widetilde G_k)=\{\rho(z)=1,\ \rho|_{G_k}\in\V_1(G_k)\}$, i.e., 
$\lambda_1\in\{1,\zeta_k\}$ and $\lambda_2=\lambda_1^{-2}$.
\end{proof}

With the opposite convention for local coefficients one gets
$(\zeta_k^{-1},\zeta_k^2)$ instead; nothing below depends on this.  The point
$(\zeta_k,\zeta_k^{-2})$ lies off the diagonal $\{\lambda_1=\lambda_2\}$, 
which is the image of the total-degree map $\nu$; this is why the one-variable 
polynomial $\Delta_{\widetilde G_k}=t-1$ does not see it.

\begin{proof}[Proof of Theorem~\ref{thm:alex}]
Parts \ref{c1} and \ref{c2} are Propositions~\ref{prop:A0} and~\ref{prop:cvar}.  
If $G$ is the fundamental group of a smooth quasi-projective variety, the
positive-dimensional irreducible components of $\V_1(G)$ are
torsion-translated subtori, by Arapura's structure theorem
\cite[Thm.~V.1.6]{Arapura-jag97}, whose proof is completed in
\cite[\S\S4--5]{Artal-Cogolludo-Matei-gt13}; and the isolated points of
$\V_1(G)$ are torsion characters \cite[Thm.~1]{Artal-Cogolludo-Matei-gt13}.
More generally, all cohomology jump loci of smooth quasi-projective
varieties are finite unions of torsion-translated subtori
\cite[Thm.~1.1]{Budur-Wang-ens15}.  We use the statement on isolated points.  
The point $(\zeta_k,\zeta_k^{-2})$ is isolated and not torsion, as $\zeta_k\in(0,1)$.  
For $\BS(k,k+1)\times\Z/4$ the same K\"unneth argument over $\C$ gives
$\V_1=\V_1(G_k)\times\{1\}$, which again contains a non-torsion isolated
point.  If $\overline H_k$ were isotopic to an algebraic curve $C$, then
$\CP^2\setminus\overline H_k$ would be diffeomorphic to the smooth
quasi-projective variety $\CP^2\setminus C$, contradicting the above.
\end{proof}

Without a bound on the degree, the characteristic varieties of Hurwitz curve
groups are as flexible as those of finitely presented groups, for which every
subvariety defined over $\Z$ occurs, up to the trivial character
\cite[Lem.~10.3]{Suciu-Yang-Zhao-pisa15},
\cite{Wang-arx13}.

\begin{proposition}
\label{prop:V1-universal}
Let $W\subseteq(\C^*)^n$, $n\ge2$, be a Zariski-closed subset defined over
$\Z$.  There is a Hurwitz curve with $n+1$ irreducible components whose
group $G$ satisfies $\V_1(G)=\iota(W\cup\{1\})$, in the coordinates of
Proposition~\ref{prop:V1-central}, where
$\iota(\lambda)=(\lambda_1,\dots,\lambda_n,(\lambda_1\cdots\lambda_n)^{-1})$.
\end{proposition}

\begin{proof}
Let $W$ be the zero set of $f_1,\dots,f_s\in\Z[t_1^{\pm1},\dots,t_n^{\pm1}]$.
For $i<j$ and $l\le s$, the vector $f_l\cdot\bigl((t_j-1)e_i-(t_i-1)e_j\bigr)$
is orthogonal to $(t_1-1,\dots,t_n-1)$, so by Lyndon's lemma, as in the proof
of \cite[Lem.~10.3]{Suciu-Yang-Zhao-pisa15}, it is the vector of abelianized
Fox derivatives of some $w_{ijl}\in F_n'$.  Let $G_W$ be the C-group with
generators $x_1,\dots,x_n$ and relations $x_k=w_{ijl}^{-1}x_kw_{ijl}$ for all
$k,i,j,l$.  Its abelianization is $\Z^n$, with one component for each
generator.  As $w_{ijl}\in F_n'$, the row of the Alexander matrix given by
the relation $x_k=w^{-1}x_kw$ is $t_k^{-1}(t_k-1)$ times the vector of $w$.
For $\rho\ne1$, choose $k$ with $\rho(t_k)\ne1$; the vectors
$(\rho(t_j)-1)e_i-(\rho(t_i)-1)e_j$ then span a space of dimension $n-1$, so
the Alexander matrix at $\rho$ has rank $n-1$ unless $f_1,\dots,f_s$ all
vanish at $\rho$, in which case it is zero.  Hence $\V_1(G_W)=W\cup\{1\}$.
By Proposition~\ref{prop:product} and \cite[Thm.~7.3]{Kulikov-usp07},
$G_W\times\Z$ is the group of a Hurwitz curve with meridians $x_1,\dots,x_n$
and $l_n^{-1}z$, and $\V_1(G_W\times\Z)=\V_1(G_W)\times\{1\}$ by the
K\"unneth argument in the proof of Proposition~\ref{prop:cvar}; in the
coordinates of Proposition~\ref{prop:V1-central} this is
$\iota(W\cup\{1\})$.
\end{proof}

For instance, $W=\{t_1+t_2+1=0\}$ gives a Hurwitz curve whose characteristic variety
contains a curve that is not a translated subtorus.

\subsection{Link groups}
\label{subsec:linkgroups}
By \cite[Prop.~3.3]{Suciu-borro}, the multivariable Alexander polynomial of the link
$L_n$ of \S\ref{subsec:Ln}, in the variables $t_1,t_2,t_3$ of the meridians
$x_1,x_2,x_3$ of $\Gamma_n$, that is, with the components labeled by the
strings of $\gamma_n$, is
\begin{equation}
\label{eq:DeltaLn}
\Delta_{L_n}(t_1,t_2,t_3)\doteq(t_1-1)(t_2-1)(t_3-1)
\Bigl(\sum_{k=0}^{n-1}(t_1t_2)^k\Bigr)^2 .
\end{equation}
(In \cite{Suciu-borro} the components are labeled instead by the strands of the
conjugate pure braid $\beta_n=[\sigma_1^2,\sigma_2^{2n}]$, which turns $t_1t_2$
into $t_2t_3$.)
With respect to its meridional C-structure, $\Gamma_n$ thus fails, besides
centrality (Proposition~\ref{prop:linking}), the Alexander-invariant
constraint~\ref{k2} of \S\ref{subsec:other}: its one-variable Alexander
polynomial $(t-1)\Delta_{L_n}(t,t,t)\doteq
(t-1)^4\bigl((t^{2n}-1)/(t^2-1)\bigr)^2$ has the root $t=1$ with
multiplicity $4$, whereas a curve with $3$ components would require $2$.

\begin{proposition}
\label{prop:multivar-fails}
For every $n\ge1$, the first characteristic varieties of $\Gamma_n\times\Z$
and of $\Gamma_n\times\Z/q$, $q\ge1$, are finite unions of
torsion-translated subtori, yet none of these groups is the fundamental 
group of a smooth quasi-projective variety.
\end{proposition}

\begin{proof}
For a link with at least two components, $\V_1=V(\Delta_L)\cup\{1\}$
\cite[\S\S6.1, 7.1]{Dimca-Papadima-Suciu-imrn08}; by \eqref{eq:DeltaLn}, $\V_1(\Gamma_n)$
is the union of the subtori $\{t_i=1\}$ and, for $n\ge2$, of the torsion-translated 
subtori $\{t_1t_2=\xi\}$, $\xi^n=1\neq\xi$.  By the K\"unneth argument of
Proposition~\ref{prop:cvar}, $\V_1(\Gamma_n\times\Z)=\V_1(\Gamma_n)\times\{1\}$;
the same argument, using $H^1(\Z/q;\C_\chi)=H^0(\Z/q;\C_\chi)=0$ for
$\chi\ne1$ and $H^1(\Z/q;\C)=0$, gives
$\V_1(\Gamma_n\times\Z/q)=\V_1(\Gamma_n)\times\{1\}$.  In all cases the
components $\{t_2=1\}\times\{1\}$ and $\{t_3=1\}\times\{1\}$ are distinct
and meet along the positive-dimensional subtorus
$\{t_2=t_3=1\}\times\{1\}$.  This violates the relative-position property of
components of $\V_1$ of quasi-projective groups
\cite[Thm.~A]{Dimca-Papadima-Suciu-imrn08}, 
\cite[Thm.~C(2)]{Dimca-Papadima-Suciu-duke09}, exactly as in the 
proof that $\Gamma_n$ itself is not quasi-projective \cite[Prop.~4.7]{Suciu-borro}.
\end{proof}

Thus, unlike for the groups $\widetilde G_k$, the
failure of quasi-projectivity of $\Gamma_n\times\Z$ is invisible to the
torsion conditions of \cite{Arapura-jag97, Artal-Cogolludo-Matei-gt13}, and is
detected instead by the finer relative-position condition.

\begin{corollary}
\label{cor:Ln-nonalg}
The curves $\overline H_{\gamma_n}$ of Theorem~\ref{thm:twisted} are not isotopic
to algebraic curves.
\end{corollary}

\begin{proof}
By Theorem~\ref{thm:twisted}, $\pi_1(\CP^2\setminus\overline H_{\gamma_n})\cong
\Gamma_n\times\Z/4$, which is not quasi-projective by
Proposition~\ref{prop:multivar-fails}; conclude as in the proof of
Theorem~\ref{thm:alex}.
\end{proof}

\section{An algebraic model}
\label{sec:model}

Kulikov's doubling formula has a simple geometric meaning.  The factorization
$d(s)$ is the braid monodromy of the preimage of a curve with braid
monodromy $s$ under the map $\psi(u,v)=(u^2-\alpha_0,v^2-\beta_0)$, and the
twisted doubling $d(\lambda(\gamma)s,s,s,s)$ arises from that preimage by
twisting one quarter of it along a circle (Lemma~\ref{lem:geomdouble}).
Applying this twice to $n$ concurrent lines gives Theorem~\ref{thm:model},
which describes the curves of Theorems~\ref{thm:explicit}
and~\ref{thm:twisted} as twists of explicit unions of quartics in a pencil.

\subsection{Lifting braids through a square root}
\label{subsec:lift}
Let $\beta_0\ge0$ and let $V_1<\dots<V_m$ be real numbers greater than
$-\beta_0$.  Let $B_m(\C\setminus\{-\beta_0\})$ be the braid group of $m$
points in $\C\setminus\{-\beta_0\}$, based at $\{V_1,\dots,V_m\}$.  It contains
the half-twists $\sigma_1,\dots,\sigma_{m-1}$ exchanging consecutive points
along the real axis, which generate a copy of $B_m$, and the braid
$\tau_1$ in which $V_1$ turns once counterclockwise around $-\beta_0$; put
$\tau_j=(\sigma_{j-1}\cdots\sigma_1)\,\tau_1\,(\sigma_{j-1}\cdots\sigma_1)^{-1}$.
Viewing $-\beta_0$ as a fixed additional string, the kernel $F$ of the map
forgetting that string is a free group of rank $m$, with basis
$\tau_1,\dots,\tau_m$: it is the point-pushing subgroup of the Birman exact
sequence \cite[Thm.~1.4]{Birman-74}.  Thus
$B_m(\C\setminus\{-\beta_0\})=F\rtimes B_m$, and the product
$\tau_1\cdots\tau_m$, in which all the points turn once around $-\beta_0$,
commutes with $B_m$, since $\Delta_m^2\tau_1\cdots\tau_m$ is the full
twist of all $m+1$ strings.  For $1\le i<m$,
\begin{equation}
\label{eq:tauconj}
\sigma_i\tau_i\sigma_i^{-1}=\tau_{i+1},\qquad
\sigma_i\tau_{i+1}\sigma_i^{-1}=\tau_{i+1}^{-1}\tau_i\tau_{i+1},\qquad
\sigma_i\tau_j\sigma_i^{-1}=\tau_j\quad(j\neq i,i+1).
\end{equation}
The first identity is the definition of $\tau_{i+1}$.  For the third, note
that $\sigma_i$ commutes with $\tau_1$ if $i\ge2$, their supports being
disjoint, hence with $\tau_j$ for $j<i$; for $j>i+1$, the braid relations
give $\sigma_i(\sigma_{j-1}\cdots\sigma_1)=(\sigma_{j-1}\cdots\sigma_1)
\sigma_{i+1}$, and $\sigma_{i+1}$ commutes with $\tau_1$.  The second
identity follows from the other two, since $\sigma_i$ commutes with
$\tau_1\cdots\tau_m$.  Thus $\lambda(\sigma_i)$ acts on the sequence
$(\tau_1,\dots,\tau_m)$ as the Hurwitz move at the $i$-th place.  Since
simultaneous conjugation commutes with Hurwitz moves, induction on the
length of $h$ gives
\begin{equation}
\label{eq:tauhurwitz}
\lambda(h)(\tau_1,\dots,\tau_m)\ \text{is Hurwitz equivalent to}\
(\tau_1,\dots,\tau_m)\quad\text{for every }h\in B_m .
\end{equation}
By \eqref{eq:tauconj}, conjugation by $\sigma_i$ also permutes the
conjugacy classes of $\tau_1,\dots,\tau_m$ in $F$; as $F$ is normal, every
conjugate of $\tau_1$ in $F\rtimes B_m$ is conjugate in $F$ to some
$\tau_k$.

The preimages of $V_1,\dots,V_m$ under $q(v)=v^2-\beta_0$ are the $2m$ real
numbers $\pm\sqrt{V_j+\beta_0}$; ordering them, we identify them with the
standard configuration of $B_{2m}$.  Lifting motions through $q$ defines a
homomorphism $\Phi\colon B_m(\C\setminus\{-\beta_0\})\to B_{2m}$.

\begin{lemma}
\label{lem:lift}
$\Phi(b)=\Delta_{m,0}\,b\,\Delta_{m,0}^{-1}\cdot\sh(b)$ for $b\in B_m$, and
$\Phi(\tau_j)=\bar c_{m,m+1-j}$ for $1\le j\le m$.  In particular, $\Phi$
maps the sequence $(\tau_1,\dots,\tau_m)$ termwise to the factorization
$\bar r_m=\bar c_{m,m}\cdot\bar c_{m,m-1}\cdots\bar c_{m,1}$.
\end{lemma}

\begin{proof}
The map $V\mapsto\sqrt{V+\beta_0}$ takes the half-plane $\Re(V+\beta_0)>0$
orientation-preservingly into the right half-plane and preserves the order
of real points, while $V\mapsto-\sqrt{V+\beta_0}$ takes it into the left
half-plane and reverses that order.  Hence
$\Phi(\sigma_j)=\sigma_{m-j}\sigma_{m+j}=\Delta_{m,0}\sigma_j\Delta_{m,0}^{-1}\cdot
\sh(\sigma_j)$, which gives the first assertion.  When $V_1+\beta_0$ turns
once around $0$, its two square roots, which are the strings $m$ and $m+1$,
are exchanged by a counterclockwise half-turn, so $\Phi(\tau_1)=\sigma_m=\bar
c_{m,m}$.  As $\Phi(\sigma_j)=\sigma_{m,m-j}$, the formula for $\Phi(\tau_j)$
follows from the definition of $\bar c_{m,i}$ in \S\ref{subsec:doubling}.
\end{proof}

In the braid $\Delta_m^2\tau_1\cdots\tau_m$ all the points $V_j$ turn once
around $-\beta_0$, together; its lift turns the $2m$ strings by a half-turn.
So Lemma~\ref{lem:lift} gives a geometric proof of Kulikov's identity
\cite[Lem.~2.1]{Kulikov-izv04a}, in which $\bar r_m$ stands for its product,
\begin{equation}
\label{eq:halftwist}
\Delta_{2m}=\Delta_{m,0}^2\,\Delta_{m,m}^2\,\bar r_m .
\end{equation}
Note also that conjugation by $\Delta_{2m}$, which maps $\sigma_i$ to $\sigma_{2m-i}$,
fixes every $\sigma_{m,j}$ and $\sigma_m$, hence every $\bar c_{m,i}$.

\begin{lemma}
\label{lem:artin}
Let $F$ be a free group with basis $y_1,\dots,y_m$, and let
$g_1,\dots,g_m\in F$ be elements, each conjugate to one of
$y_1,\dots,y_m$, such that $g_1\cdots g_m=y_1\cdots y_m$.  Then
$(g_1,\dots,g_m)$ is Hurwitz equivalent to $(y_1,\dots,y_m)$.
\end{lemma}

\begin{proof}
This is Artin's characterization of braid automorphisms of a free group,
\cite{Artin-ann47}, \cite[Thm.~1.9]{Birman-74}, which is stated in
\cite{Birman-74} for endomorphisms; we translate it.  Abelianizing $F$, the
basis elements to which $g_1,\dots,g_m$ are conjugate are $y_1,\dots,y_m$
in some order.  Writing endomorphisms on the right, as in \cite{Birman-74},
let $\varphi$ be the endomorphism of $F$ with $(y_j)\varphi=g_j$.  It satisfies
the two conditions of \cite[Thm.~1.9]{Birman-74}, so it is a product of the
automorphisms $\sigma_\nu^{\pm1}$ of \cite[(1-14)]{Birman-74}, which are
given by the formulas \eqref{eq:action} acting on the right, a product
$\sigma_\nu\varphi'$ meaning that $\sigma_\nu$ is applied first.  Now
$(y_\nu)\sigma_\nu\varphi'=g'_\nu g'_{\nu+1}(g'_\nu)^{-1}$ and
$(y_{\nu+1})\sigma_\nu\varphi'=g'_\nu$, where $g'_j=(y_j)\varphi'$, while
$\sigma_\nu^{-1}$ gives $g'_{\nu+1}$ and $(g'_{\nu+1})^{-1}g'_\nu g'_{\nu+1}$;
in both cases the other entries are unchanged.  So passing from $\varphi'$ to
$\sigma_\nu^{\pm1}\varphi'$ applies a Hurwitz move to the sequence of images
of $y_1,\dots,y_m$, and induction on the number of factors of $\varphi$ shows
that $(g_1,\dots,g_m)$ is Hurwitz equivalent to $(y_1,\dots,y_m)$.
\end{proof}
We shall apply Lemma~\ref{lem:artin} to the point-pushing subgroup $F$
itself, with the basis $\tau_1,\dots,\tau_m$ and the multiplication of
$B_m(\C\setminus\{-\beta_0\})$; no identification of $F$ with the
fundamental group of a punctured plane, which would involve orientation
and order conventions, is needed.

\begin{lemma}
\label{lem:framing}
For $h\in B_m$ and $t_1,\dots,t_4\in\cS_m$ with $\alpha(t_i)=\Delta_m^2$,
\[
d\bigl(\lambda(h)t_1,\dots,\lambda(h)t_4\bigr)=
\lambda\bigl(h\,\sh(h)\bigr)\,d(t_1,\dots,t_4)\quad\text{in }\cS_{2m}.
\]
\end{lemma}

\begin{proof}
Since $\sh(h)$ commutes with the factors of $t_1,t_2$, and $h$ with those of
$\sh(t_3),\sh(t_4)$, the two factorizations differ only in their last
$2m$ factors, and it suffices to show that $\lambda(h\,\sh(h))\tilde r_m$
is Hurwitz equivalent to $\tilde r_m$.  Put $h'=\Delta_mh\Delta_m^{-1}$.  Then
$\Delta_{m,m}\,h\,\sh(h)\,\Delta_{m,m}^{-1}=h\,\sh(h')=\Phi(h')$, and since
$\tilde r_m=\lambda(\Delta_{m,m}^{-1})\bar r_m$, Lemma~\ref{lem:lift} gives
\[
\lambda\bigl(h\,\sh(h)\bigr)\tilde r_m
=\lambda(\Delta_{m,m}^{-1})\,\Phi\bigl(\lambda(h')(\tau_1,\dots,\tau_m)\bigr).
\]
As $h'\in B_m$, \eqref{eq:tauhurwitz} shows that
$\lambda(h')(\tau_1,\dots,\tau_m)$ is Hurwitz equivalent to
$(\tau_1,\dots,\tau_m)$, whose image under the homomorphism $\Phi$ is
$\bar r_m$.
\end{proof}

\subsection{Doubling by pullback}
\label{subsec:pullback}
\begin{figure}[ht]
\centering
\begin{tikzpicture}[tk/.style={line width=0.5pt},bp/.style={line width=0.5pt,fill=white},
  bd/.style={line width=0.5pt},lp/.style={line width=0.4pt}]
\begin{scope}[x=2.5cm,y=2.5cm]
\draw[black!50] (-0.5,0) -- (1,0);
\filldraw[fill=black!12,bd] (0.68,0) circle (0.32);
\draw[bd] (0.25,0) circle (0.75);
\draw[lp] (1,0) arc[start angle=0,end angle=-72,radius=0.535]
  .. controls (0.345,-0.602) and (0,-0.3) .. (0,-0.08);
\draw[lp,->] (0,-0.08) arc[start angle=-90,end angle=100,radius=0.08];
\draw[lp] (100:0.08) arc[start angle=100,end angle=270,radius=0.08];
\foreach \w in {0.14,0.215,0.29} \draw[tk] (\w,-0.035) -- ++(0,0.07);
\foreach \x/\y in {0.671/0.105,0.585/-0.044,0.707/-0.101} \fill (\x,\y) circle (1.5pt);
\draw[bd] (-0.025,-0.025) -- (0.025,0.025) (-0.025,0.025) -- (0.025,-0.025);
\draw[bp] (1,0) ++(-1.5pt,-1.5pt) rectangle ++(3pt,3pt);
\node[above right] at (0.98,0) {$U_b$};
\node at (0.68,0.215) {$\Omega_0$};
\node at (-0.47,0.68) {$\Omega_1$};
\node[below left] at (-0.05,-0.05) {$-\alpha_0$};
\node at (0.4,-0.36) {$\delta$};
\end{scope}
\node at (0.625,-2.3) {(a)};
\draw[->] (5.2,0) -- (3.4,0) node[midway,above] {$u\mapsto u^2-\alpha_0$};
\def\sqdisk#1{plot[domain=0:355,samples=72,smooth cycle,variable=\t]
  ({atan2(#1*sin(\t),0.68+#1*cos(\t))/2}:{sqrt(veclen(0.68+#1*cos(\t),#1*sin(\t)))})}
\begin{scope}[shift={(8.95,0)},x=2.75cm,y=2.75cm]
\draw[black!50] (-1,0) -- (1,0);
\foreach \r in {0,180}{
  \begin{scope}[rotate=\r]
  \fill[black!12] \sqdisk{0.176};
  \fill[black!28,even odd rule] \sqdisk{0.32} \sqdisk{0.176};
  \draw[bd] \sqdisk{0.32};
  \draw[line width=0.3pt] \sqdisk{0.176};
  \foreach \x/\y in {0.822/0.064,0.765/-0.029,0.843/-0.06} \fill (\x,\y) circle (1.5pt);
  \foreach \w in {0.14,0.215,0.29} \draw[tk] ({sqrt(\w)},-0.035) -- ++(0,0.07);
  \end{scope}}
\draw[bd] plot[domain=0:358,samples=180,smooth cycle,variable=\p]
  ({\p}:{sqrt(0.25*cos(2*\p)+sqrt(0.5625-(0.25*sin(2*\p))^2))});
\draw[lp,->] plot[smooth] coordinates {
  (1.0000,0.0000) (0.9999,-0.0085) (0.9996,-0.0170) (0.9991,-0.0255) (0.9975,-0.0425) (0.9951,-0.0594)
  (0.9911,-0.0803) (0.9845,-0.1052) (0.9763,-0.1297) (0.9644,-0.1576) (0.9480,-0.1885) (0.928,-0.218)
  (0.903,-0.250) (0.874,-0.279) (0.838,-0.307) (0.799,-0.331) (0.752,-0.350) (0.700,-0.362)
  (0.649,-0.367) (0.593,-0.365) (0.536,-0.356) (0.479,-0.340) (0.422,-0.316) (0.370,-0.285)
  (0.324,-0.241) (0.291,-0.185) (0.272,-0.123) (0.249,-0.061) (0.214,-0.005) (0.169,0.046)
  (0.117,0.089) (0.061,0.126) (0.001,0.160) (-0.059,0.190) (-0.117,0.217) (-0.177,0.244)
  (-0.239,0.269)};
\draw[lp] plot[smooth] coordinates {
  (-0.239,0.269) (-0.300,0.293) (-0.358,0.313) (-0.419,0.332) (-0.478,0.347) (-0.534,0.358)
  (-0.592,0.365) (-0.647,0.367) (-0.701,0.362) (-0.751,0.350) (-0.798,0.331) (-0.838,0.307)
  (-0.874,0.279) (-0.903,0.250) (-0.928,0.218) (-0.9479,0.1885) (-0.9644,0.1576) (-0.9763,0.1297)
  (-0.9845,0.1052) (-0.9911,0.0803) (-0.9951,0.0594) (-0.9975,0.0425) (-0.9991,0.0255) (-0.9996,0.0170)
  (-0.9999,0.0085) (-1.0000,0.0000)};
\draw[line width=1.3pt] (0.86,0) -- (1,0);
\draw[bd] (-0.025,-0.025) -- (0.025,0.025) (-0.025,0.025) -- (0.025,-0.025);
\draw[bp] (1,0) ++(-1.5pt,-1.5pt) rectangle ++(3pt,3pt);
\draw[bp] (-1,0) ++(-1.5pt,-1.5pt) rectangle ++(3pt,3pt);
\node[above right] at (0.98,0) {$u_b$};
\node[below left] at (-0.98,0) {$-u_b$};
\node at (0.66,0.3) {$\Omega$};
\node at (-0.8,-0.3) {$-\Omega$};
\draw[line width=0.3pt] (1.06,0.28) -- (0.905,0.118);
\node at (1.1,0.33) {$A$};
\node[below] at (0,-0.02) {$0$};
\node at (-0.84,0.62) {$\widetilde\Omega_1$};
\node at (0.36,-0.43) {$\tilde\delta$};
\end{scope}
\node at (8.95,-2.3) {(b)};
\begin{scope}[shift={(0,-4.2)}]
\fill[black!12] (1.4,-1.2) rectangle (4.2,1.2);
\fill[black!28] (4.2,-1.2) rectangle (7.6,1.2);
\draw[line width=0.3pt] (4.2,-1.2) -- (4.2,1.2);
\draw[bd] (7.6,-1.2) -- (7.6,1.2);
\draw[black!60,dashed] (1.4,0) -- (7.6,0);
\foreach \y in {0.3,0.6,0.9}{
  \draw[bd] (1.4,\y) -- (5.2,\y) (6.6,\y) -- (7.6,\y);
  \draw[bd] (1.4,-\y) -- (7.6,-\y);}
\draw[bd,fill=white] (5.2,0.18) rectangle (6.6,1.02);
\node at (5.9,0.6) {$\gamma$};
\draw[line width=1.3pt] (1.4,-1.2) -- (7.6,-1.2);
\draw[bp] (7.6,-1.2) ++(-1.5pt,-1.5pt) rectangle ++(3pt,3pt);
\node[left] at (1.4,0.6) {$D^+$};
\node[left] at (1.4,-0.6) {$D^-$};
\node[right] at (7.7,-0.08) {$\Re v=0$};
\draw[line width=0.4pt] (7.75,0.13) -- ++(0.1,0) -- ++(0,0.94) -- ++(-0.1,0);
\node[right] at (7.85,0.6) {$K$};
\node[below] at (2.8,-1.25) {$\Omega\setminus A$};
\node[below] at (5.9,-1.25) {$A$};
\node[below] at (7.6,-1.25) {$u_b$};
\end{scope}
\node at (4.5,-6.1) {(c)};
\end{tikzpicture}
\caption{Doubling by pullback, schematically, for $m=3$.
(a) The critical values of $\pr|_C$ (dots) lie in $\operatorname{int}\Omega_0$;
the $U$-coordinates of the points of $C\cap\{V=-\beta_0\}$ (ticks) and
$-\alpha_0$ lie in $\operatorname{int}\Omega_1\setminus\Omega_0$.
(b) Their preimages under $u\mapsto u^2-\alpha_0$, with a collar $A$ of
$\partial\Omega$ shaded darker.  The ticks are now critical values where $D$ is
tangent to the pencil at $v=0$, the point $0$ is not a critical value, and
$\delta$ lifts to the path $\tilde\delta$ from $u_b$ to $-u_b$.  Over
$\Omega\cup-\Omega$, $D$ consists of four copies of $C\cap\pr^{-1}(\Omega_0)$,
one for each choice of signs of $\Re u$ and $\Re v$; they carry the four inputs
of the doubling formula, and the tangencies at $v=0$ its two copies of
$\tilde r_m$.
(c) The twisted curve $D^\gamma$, for a braid $\gamma\in B_m$, over the thick
segment in (b), with $u$ horizontal, $\Re v$ vertical, and the shading of (b).
It agrees with $D$ except over $A$, where diffeomorphisms $h_u$, supported in a
disk $K\subset\{\Re v>0\}$ containing the fibers of $D^+$, move the points of
$D^+$ so as to realize $\gamma$ (the box); the points of $D^-=-D^+$ stay put.
This conjugates by $\gamma$ the input carried by $D^+$
(Lemmas~\ref{lem:twist} and~\ref{lem:geomdouble}).}
\label{fig:pullback}
\end{figure}
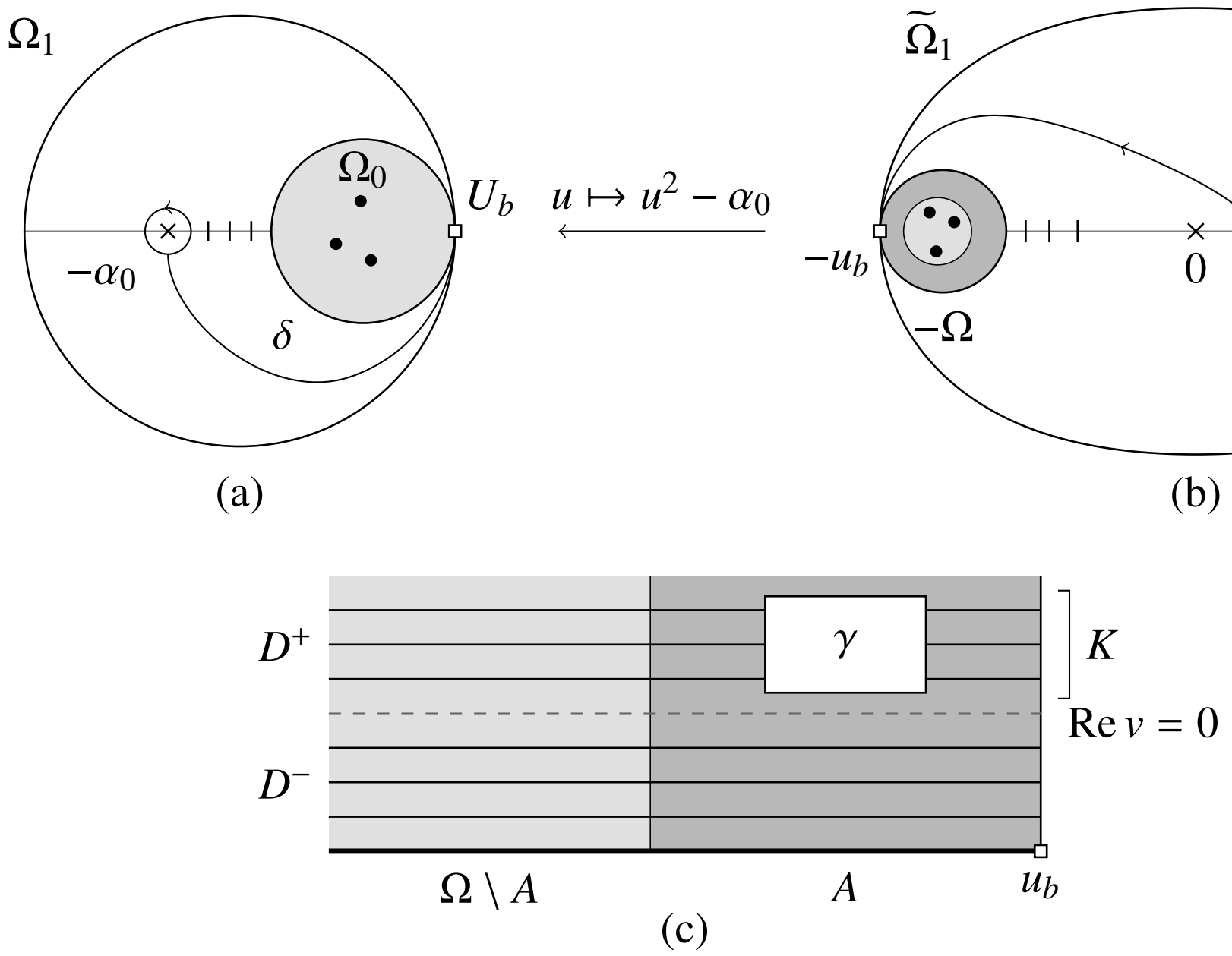

Let $C\subset\C^2$ be an algebraic curve, with coordinates $(U,V)$, of degree $m$ with
respect to the projection $\pr(U,V)=U$ and transverse to the line at infinity.  
Fix $\alpha_0,\beta_0\ge0$, a closed round disk $\Omega_0\subset\C$ 
and a real point $U_b\in\partial\Omega_0$ such that:
\begin{enumerate}[label=(\alph*),itemsep=1pt]
\item \label{pb1}
the interior of $\Omega_0$ contains all critical values of $\pr|_C$, and
$\Re(U+\alpha_0)>0$ on $\Omega_0$;
\item \label{pb2}
$\Re(V+\beta_0)>0$ at every point $(U,V)\in C$ with $U\in\Omega_0$, and the
$m$ points of $C$ over $U_b$ are real;
\item \label{pb3}
$C$ meets the line $\{V=-\beta_0\}$ transversally, in $m$ points with
distinct $U$-coordinates, none of which equals $-\alpha_0$.
\end{enumerate}
Let $s\in\cS_m$ be the braid monodromy factorization of $C$ over $\Omega_0$,
with respect to the base point $U_b$, the fiber over $U_b$ being identified
with the standard configuration by the real order.  Let
$D=\psi^{-1}(C)$, where $\psi(u,v)=(u^2-\alpha_0,v^2-\beta_0)$; let $\Omega$ be
the component of $\{u\mid u^2-\alpha_0\in\Omega_0\}$ in the half-plane
$\Re u>0$, and put $u_b=\sqrt{U_b+\alpha_0}\in\partial\Omega$.  By~\ref{pb1},
$u\mapsto u^2-\alpha_0$ maps $\Omega$ diffeomorphically onto $\Omega_0$.
By~\ref{pb2}, over $\Omega$ the curve $D$ is the disjoint union of its parts
$D^+$ and $D^-$ in the half-spaces $\Re v>0$ and $\Re v<0$, and $\psi$
maps $D^+$ diffeomorphically onto $C\cap\pr^{-1}(\Omega_0)$; so the boundary
of $D^+$ over $\partial\Omega$ is the closed braid of the full twist.  The
fiber $D_{u_b}=\{v\mid(u_b,v)\in D\}$ consists of $2m$ real numbers,
the $m$ positive ones forming the fiber of $D^+$.  The construction is
sketched in Figure~\ref{fig:pullback}.

\smallskip\noindent\emph{Twisting.}
Choose a closed collar $A\subset\Omega$ of $\partial\Omega$ containing no
critical value of $\pr|_D$; this is possible, since by~\ref{pb1}
and~\ref{pb2} the critical values of $\pr|_D$ in $\Omega$ are the preimages
of those of $\pr|_C$, which lie in the interior of $\Omega_0$.  Choose also
a closed round disk $K$ in the half-plane $\{\Re v>0\}$ with
$D^+\subset\Omega\times\operatorname{int}K$, and for $u\in\Omega$ let
$F_u\subset\operatorname{int}K$ be the fiber of $D^+$ over $u$.  We
identify $B_{2m}$ with the mapping class group of $(\C,D_{u_b})$ as in
\S\ref{subsec:hurwitz}, using the real order of $D_{u_b}$.  
Let $\gamma\in B_m$.  A \emph{twisting family}\/ for $\gamma$ is a 
smooth family $(h_u)_{u\in\Omega}$ of diffeomorphisms of $\C$
supported in $\operatorname{int}K$, such that $h_u=\id$ for $u$ near
the closure of $\Omega\setminus A$, $h_u(F_u)=F_u$ for $u$ near
$\partial\Omega$, and the class of $h_{u_b}$ is $\sh(\gamma)$; that is,
$h_{u_b}$ realizes $\gamma$ on the $m$ points of $F_{u_b}$, in their real
order.  Put $\Theta(u,v)=(u,h_u(v))$, a fiber-preserving diffeomorphism of
$\pr^{-1}(\Omega)$, and
\[
D^\gamma=\bigl(D\setminus\pr^{-1}(\Omega)\bigr)\cup
\Theta\bigl(D\cap\pr^{-1}(\Omega)\bigr).
\]
Since $\Theta$ preserves $D$ over a neighborhood of $\partial\Omega$ and is
the identity over $\Omega\setminus A$, $D^\gamma$ is a Hurwitz curve which
coincides with $D$ outside $\pr^{-1}(A)$, in particular near all critical
points of $\pr|_D$.  Informally, $D^\gamma$ is obtained by cutting $D^+$
along $\pr^{-1}(\partial\Omega)$ and regluing it by the diffeomorphisms
$h_u$, $u\in\partial\Omega$, which realize $\gamma$
(Figure~\ref{fig:pullback}(c)); it is the analogue, for
braided surfaces, of a twisted fiber sum.  When $\gamma$ is a half-twist
along an arc, $D^\gamma$ is the braiding of $D$ along the annulus swept out
by this arc over $\partial\Omega$, in the sense of
\cite[\S3]{Auroux-Donaldson-Katzarkov-mathann03}; in general, writing $\gamma$
as a product of half-twists exhibits $D^\gamma$ as an iterated braiding
\cite[Rem.~3.2]{Auroux-Donaldson-Katzarkov-mathann03}.

Twisting families exist for every $\gamma$, and $D^\gamma$ does not 
depend on the choices of $A$, $K$ and $(h_u)$, up to an H-isotopy
\cite[\S5.2]{Kulikov-usp07} that is the identity near the critical points.
Indeed, a twisting family for $(A,K)$ is one for $(A',K')$ whenever
$A\subset A'$ and $K\subset K'$, any two collars contain a common one, and
any two disks lie in a common one; and a path of twisting families gives a
smooth family of Hurwitz curves $D^\gamma$, all equal to $D$ outside
$\pr^{-1}(A)$ and near $\pr^{-1}(\partial\Omega)$, hence, by fiberwise isotopy
extension, an H-isotopy.  So it suffices to show that, for fixed $A$ and
$K$, the twisting families for $\gamma$ form a non-empty connected space.
The diffeomorphisms of $\C$ supported in $\operatorname{int}K$ form a
contractible group \cite{Smale-pams59}; since configuration spaces of points
in a disk are aspherical \cite{Fadell-Neuwirth-mathscand62}, the exact sequence of
the fibration given by evaluation on a set of $m$ points shows that the
diffeomorphisms preserving that set form a space whose components are
contractible and are indexed by $B_m$.  Transporting these components 
once along $\partial\Omega$ conjugates their indices by the full twist, which is
central.  Hence the families $(h_u)$ over a neighborhood of $\partial\Omega$
with $h_u(F_u)=F_u$ and $h_{u_b}$ in the class of $\sh(\gamma)$ form a
non-empty connected space, and each of them extends over the rest of $A$,
uniquely up to homotopy, to a family that is the identity near the inner
boundary of $A$.

\begin{lemma}
\label{lem:twist}
Let $\ell$ be a loop based at $u_b$ that avoids the critical values of
$\pr|_D$, and let $b$ and $b^\gamma\in B_{2m}$ be the braid monodromies of
$D$ and $D^\gamma$ along $\ell$, with respect to the real order of
$D_{u_b}$.  If $\ell\subset\C\setminus\operatorname{int}\Omega$, then
$b^\gamma=b$; if $\ell\subset\Omega$, then
$b^\gamma=\sh(\gamma)\,b\,\sh(\gamma)^{-1}$.  Consequently, if
$\operatorname{bmf}(D)=t\cdot t'$ with respect to a geometric basis of loops
based at $u_b$ such that the loops giving $t$ lie in $\Omega$ and those
giving $t'$ lie in $\C\setminus\operatorname{int}\Omega$, then
$\operatorname{bmf}(D^\gamma)=\lambda(\sh\gamma)(t)\cdot t'$ with respect to
the same basis.
\end{lemma}

\begin{proof}
The curves $D$ and $D^\gamma$ coincide over $\C\setminus\operatorname{int}
\Omega$, and in particular $D^\gamma_{u_b}=D_{u_b}$; this gives the first
assertion.  If $\ell\subset\Omega$, then $\Theta$ is an isomorphism over
$\ell$ from the pair $(\pr^{-1}(\ell),D\cap\pr^{-1}(\ell))$ to the pair
$(\pr^{-1}(\ell),D^\gamma\cap\pr^{-1}(\ell))$, equal to $h=h_{u_b}$ over
$u_b$.  Hence, if $T$ is a monodromy diffeomorphism of $(\C,D_{u_b})$ along
$\ell$ for $D$, then $h\circ T\circ h^{-1}$ is one for $D^\gamma$.  The braid
$b$ corresponds to the class of $T$ (or of $T^{-1}$, according to the
orientation convention) under the identification above, which is an
isomorphism and maps the class of $h$ to $\sh(\gamma)$; in either case
$b^\gamma=\sh(\gamma)\,b\,\sh(\gamma)^{-1}$.  If $\ell=l\cdot c\cdot l^{-1}$,
with $c$ a small circle around a critical value, then $c\subset\Omega\setminus A$,
over which $\Theta$ is the identity; so the same argument, applied to the
transports along $l$ of the motions of the pieces of $D$ near the critical
points over $c$, shows that the local factors (\S\ref{subsec:hurwitz}) are
also conjugated by $\sh(\gamma)$.  The last assertion follows.
\end{proof}

\begin{lemma}
\label{lem:geomdouble}
Let $\Omega_1\supset\Omega_0$ be a closed disk with $\Omega_0\cap\partial
\Omega_1=\{U_b\}$ whose interior contains $-\alpha_0$ and the $U$-coordinates
of the points of $C\cap\{V=-\beta_0\}$, for instance a large round disk tangent
to $\Omega_0$ at $U_b$ when these points and $\Omega_0$ lie on the same side of
the tangent line to $\Omega_0$ at $U_b$, and let $\widetilde\Omega_1$ be its preimage 
under $u\mapsto u^2-\alpha_0$, a closed disk with $u_b\in\partial\widetilde\Omega_1$
that contains all critical values of $\pr|_D$.  With respect to the base
point $u_b$, the boundary loop of $\widetilde\Omega_1$ and the real order of
the fiber over $u_b$,
\[
\operatorname{bmf}(D^\gamma)=
\lambda(\Delta_{m,0}\Delta_{2m}^{-1})\,
d\bigl(\lambda(\Delta_m\gamma)s,\lambda(\Delta_m)s,\lambda(\Delta_m)s,
\lambda(\Delta_m)s\bigr)\quad\text{in }\cS_{2m}.
\]
Consequently
$\operatorname{bmf}(D^\gamma)$ has the type of $d(\lambda(\gamma)s,s,s,s)$, and
$\operatorname{bmf}(D)$ has the type of $d(s)$.
\end{lemma}

\begin{proof}
Put $q(u)=u^2-\alpha_0$.  The critical values of $\pr|_D$ are the preimages
under $q$ of the critical values of $\pr|_C$, at which $\psi$ is a local
diffeomorphism by~\ref{pb1} and~\ref{pb2}, and of the $U$-coordinates of the
points of $C\cap\{V=-\beta_0\}$, which give simple tangencies at points with
$v=0$, one in each fiber by~\ref{pb3}.  Choose loops based at $U_b$ as
follows: first loops $\ell_1,\dots,\ell_p$ in $\Omega_0$ around the critical
values of $\pr|_C$, forming a geometric basis of $\Omega_0$ minus these
points, so that their local factors form $s=(b_1,\dots,b_N)$; then loops
$\nu_1,\dots,\nu_m$ around the $U$-coordinates of the points of
$C\cap\{V=-\beta_0\}$; and finally a loop $\delta$ around $-\alpha_0$; the last
$m+1$ loops lying in $\Omega_1\setminus\operatorname{int}\Omega_0$, and all
$p+m+1$ loops forming a geometric basis of $\Omega_1$ minus these points.
Such loops exist.  Indeed, as both configurations consist of two closed
disks, one inside the other, whose boundaries meet in one point, an
orientation-preserving homeomorphism of $\C$ carries $\Omega_1$, $\Omega_0$
and $U_b$ to the closed upper half of the disk $\lvert U\rvert\le2$, the
sector $\{\lvert U\rvert\le1,\ \pi/3\le\arg U\le2\pi/3\}\cup\{0\}$ and $0$; by
the homogeneity of the open sets $\operatorname{int}\Omega_0$ and
$\operatorname{int}\Omega_1\setminus\Omega_0$, it may be chosen to carry the
critical values of $\pr|_C$ into the interior of the sector and the other
$m+1$ points into the region $\{\lvert U\rvert<2,\ 2\pi/3<\arg U<\pi\}$,
all $p+m+1$ points having distinct arguments and $-\alpha_0$ going to the
point of largest argument; there,
the loops running along the segments from $0$ to these points, taken in the
order of increasing argument, have the required properties.

Each of these loops other than $\delta$ lifts through $q$, from $u_b$ and
from $-u_b$, to loops $\tilde\ell_i,\tilde\nu_j$ and
$\tilde\ell_i'=-\tilde\ell_i$, $\tilde\nu_j'=-\tilde\nu_j$ around the two
preimages of the corresponding point; the $\tilde\ell_i$ lie in $\Omega$, the
$\tilde\ell_i'$ in $-\Omega$, and the $\tilde\nu_j$, $\tilde\nu_j'$ in
$\C\setminus\operatorname{int}(\Omega\cup-\Omega)$.  The loop $\delta$ lifts
to a path $\tilde\delta$ from $u_b$ to $-u_b$, also in
$\C\setminus\operatorname{int}(\Omega\cup-\Omega)$ (Figure~\ref{fig:pullback}(b)).
The $2(p+m)$ loops
\[
\tilde\ell_1,\dots,\tilde\ell_p,\ \tilde\nu_1,\dots,\tilde\nu_m,\
\tilde\delta\tilde\ell_1'\tilde\delta^{-1},\dots,\tilde\delta\tilde\ell_p'\tilde\delta^{-1},\
\tilde\delta\tilde\nu_1'\tilde\delta^{-1},\dots,\tilde\delta\tilde\nu_m'\tilde\delta^{-1}
\]
are positive meridians of the critical values of $\pr|_D$, one around
each.  Their product is the boundary loop of $\widetilde\Omega_1$: it is the
lift from $u_b$ of the square of the boundary loop of $\Omega_1$, followed by
the inverse of the lift $\tilde\delta\cdot(-\tilde\delta)$ of $\delta^2$, which
is a loop around $q^{-1}(-\alpha_0)=0$ alone, and $0$ is not a critical value.
By Lemma~\ref{lem:artin}, applied in the fundamental group of $\C$ minus the
critical values, based at $u_b$, with a geometric basis whose product is the
boundary loop of $\widetilde\Omega_1$, these loops are Hurwitz equivalent to
such a geometric basis.  All of them except the $\tilde\ell_i$ lie in
$\C\setminus\operatorname{int}\Omega$.

Along loops that avoid the critical values of $\pr|_C$ and the
$U$-coordinates of the points of $C\cap\{V=-\beta_0\}$, the fibers of $C$
lie in $\C\setminus\{-\beta_0\}$, so the braid monodromy of $C$ along such
loops lies in $B_m(\C\setminus\{-\beta_0\})$.  By~\ref{pb2}, the local factors
$b_i$ lie in $B_m$, and the germs of $D$ at its critical points over $\Omega$
lie in the half-planes $\pm\Re v>0$, where $\psi$ maps them isomorphically
onto germs of $C$ by $v\mapsto v^2-\beta_0$.  So the proof of
Lemma~\ref{lem:lift} shows that the local factors of $D$ along the
$\tilde\ell_i$ are the $2N$ braids $\Delta_{m,0}b_i\Delta_{m,0}^{-1}$ and
$\sh(b_i)$, the monodromy along each $\tilde\ell_i$ being the image under
$\Phi$ of that along $\ell_i$; write $\Phi(s)$ for this sequence.
The monodromy $n_j$ along $\nu_j$ is a conjugate of $\tau_1$, since near
the corresponding point of $C\cap\{V=-\beta_0\}$ one point of the fiber turns
once around $-\beta_0$; so, by~\eqref{eq:tauconj}, $n_j$ is conjugate in $F$
to some $\tau_k$.  Moreover $n_1\cdots n_m=\tau_1\cdots\tau_m$, since
$b_1\cdots b_N=\Delta_m^2$, $\delta$ has trivial monodromy, and the monodromy
of $C\cup\{V=-\beta_0\}$ along the boundary loop of $\Omega_1$ is the full
twist $\Delta_m^2\tau_1\cdots\tau_m$ of $m+1$ strings: by~\ref{pb3}, the line
$\{V=-\beta_0\}$ meets $C$ in $m$ affine points, so its point at infinity is
not on $\overline C$, and $C\cup\{V=-\beta_0\}$ is a curve of degree $m+1$,
transverse to the line at infinity, whose critical values all lie in the
interior of $\Omega_1$.  By Lemma~\ref{lem:artin},
applied in $F$ with the basis $\tau_1,\dots,\tau_m$, and Lemma~\ref{lem:lift},
the $\tilde\nu_j$ contribute a factorization Hurwitz equivalent to $\bar r_m$.
The path $\tilde\delta$ lifts the loop $\delta$, along which the points of $C$
move by a null-homotopic motion in $\C\setminus\{-\beta_0\}$, so its braid is
trivial; and the loops conjugated by $\tilde\delta$ have the same monodromies
and local factors as the $\tilde\ell_i$ and $\tilde\nu_j$, by the symmetry
$u\mapsto-u$ of $D$.  Hence $\operatorname{bmf}(D)$ is Hurwitz equivalent to
$\bigl(\Phi(s),\bar r_m,\Phi(s),\bar r_m\bigr)$.  For $D^\gamma$,
Lemma~\ref{lem:twist} conjugates the local factors along the $\tilde\ell_i$
by $\sh(\gamma)$ and leaves the others unchanged; as $\sh(\gamma)$ commutes
with $\Delta_{m,0}b_i\Delta_{m,0}^{-1}$, this replaces the factors $\sh(b_i)$ in
the first block by $\sh(\gamma b_i\gamma^{-1})$.

Conjugate by $\Delta_{m,0}^{-1}$.  The first block becomes
$\bigl(s,\sh(\lambda(\gamma)s)\bigr)$, the second $(s,\sh(s))$, and $\bar r_m$
becomes $r'=\lambda(\Delta_{m,0}^{-1})\bar r_m$.  By~\eqref{eq:halftwist},
the product of $r'$ is $\Delta_{m,0}^{-2}\Delta_{m,m}^{-2}\,
\Delta_{m,0}^{-1}\Delta_{2m}\Delta_{m,0}$, and conjugation by it maps $b\,\sh(c)$
to $c\,\sh(b)$ for $b,c\in B_m$; so moving the first copy of $r'$ to the
right past the second block leaves that block unchanged.  Conjugating by
$\Delta_{2m}$ exchanges the two halves of the strings, acting by
$b\mapsto\Delta_mb\Delta_m^{-1}$, and maps $r'$ to
$\lambda(\Delta_{m,m}^{-1}\Delta_{2m})\bar r_m=\tilde r_m$.  After permuting
factors with disjoint supports, we obtain
$d(\lambda(\Delta_m\gamma)s,\lambda(\Delta_m)s,\lambda(\Delta_m)s,
\lambda(\Delta_m)s)$.  The last assertion follows from
Lemma~\ref{lem:framing} with $h=\Delta_m$.
\end{proof}

\subsection{Pencils of quartics}
\label{subsec:pencil}
Recall the curves $D_n$, the region $\Omega$ and the point $u_b$ of
\S\ref{subsec:intro-model}, and put $s_n=d(\Delta_n^2)$ and
$\omega=\Delta_{n,0}\Delta_{2n}^{-1}\in B_{2n}$; so $s_3$ is the factorization $s$
of Table~\ref{tab:s}.  For $\gamma\in B_{2n}$ let $D_n^\gamma$ be obtained
from $D_n$ by the twist of \S\ref{subsec:pullback}, with $m=2n$.

\begin{proof}[Proof of Theorem~\ref{thm:model}]
For $\lambda\ne0,\infty$, the member $Q_n(v)=\lambda Q_1(u)$ of the pencil
meets the line at infinity transversally, in the four points with
$v^4=\lambda u^4$, and an affine point of it is singular only where
$Q_1'(u)=Q_n'(v)=0$.  The critical values of $Q_1$ are $Q_1(0)=3$ and
$Q_1(\pm\sqrt2)=-1$, those of $Q_n$ are $Q_n(0)=n(n+2)$ and
$Q_n(\pm\sqrt{n+1})=-1$, so the members with
$\lambda\notin\{0,\infty,1,-1/3,n(n+2)/3,-n(n+2)\}$ are smooth; 
this includes $\lambda_i=i+1$, $1\le i\le n$.  Two members
meet only at the base points, where $Q_1$ and $Q_n$ vanish simply; the
member with parameter $\lambda$ has slope $\lambda Q_1'(u)/Q_n'(v)$ there, so the
$n$ members are pairwise transverse and the base points are ordinary
$n$-fold points.

Let $L\subset\C^2$ be the union of the $n$ lines $Y-1=\lambda_i(X-1)$, and
$C_n=\{\prod_i(V^2-1-\lambda_i(U^2-1))=0\}$.  Then $C_n=\psi_0^{-1}(L)$ with
$\psi_0(U,V)=(U^2,V^2)$, and $D_n=\psi^{-1}(C_n)$ with
$\psi(u,v)=(u^2-2,v^2-n-1)$, since $Q_1(u)=(u^2-2)^2-1$ and
$Q_n(v)=(v^2-n-1)^2-1$.

First apply Lemma~\ref{lem:geomdouble} to $L$, with $\alpha_0=\beta_0=0$,
$\Omega_0=\{\lvert X-1\rvert\le\rho\}$, $\rho=1/(n+2)$, $X_b=1+\rho$ and
$\Omega_1=\{\lvert X\rvert\le1+\rho\}$.  The
only critical value is $X=1$, with factor $\Delta_n^2$; over $\Omega_0$ we
have $\Re Y\ge1-(n+1)\rho>0$; the points over $X_b$ are real; and the lines
meet $Y=0$ transversally at the distinct points $X=1-1/\lambda_i\notin\{0\}$,
which lie, with $0$, in the interior of $\Omega_1$.
Since $\lambda(\Delta_n)$ fixes $\Delta_n^2$, the lemma, with $\gamma=1$, gives
$\operatorname{bmf}(C_n)=\lambda(\omega)\,s_n$ with respect to the base point
$U_b=\sqrt{1+\rho}=r_n$, the boundary loop of the preimage
$\{\lvert U\rvert\le r_n\}$ of $\Omega_1$, and the real order.  As this disk
contains all critical values of $\pr|_{C_n}$, $\lambda(\omega)s_n$ is the
braid monodromy factorization of $C_n$ over $\{\lvert U\rvert\le r_n\}$
with respect to $U_b$.

Next apply Lemma~\ref{lem:geomdouble} to $C_n$, with $m=2n$, $\alpha_0=2$,
$\beta_0=n+1$, $\Omega_0=\{\lvert U\rvert\le r_n\}$, and any $\Omega_1$ as in
the lemma, for instance a large round disk with a thin open sector removed,
with vertex $U_b$ and bisector $U_b+\R_{>0}e^{\pi i/4}$ (a round disk tangent
to $\Omega_0$ at $U_b$ will not do, as $C_n$ meets $V=-n-1$ also at real
points $U>r_n$).  The critical values
$\pm1$ and $\pm\sqrt{1-1/\lambda_i}$ of $\pr|_{C_n}$ lie in the interior
of $\Omega_0$, as $r_n>1$, and
$\Re(U+2)\ge2-r_n>0$ there.  For $(U,V)\in C_n$ with $U\in\Omega_0$,
\[
\lvert V\rvert^2\le1+(n+1)(r_n^2+1)<(n+1)^2,
\]
so $\Re(V+n+1)>0$; the points of $C_n$ over $U_b=r_n$ are the real numbers
$\pm\sqrt{1+\lambda_i\rho}$.  Finally $C_n$ meets $V=-n-1$ where
$U^2=1+n(n+2)/\lambda_i$, transversally, in $2n$ points with distinct
$U$-coordinates; none of them is $-2$, since $\lambda_i\ne n(n+2)/3$.  With
$s=\operatorname{bmf}(C_n)=\lambda(\omega)s_n$, the lemma shows that
$\operatorname{bmf}(D_n^\gamma)$ has the type of
\[
d\bigl(\lambda(\gamma\omega)s_n,\lambda(\omega)s_n,
\lambda(\omega)s_n,\lambda(\omega)s_n\bigr)
=\lambda\bigl(\omega\,\sh(\omega)\bigr)\,
d\bigl(\lambda(\omega^{-1}\gamma\omega)s_n,s_n,s_n,s_n\bigr),
\]
by Lemma~\ref{lem:framing} with $h=\omega$.  For $\gamma=1$ this gives
part~\ref{mod1}, and in general part~\ref{mod2}.

For part~\ref{mod3}, recall that $\overline H_k$ and $\overline H_\beta$, for
$\beta\in B_m$ as in Theorem~\ref{thm:twisted}, have braid monodromy
factorizations
\[
T_k=d(\lambda(b_k)s_3,s_3,s_3,s_3),\qquad
U_\beta=d(\lambda(\beta)s_{m+1},s_{m+1},s_{m+1},s_{m+1}),
\]
which by part~\ref{mod2}, with $\gamma=\omega b_k\omega^{-1}$ and
$\gamma=\omega\beta\omega^{-1}$, are also those of
$D_3^{\,\omega b_k\omega^{-1}}$ and $D_{m+1}^{\,\omega\beta\omega^{-1}}$, up to
type.  The critical points of all these curves are ordinary multiple points
and simple tangencies, so part~\ref{mod3} follows from \S\ref{subsec:hurwitz}.
The separations cannot be omitted, as $D_n^\gamma$ has several critical
points in some fibers.  Its critical points are those of $D_n$: the $16$ base
points and the $12n$ simple tangencies, at the points with $Q_n'(v)=0$.  The
fibers of $\pr$ that contain more than one of them are
those over $u=\pm1,\pm\sqrt3$, each containing four $n$-fold points, and those
over the $4n$ solutions of $Q_1(u)=-1/\lambda_i$, each containing two
tangencies, at $v=\pm\sqrt{n+1}$; the other $4n$ tangencies, at $v=0$ over
the solutions of $Q_1(u)=n(n+2)/\lambda_i$, lie in distinct fibers.
\end{proof}

\begin{remark}
\label{rem:twistbands}
For $n=3$ and $k=1,2$, the braids $\omega b_k\omega^{-1}$ of
Theorem~\ref{thm:model}\ref{mod3} are again bands, acting on the six positive
points of $D_3$ over $u_b$, numbered from left to right in the $v$-line,
that is, from bottom to top in Figure~\ref{fig:D3}(a), where they are the
points of $D_3$ on the right edge of the shaded strip.  A direct
computation with \eqref{eq:action} shows that
\[
\omega b_1\omega^{-1}=g_1'\,\sigma_1\,g_1'^{-1},\qquad
g_1'=\sigma_5^{-1}\sigma_4\sigma_5^{-1}\sigma_4\sigma_3^{-1}\sigma_2^{-1},
\]
with endpoint loops $g_1'(x_1)=x_1$ and $g_1'(x_2)=c\,x_5\,c^{-1}$, where
$c=x_4x_6x_4^{-1}x_6^{-1}$, and that
\[
\omega b_2\omega^{-1}=g_2'\,\sigma_5\,g_2'^{-1},\qquad
g_2'=\sigma_4\sigma_1\sigma_2\sigma_3^{-1}\sigma_2^{-1}
\sigma_1^{-3}\sigma_2^{-1}\sigma_3\sigma_4,
\]
with $g_2'(x_6)=x_6$ and $g_2'(x_5)=q\,x_1x_2x_1^{-1}q^{-1}$, where
$q=p\,x_1x_3^{-1}x_2^{-1}x_1^{-1}\,p$ and $p=x_4x_5^{-1}x_4^{-1}$.  Thus
$\overline H_1$ and $\overline H_2$ are obtained from $D_3$ by twisting along
bands that join the first and fifth, respectively the second and sixth, of
these points.  As $Q_1(u_b)=r_3^2-1=1/5$, the positive points over $u_b$ of
the quartic with parameter $\lambda$ are $v=\bigl(4\pm(1+\lambda/5)^{1/2}
\bigr)^{1/2}$, so from left to right the six points lie on the quartics with
$\lambda=4,3,2,2,3,4$.  In both cases the band therefore joins the quartics
with $\lambda=3$ and $\lambda=4$, and its permutation fixes the two points of
the quartic with $\lambda=2$, so that the twist reglues this quartic to
itself.  As $\overline H_k$ has two components, of degrees $8$ and $4$
(Theorem~\ref{thm:explicit}\ref{a3}), the quartic with $\lambda=2$ becomes the
component of degree $4$, and the other two merge into the component of
degree $8$.
The words
$g_1'$ and $g_2'$, of lengths $6$ and $11$, are the shortest found by a
computer search.
\end{remark}

\begin{remark}
\label{rem:model}
By \eqref{eq:Gdouble}, $\pi_1(\C^2\setminus D_n)\cong G(s_n)\cong
G(\Delta_n^2)\cong F_{n-1}\times\Z$.  Thus a single twist of $D_3$ along a
circle replaces the quasi-projective group $F_2\times\Z$ by the
non-Hopfian group $\BS(2,3)\times\Z$, and the two quartics of
Remark~\ref{rem:twistbands} by the component of degree $8$ and geometric
genus $5$ of $\overline H_2$.  Similarly $\overline H_{\gamma_n}$ is a twist 
of the union $D_4$ of four quartics, and the twist does not merge components, as
$\gamma_n$ is pure.
\end{remark}

\section{On the minimal degree}
\label{sec:mindeg}

Since $\BS(k,k+1)\times\Z$ and $\pi_1(S^3\setminus L_n)\times\Z$ are not
quasi-projective (Theorem~\ref{thm:alex} and
Proposition~\ref{prop:multivar-fails}), any Hurwitz curve with one of these
groups is non-algebraic.  Isotopy theorems that would exclude low degrees
are scarce: symplectic isotopy is known for smooth curves of degree at most
$17$ \cite{Siebert-Tian-ann05} and for rational cuspidal curves of degree at
most $5$ \cite{Golla-Starkston-comp22}, but not for curves with the
singularities that occur here.  We prove one for curves of degree $3$ with
two components (Proposition~\ref{prop:deg3}), which gives the lower bound
$4$; otherwise our lower bounds are group-theoretic.

In \S\ref{subsec:mindeg-BS} we determine the numerical invariants forced on
a curve with group $\BS(k,k+1)\times\Z$, show that degree $12$ is the smallest
reached by twisted doublings of iterated doublings of full twists, classify
the two-component Hurwitz curves of degree $3$, which shows that the degree is
at least $4$ and excludes doublings of degree $6$, and discuss a natural
attempt in degree $10$.  In \S\ref{subsec:mindeg-links}
we show that, for a link with vanishing linking numbers, the components of
any Hurwitz curve realizing its group times $\Z$ meet like members of a
pencil; this bounds the degree from below and shows that, with its
C-structure, Theorem~\ref{thm:twisted} is optimal among Kulikov doublings
when $\beta$ is pure.

\subsection{Baumslag--Solitar groups}
\label{subsec:mindeg-BS}

\begin{proposition}
\label{prop:degconstraints}
Let $\overline H$ be a Hurwitz curve of degree $m$ with
$\pi_1(\C^2\setminus H)\cong\BS(k,k+1)\times\Z$, $k\ge1$.  Then $\overline H$
has two components $A,B$.  Write their meridian classes in the basis
$[t],[z]$ of $H_1=\Z^2$ as $\epsilon_A=(u_A,v_A)$ and
$\epsilon_B=(u_B,v_B)$.  Then $u_Au_B<0$, and with
$h=\lvert u_A\rvert+\lvert u_B\rvert\ge2$ there is $n\ge1$ such that
\[
(\deg A,\deg B)=n\,(\lvert u_B\rvert,\lvert u_A\rvert),\qquad m=nh,\qquad
l_m=z^{\pm n}.
\]
The torsion subgroup of the Alexander invariant $A_0$ (for the total-degree
map) is cyclic of order $(k+1)^h-k^h$, and the meridians act on it with
order $h$.
\end{proposition}

\begin{proof}
The meridian classes form a basis of $H_1$ (Proposition~\ref{prop:H1}), so
there are two components and $u_A,u_B$ are coprime.  The center of
$\BS(k,k+1)\times\Z$ is $\la z\ra$, so $l_m=z^e$ for some $e$ by
Lemma~\ref{lem:central}.  In $H_1$ this reads
$\deg A\cdot\epsilon_A+\deg B\cdot\epsilon_B=(0,e)$, whence
$\deg A\cdot u_A=-\deg B\cdot u_B$.  The coefficient $u_B$ cannot vanish,
since then $u_A=\pm1$ and $\deg A=0$; likewise $u_A\ne0$.  This gives the
degrees, and
$e=n(\lvert u_B\rvert v_A+\lvert u_A\rvert v_B)=\pm n\det(\epsilon_A,\epsilon_B)=\pm n$.

Let $\nu$ be the total-degree map, so $\nu(\epsilon_A)=\nu(\epsilon_B)=1$
and $\nu(z)=m/e=\pm h$.  Since $\nu$ is onto, $\tau_0=\nu(t)$ is prime to
$h$.  The kernel $K$ of $\BS(k,k+1)\to\Z$, $t\mapsto1$, $a\mapsto0$, has
$H_1(K)=A_0(G_k)\cong\Z[1/k(k+1)]$, with $t$ acting by multiplication by
$\varrho=((k+1)/k)^{\pm1}$ \cite[Ex.~3.20]{Kulikov-izv08}.  Now
$\ker\nu=K\rtimes\la w\ra$ with $w=t^hz^{\mp\tau_0}$ when $\nu(z)=\pm h$, so
$A_0=H_1(\ker\nu)\cong H_1(K)/(\varrho^h-1)\oplus\Z$.  The first summand is the
derived module $A_h\bigl(A_0(G_k)\bigr)\cong\Z/((k+1)^h-k^h)$
\cite[Ex.~2.14]{Kulikov-izv08}, as $(k+1)^h-k^h$ is prime to $k(k+1)$.  A meridian of
$A$ acts by $\varrho^{u_A}$, of order $h/\gcd(u_A,h)=h$, because
$0<(k+1)^j-k^j<(k+1)^h-k^h$ for $0<j<h$; similarly for $B$.
\end{proof}

For the curves of
Theorem~\ref{thm:explicit}, $(u_A,u_B)=(1,-2)$ and $(n,h)=(4,3)$; the torsion
is $\Z/(3k^2+3k+1)$, as in Theorem~\ref{thm:alex}.  The torsion is an invariant
of the curve, so realizations with different $h$ are distinguished by it.

Degree $12$ is optimal for the constructions used in this paper.  For $j\ge0$
write $d^{\,0}(s)=s$ and $d^{\,j+1}(s)=d\bigl(d^{\,j}(s)\bigr)$.

\begin{proposition}
\label{prop:min12}
Let $p\ge2$ and $j\ge0$, put $M=p2^j$ and $s=d^{\,j}(\Delta_p^2)\in\cS_M$, and
let $T=d\bigl(\lambda(\beta_1)s,\lambda(\beta_2)s,\lambda(\beta_3)s,
\lambda(\beta_4)s\bigr)$ with $\beta_1,\dots,\beta_4\in B_M$.  If
$G(T)\cong\BS(k,k+1)\times\Z$ for some $k\ge1$, then $T$ has degree
$2M\ge12$.
\end{proposition}

\begin{proof}
By \eqref{eq:Gdouble} and \eqref{eq:Nconj}, $G(T)\cong
F_M/\beta_1(N(s))\cdots\beta_4(N(s))$, a quotient of
$F_M/\beta_1(N(s))\cong G(s)$, and $G(s)\cong G(\Delta_p^2)=\la
x_1,\dots,x_p\mid[x_i,x_1\cdots x_p]\ (i\le p)\ra\cong F_{p-1}\times\Z$.  The
group $\BS(k,k+1)\times\Z$ is not abelian, and its abelianization is
$\Z^2$.  Suppose that $2M<12$, that is, $M\le5$.  If $j=0$, then $s=\Delta_M^2$
is central, so $\lambda(\beta_i)s=s$ and $G(T)\cong F_{M-1}\times\Z$, which
is abelian for $M=2$ and has abelianization $\Z^M\not\cong\Z^2$ for $M\ge3$.
If $j\ge1$, then $p2^j\le5$ forces $p=2$ and $j=1$, and $G(T)$ is a
quotient of $G(\Delta_2^2)\cong\Z^2$, hence abelian.
\end{proof}

Proposition~\ref{prop:min12} covers Kulikov's band lemma, which uses
$d\bigl(\lambda(b)(d(s)),d(s),d(s),d(s)\bigr)$, as well as
Lemma~\ref{lem:band} and the twisted doublings of \S\ref{sec:twisted}, applied
to iterated doublings of full twists; by Theorem~\ref{thm:explicit}, with
$p=3$, $j=1$, $\beta_1=b_k$ and $\beta_2=\beta_3=\beta_4=1$, the bound is
attained.  It says nothing about curves of lower degree obtained in other
ways; see Corollary~\ref{cor:deg3} and Question~\ref{q:mindegree}.

Curves of degree $3$ with two components can be classified completely.  We
first classify the relevant factorizations.

\begin{lemma}
\label{lem:deg3}
Let $B_{1,2}\subset B_3$ be the stabilizer of the first string, put
$a=\sigma_2$, $b=\sigma_1^2$ and $\theta=ab$, and let $s$ be a factorization
of $\Delta_3^2$ whose factors lie in $B_{1,2}$ and are conjugate in $B_3$
to $\sigma_1$, $\sigma_1^2$, $\sigma_1^4$ and $\Delta_3$, the numbers of
factors in these four classes being
\[
\text{(I)}\ \ 2,2,0,0,\qquad \text{(II)}\ \ 2,0,1,0,\qquad
\text{(III$_1$)}\ \ 1,1,0,1,\qquad\text{or}\qquad \text{(III$_2$)}\ \ 0,0,0,2.
\]
Then $s$ is Hurwitz equivalent to $\lambda(g)t$ for some $g\in B_{1,2}$, where
$t$ is $(a,a,b,aba^{-1})$, $(b^{-1}ab,a,b^2)$, $(a,b,\theta)$, or
$(\theta,\theta)$, respectively.  In each case $G(s)\cong\Z^2$.
\end{lemma}

\begin{proof}
The braid relation shows that
$\psi(\sigma_1)=\left(\begin{smallmatrix}1&0\\-1&1\end{smallmatrix}\right)$ and
$\psi(\sigma_2)=\left(\begin{smallmatrix}1&1\\0&1\end{smallmatrix}\right)$
define a homomorphism $\psi\colon B_3\to\SL_2(\Z)$.  By
\cite[Ch.~I, \S1.5.2--1.5.3]{Serre-80}, $B_3=\la\Delta_3,\sigma_1\sigma_2\mid
\Delta_3^2=(\sigma_1\sigma_2)^3\ra$, and $\SL_2(\Z)$ is the amalgam, over
$\{\pm I\}$, of the cyclic groups of orders $4$ and $6$ generated by
$J=\psi(\Delta_3)=\left(\begin{smallmatrix}0&1\\-1&0\end{smallmatrix}\right)$ and
$\psi(\sigma_1\sigma_2)^{-1}$.  So $\psi$ is onto, with kernel the normal
closure of $\Delta_3^4$, which is $\la\Delta_3^4\ra$; as conjugate braids
have the same exponent sum, two conjugate braids with the same image under
$\psi$ are equal.  Reduction modulo $2$ kills $\sigma_1^2$, hence the pure
braid group, and maps $B_3$ onto $\SL_2(\Z/2)$, a group of order $6$; the
induced isomorphism from $S_3$ maps the transposition $(2\,3)$, the image of
$\sigma_2$, to an upper triangular matrix.  Hence
$B_{1,2}=\psi^{-1}(\Gamma_0(2))$, where $\Gamma_0(2)\subset\SL_2(\Z)$ 
consists of the matrices with even lower left entry.

For a primitive vector $u\in\Z^2$ put $T_u=I+uu^{\top}J$, which depends only
on $\pm u$.  Then $\psi(\sigma_2)=T_{e_1}$, $\psi(\sigma_1)=T_{e_2}$, and
$MT_uM^{-1}=T_{Mu}$ for $M\in\SL_2(\Z)$, as $M^{\top}JM=J$.  Since
$u^{\top}Jv=\det(u,v)$, we have $(uu^{\top}J)^2=0$ and
$\tr(uv^{\top}J)=-\det(u,v)$, so for all integers $k,l$
\begin{equation}
\label{eq:trace}
T_u^kT_v^l=I+\bigl(k\,uu^{\top}+l\,vv^{\top}+kl\det(u,v)\,uv^{\top}\bigr)J,
\end{equation}
and $\tr(T_u^kT_v^l)=2-kl\det(u,v)^2$.
Thus a factor conjugate to $\sigma_1^k$ has image $T_u^k$ for a primitive
vector $u$, its \emph{vector}, which is determined up to sign and determines
the factor.  If $A\in B_{1,2}$ is conjugate to $\sigma_1$, then
$T_u\in\Gamma_0(2)$, so the second coordinate of its vector $u$ is even;
hence every $M\in\SL_2(\Z)$ with first column $u$ is $\psi(g)$ for some
$g\in B_{1,2}$, and then $g^{-1}Ag=a$.

Hurwitz moves permute the classes of the factors, so we may order the factors
as in $t$.  We write $A,A_i$, $B_i$, $C$, and $\Theta,\Theta_i$ for factors
conjugate to $\sigma_1$, $\sigma_1^2$, $\sigma_1^4$, and $\Delta_3$, 
and $u,u_i$, $w_i$, and $w$ for the vectors of the first three kinds.  Note that
$\psi(\Delta_3^2)=-I$, that $\theta=\sigma_1^{-1}\Delta_3\sigma_1$, so
that $\theta^2=\Delta_3^2$, and that
$\tr(-M^{-1})=-\tr M$ for $M\in\SL_2(\Z)$.

(III$_2$) If $s=(\Theta_1,\Theta_2)$, then
$\Theta_2=\Theta_1^{-1}\Delta_3^2=\Theta_1$, since $\Theta_1^2$ is a
conjugate of $\Delta_3^2$.  If $\Theta_1=h\theta h^{-1}$, then the
permutation of $h$ commutes with $(2\,3)$, the permutation of $\theta$ and
of $\Theta_1$, so $h\in B_{1,2}$ and $s=\lambda(h)t$.

(III$_1$) If $s=(A,B_1,\Theta)$, then $\psi(AB_1)=-\psi(\Theta)^{-1}$ has
trace $0$, since $\psi(\Theta)$ is conjugate to $J$, and \eqref{eq:trace} gives
$\det(u,w_1)=\pm1$.  Choosing the sign of $w_1$, the matrix $(u\ w_1)$ with
columns $u$ and $w_1$ is $\psi(g)$ for some $g\in B_{1,2}$; then
$g^{-1}Ag=a$, $g^{-1}B_1g=b$ and $g^{-1}\Theta g=(ab)^{-1}\Delta_3^2=\theta$.

(II) If $s=(A_1,A_2,C)$, then $\psi(A_1A_2)=-T_w^{-4}=-I+4ww^{\top}J$, and
\eqref{eq:trace} gives $d=\det(u_1,u_2)=\pm2$ and, after multiplication by
$J^{-1}=-J$, $u_1u_1^{\top}+u_2u_2^{\top}+d\,u_1u_2^{\top}=2J+4ww^{\top}$.  As
$d^2=4$, the symmetric part of the left side is
$(u_1+\frac d2u_2)(u_1+\frac d2u_2)^{\top}$, so $u_1+\frac d2u_2=\pm2w$, and
then $d=\pm2\det(w,u_2)$, so $\det(u_2,w)=\pm1$.  Conjugating as in
(III$_1$), we may assume that $u_2=e_1$ and $w=e_2$, that is, $A_2=a$ and
$C=b^2$.  Then $u_1=\pm2e_2-\frac d2e_1$ and $d=\det(u_1,e_1)=\mp2$, so
$u_1=\pm(e_1+2e_2)=\pm\psi(b)^{-1}e_1$ and $A_1=b^{-1}ab$.

(I) If $s=(A_1,A_2,B_1,B_2)$, then $\psi(A_1A_2)=-\psi(B_1B_2)^{-1}$, and
comparing traces by \eqref{eq:trace} gives
$\det(u_1,u_2)^2+4\det(w_1,w_2)^2=4$.  So either $w_1=\pm w_2$ or
$u_1=\pm u_2$.  In the first case $B_1=B_2$, and $(A_1,A_2,B_1^2)$ is as in
(II); so we may assume that it is $(b^{-1}ab,a,b^2)$, and then $B_1=b$, as
$T_{w_1}^4=T_{e_2}^4$ forces $w_1=\pm e_2$.  Two Hurwitz moves at the second
place now give
\[
(b^{-1}ab,a,b,b)\mapsto(b^{-1}ab,b,b^{-1}ab,b)\mapsto
(b^{-1}ab,b^{-1}ab,b^{-1}a^{-1}bab,b)=\lambda(b^{-1}a^{-1})\,t.
\]
In the second case $A_1=A_2$, and we may assume that both are $a$.  Then
$\psi(B_1B_2)=\psi(a^{-2}\Delta_3^2)=-T_{e_1}^{-2}$, and \eqref{eq:trace}
gives $d=\det(w_1,w_2)=\pm1$ and
$w_1w_1^{\top}+w_2w_2^{\top}+2d\,w_1w_2^{\top}=J+e_1e_1^{\top}$, whose
symmetric part gives $w_1+dw_2=\pm e_1$.  Hence $d=\pm\det(e_1,w_2)$, so
that, up to sign, $w_2=\psi(a^k)e_2=(k,1)^{\top}$ for some $k\in\Z$, and then
$w_1=\pm(e_1-w_2)=\mp\psi(a^{k-1})e_2$.  So $B_1=a^{k-1}ba^{1-k}$,
$B_2=a^kba^{-k}$ and $s=\lambda(a^{k-1})t$.

Finally, $G(s)\cong G(t)$ by \eqref{eq:Nconj}.  By \eqref{eq:action}, the
relations of $G(t)$ include $x_2=x_3$, given by $a$, or by $\theta$, which maps
$x_3$ to $x_2$, and then $[x_1,x_2]=1$, given by $b$, by $\theta$, or by
$b^{-1}ab$.  Conversely, all the relations hold in $\Z^2=\la x_1,x_2\ra$, with
$x_3\mapsto x_2$, since the permutation of every factor fixes $1$.  So
$G(s)\cong\Z^2$.
\end{proof}

\begin{proposition}
\label{prop:deg3}
Let $\overline H=C\cup L$ be a Hurwitz curve of degree $3$, where $C$ is
irreducible of degree $2$ and $L$ has degree $1$.  Then $C$ is smooth and is
simply tangent to exactly two lines of the pencil, and
$\operatorname{bmf}(\overline H)$ has the type of the braid monodromy
factorization of an algebraic curve, the union of a smooth conic and a line
that meet transversally, off the points of tangency of the conic with lines of
the pencil (I), or are tangent (II), or such that the line passes through one
(III$_1$) or both (III$_2$) of these points.  In particular
$\pi_1(\C^2\setminus H)\cong\Z^2$.  If the critical points of $\overline H$ lie
in distinct fibers, then $\overline H$ is H-isotopic to such an algebraic curve;
in general, every separation of $\overline H$ is.
\end{proposition}

\begin{proof}
Choose the identification of the base fiber so that $L$ is the first string.
The factors of $s=\operatorname{bmf}(\overline H)$ are conjugate to braid
monodromies of analytic germs, hence positive, and they lie in $B_{1,2}$, as
$L$ is a component.  Forgetting the first string gives the braid monodromy
factorization of $C$, a factorization of $\sigma_1^2$ in $B_2$ into
non-negative powers of $\sigma_1$, not all even since $C$ is irreducible; so
its non-trivial factors are two copies of $\sigma_1$.  Since the exponent
sum of the braid monodromy of a germ with $k$ strings and Milnor number $\mu$
is $\mu+k-1$, by the Riemann--Hurwitz formula for its Milnor fiber, it 
follows that $C$ is smooth and that it is tangent to exactly
two lines of the pencil, simply.  As $C\cdot L=2$, the local intersection
numbers are positive and $L$ is transverse to the fibers, $C\cap L$ consists
of two transverse points off the points of tangency (I), or of one point of
contact $2$ (II), or contains one (III$_1$) or both (III$_2$) of the points of
tangency, at which $C$ and $L$ are transverse.  So the germs of $\overline H$
at its critical points are those of $v^2=u$, $v^2=u^2$, $v^2=u^4$ and
$v(v^2-u)=0$, with braid monodromies $\sigma_1$, $\sigma_1^2$, $\sigma_1^4$ and
$\Delta_3$, the last one turning three collinear points by a half-turn, in the
numbers of Lemma~\ref{lem:deg3}; by the lemma, $s$ has the type of the
corresponding factorization $t$.

Let $C$ be a smooth conic, $p\notin C$ the center of the pencil, and
$q_1,q_2$ the points of $C$ at which it is tangent to lines through $p$.  A
general line, the tangent line at a general point of $C$, a general line
through $q_1$, and the line through $q_1$ and $q_2$, which is the polar line
of $p$ and does not contain $p$, give algebraic curves of types I, II,
III$_1$ and III$_2$ whose critical points lie in distinct fibers; by what we
have shown, their braid monodromy factorizations have the types of the
corresponding factorizations $t$.  The group is $\Z^2$ by
Theorem~\ref{thm:zvk} and Lemma~\ref{lem:deg3}.  The germs $v^2=u^j$ are of
inseparable type by \cite[Lem.~5.5]{Kulikov-usp07}, and so is $v(v^2-u)=0$, by
the proof of that lemma: an element of $F_3$ fixed by $\Delta_3$ is fixed by
$\Delta_3^2$, which acts by conjugation by $l_3$ \eqref{eq:fulltwist}, so it lies
in the centralizer $\la l_3\ra$ of $l_3$.  The last assertion therefore follows
from \cite[Thm.~5.2]{Kulikov-usp07}, as in \S\ref{subsec:hurwitz}.
\end{proof}

\begin{corollary}
\label{cor:deg3}
Let $k\ge1$.
\begin{enumerate}[label=(\alph*),itemsep=1pt]
\item \label{g1}
Every Hurwitz curve $\overline H$ with
$\pi_1(\C^2\setminus H)\cong\BS(k,k+1)\times\Z$ has degree at least $4$.
\item \label{g2}
If $s_1,\dots,s_4\in\cS_3$ are braid monodromy factorizations of Hurwitz
curves, then the group of $d(s_1,s_2,s_3,s_4)$ is not $\BS(k,k+1)\times\Z$.
\end{enumerate}
\end{corollary}

\begin{proof}
The group $\BS(k,k+1)\times\Z$ is not abelian, and its abelianization is
$\Z^2$.  \ref{g1} By Proposition~\ref{prop:degconstraints}, $\overline H$ has
two components, of degrees $n\lvert u_B\rvert$ and $n\lvert u_A\rvert$, where
$\lvert u_A\rvert+\lvert u_B\rvert\ge2$.  In degree $2$ both components are
lines, the braid monodromy factorization is $(\sigma_1^2)$, and the group is
$\Z^2$.  In degree $3$ we have $n=1$ and $h=3$, so the components have
degrees $2$ and $1$, and the group is $\Z^2$ by
Proposition~\ref{prop:deg3}.

\ref{g2} By \eqref{eq:Gdouble}, $G(T)$, where $T=d(s_1,s_2,s_3,s_4)$, is a
quotient of each $G(s_i)$.  If $G(T)\cong\BS(k,k+1)\times\Z$, then $T$ has two
components.  If all factors of the $s_i$ were pure, then by \eqref{eq:cmi}
the components of $T$ would be given by the pairs of strings $\{i,3+i\}$,
$i=1,2,3$; so some $s_i$ has a non-pure factor, and is the braid monodromy
factorization of an irreducible curve or of a conic and a line.  In the
first case $H_1(G(s_i))\cong\Z$, so that $H_1(G(T))$ is cyclic; in the second
$G(s_i)\cong\Z^2$ by Proposition~\ref{prop:deg3}, so that $G(T)$ is abelian.
\end{proof}

\begin{remark}
\label{rem:deg3}
For curves of type I, which are nodal, Proposition~\ref{prop:deg3} confirms,
in this simple case and when the critical points lie in distinct fibers,
Kulikov's conjecture that nodal Hurwitz curves are H-isotopic to algebraic
curves \cite[Conj.~5.1]{Kulikov-usp07}.  Factorizations in $B_3$ are
usually studied through $\PSL_2(\Z)$ \cite[Ch.~10]{Degtyarev-12}, where the
factorizations into conjugates of $\sigma_1$ are known to be unique, up to
Hurwitz moves, for the powers of the full twist and, more generally, for
positive braids \cite[Thms.~10.18, 10.19]{Degtyarev-12}, and those of
length $2$ are classified \cite[Thms.~10.27, 10.30]{Degtyarev-12};
Lemma~\ref{lem:deg3} also involves conjugates of $\sigma_1^2$, $\sigma_1^4$,
and $\Delta_3$.
\end{remark}

\begin{remark}
\label{rem:deg10}
The most natural degree-$10$ candidate, a single band over a factorization
of degree $5$ with group $F_2\times\Z$, fails for every band that joins the
two conics and satisfies the hypothesis of Lemma~\ref{lem:fixband}.  The conics
$C_\lambda=\{yz=x^2+\lambda y^2\}$ form a pencil with a single base point,
of contact order $4$, and for $\lambda_1\ne\lambda_2$ the complement of
$C_{\lambda_1}\cup C_{\lambda_2}\cup\{y=0\}$ in $\CP^2$ is
$\C\times(\C\setminus\{\lambda_1,\lambda_2\})$; so the affine group, a central
extension of $F_2$ by $\la l_5\ra$, is $F_2\times\Z$.  A factorization of
$\Delta_5^2$ with this group is $s_5=(b_0,b_1,D,b_3,b_4)$, where
$D=\Delta_3^4\sigma_1^4\in B_3\subset B_5$ is the braid monodromy of the germ
$y(y-x^2-\lambda_1x^4)(y-x^2-\lambda_2x^4)$ and $b_i=g_i\sigma_{j_i}g_i^{-1}$ with
\begin{gather*}
(g_0,j_0)=(\sigma_1\sigma_2\sigma_3\sigma_1\sigma_2,4),\qquad 
(g_1,j_1)=(\sigma_1\sigma_2\sigma_1\sigma_2,3),\\
(g_3,j_3)=(\sigma_1^{-1}\sigma_2^{-1}\sigma_1^{-1}\sigma_2^{-1},3),\qquad
(g_4,j_4)=(\sigma_1^{-1}\sigma_2^{-1}\sigma_3^{-1}\sigma_1^{-1},4).
\end{gather*}
A direct computation with \eqref{eq:action} shows that
$\alpha(s_5)=\Delta_5^2$ and that $G(s_5)\cong F(u_1,u_2)\times\la z\ra$ via
$x_1,x_5\mapsto u_1$, $x_2\mapsto u_2$, $x_4\mapsto u_1u_2u_1^{-1}$,
$x_3\mapsto(u_1u_2)^{-2}z$, with components $C_1=\{1,5\}$, $C_2=\{2,4\}$ and
$L=\{3\}$.  A band $b=g\,\sigma_1g^{-1}$ whose endpoint loops $y_1,y_2$ have
colors $C_1$ and $C_2$, and which satisfies the hypothesis of
Lemma~\ref{lem:fixband}, gives a curve of degree $10$ with group
$G(s_5)/\nc{y_1^{-1}y_2}\cong\la u_1,u_2\mid u_2=wu_1w^{-1}\ra\times\Z$, which is
$\BS(2,3)\times\Z$ if $w=(u_1^{-1}u_2)^2$.  This never happens.  Indeed,
$\BS(2,3)$ maps onto the centerless group $K_{13}$ of affine maps
$x\mapsto5^ix+c$ of $\Z/13$, by $a\mapsto x+1$ and $t\mapsto5x$.  An
epimorphism of $G(s_5)/\nc{y_1^{-1}y_2}$ onto $K_{13}$ kills $z$ and maps
$g(x_1),\dots,g(x_5)$ to a tuple in the Hurwitz orbit of
$(\kappa_1,\kappa_2,(\kappa_1\kappa_2)^{-2},\kappa_1\kappa_2\kappa_1^{-1},
\kappa_1)$, colored $(C_1,C_2,L,C_2,C_1)$, for a conjugate generating pair
$(\kappa_1,\kappa_2)$, whose first two entries are equal and have colors $C_1$
and $C_2$; a computer enumeration of these Hurwitz orbits, six in number and each of
size $3120$, finds no such tuple.  Bands through $L$, bands that fail the hypothesis of
Lemma~\ref{lem:fixband}, and twisted copies of $s_5$ remain open.
\end{remark}

\subsection{Link groups}
\label{subsec:mindeg-links}
For link groups with vanishing linking numbers the lower central series
gives strong restrictions on the configuration of the components.

\begin{theorem}
\label{thm:concurrency}
Let $L\subset S^3$ be a link with $\mu\ge2$ components and vanishing pairwise
linking numbers, let $\Gamma=\pi_1(S^3\setminus L)$, and let $\overline H$ be
a Hurwitz curve with $\pi_1(\C^2\setminus H)\cong\Gamma\times\Z$.  Then
$\overline H$ has $\mu+1$ irreducible components $\overline
H_0,\dots,\overline H_\mu$, of degrees $d_0,\dots,d_\mu$, and for every point
$P\in\overline H$ there is a rational number $\kappa_P\ge0$ such that
\[
(\overline H_i\cdot\overline H_j)_P=\kappa_P\,d_id_j\qquad(i\neq j).
\]
In particular, a point at which two components meet lies on all of them,
and $\sum_P\kappa_P=1$.  More precisely, each local branch of $\overline H_i$
at $P$ meets the components $\overline H_j$, $j\neq i$, with local
intersection numbers proportional to $d_j$.
\end{theorem}

\begin{proof}
Let $G=\Gamma\times\Z$ and $V=H_1(G;\Q)$.  By
Proposition~\ref{prop:H1}, $V$ has a basis $\epsilon_0,\dots,\epsilon_\mu$
given by the classes of meridians of the components, so there are $\mu+1$
components, and the class of $l_m$ is $D=\sum_jd_j\epsilon_j$.  Let $R$ be
the kernel of the surjection $\bigwedge^2V\to(\gamma_2G/\gamma_3G)\otimes\Q$,
$a\wedge b\mapsto[\tilde a,\tilde b]$.  As the linking numbers vanish,
$\Gamma/\gamma_3\Gamma$ is free nilpotent of class $2$ on $\mu$ generators
\cite{Milnor-proc57} (see the proof of Proposition~\ref{prop:linking}), so
$(\gamma_2G/\gamma_3G)\otimes\Q=(\gamma_2\Gamma/\gamma_3\Gamma)\otimes\Q$
has dimension $\binom{\mu}{2}$ and $\dim R=\binom{\mu+1}{2}-\binom{\mu}{2}=\mu$.
Since $l_m$ is central, $R$ contains $V\wedge D$, which has dimension $\mu$;
hence $R=V\wedge D$.

Now let $P$ be a critical point of $\pr|_{\overline H}$, and let
$u=hbh^{-1}$ be the corresponding factor of $\operatorname{bmf}(\overline H)$
(\S\ref{subsec:hurwitz}), where $b$ is the braid monodromy of the germ
$(\overline H,P)$, a braid on $k$ consecutive strings.  
The closure of $b$ is the link $L_P$ of the germ:
its components are the local branches of $\overline H$ at $P$, and their
linking numbers $\ell_{aa'}$ are the local intersection numbers
$(a\cdot a')_P$.
Since $h\bigl(x^{-1}b(x)\bigr)=h(x)^{-1}u\bigl(h(x)\bigr)$, the map
$x\mapsto h(x)$ on the $k$ strings induces a homomorphism
$\phi\colon\pi_1(S^3\setminus L_P)=\la x_1,\dots,x_k\mid x_j=b(x_j)\ra\to G$
taking a meridian $y_a$ of a branch $a$ to a meridian of the component
$\overline H_{c(a)}$ containing it.  In $\pi_1(S^3\setminus L_P)$, $y_a$
commutes with its longitude, which is congruent modulo the commutator
subgroup to $\prod_{a'\neq a}y_{a'}^{\ell_{aa'}}$; hence
$[y_a,\prod_{a'\ne a}y_{a'}^{\ell_{aa'}}]$ lies in the third term of the lower
central series, and applying $\phi$ we get
\[
\epsilon_{c(a)}\wedge w_a\in R,\qquad
w_a=\sum_{j\neq c(a)}(a\cdot\overline H_j)_P\,\epsilon_j ,
\]
branches of the same component contributing nothing.  Put $i=c(a)$ and
suppose $w_a\neq0$.  Then $\epsilon_i\wedge w_a$ is a non-zero decomposable
element of $R=V\wedge D$, so $D$ lies in the plane spanned by $\epsilon_i$
and $w_a$, say $D=\beta\epsilon_i+\gamma w_a$.  Since $\mu\ge2$ and all $d_j$
are positive, $\gamma\neq0$, and $(a\cdot\overline H_j)_P=d_j/\gamma$ for all
$j\neq i$.  This proves the last assertion.  Summing over the branches of
$\overline H_i$ at $P$ gives $(\overline H_i\cdot\overline H_j)_P=A_{i,P}d_j$
with $A_{i,P}\ge0$, and by symmetry $A_{i,P}/d_i=A_{j,P}/d_j=:\kappa_P$.
Finally, distinct components meet only at points of $Z$, positively, so
$\sum_P(\overline H_i\cdot\overline H_j)_P=d_id_j$.
\end{proof}

For the curves $\overline H_{\gamma_n}$ of Theorem~\ref{thm:twisted},
$\kappa_P=1/16$ at each of the $16$ quadruple points.

\begin{corollary}
\label{cor:link-degree}
Let $L$ and $\Gamma$ be as in Theorem~\ref{thm:concurrency}, and suppose
moreover that $L$ is non-split.
\begin{enumerate}[label=(\roman*),itemsep=1pt]
\item \label{d1}
Every Hurwitz curve with group $\Gamma\times\Z$ has degree at least
$\mu+2$.
\item \label{d2}
If the meridians of such a curve go to conjugates of the meridians
$x_1,\dots,x_\mu$ of $L$ and of $(x_1\cdots x_\mu)^{-1}z$, as in
Theorem~\ref{thm:twisted}, then all $\mu+1$ components have the same degree
$e\ge2$, and the degree is $(\mu+1)e$.
\end{enumerate}
\end{corollary}

\begin{proof}
\ref{d1} The degree is at least the number $\mu+1$ of components.  In degree
$\mu+1$ every component is a section, so $\kappa_P\in\Z$ for all $P$ and
there is a single point $P$ with $\kappa_P=1$, at which all $\mu+1$ sections
meet pairwise transversally.  Sections have no tangencies with the lines of
the pencil, so $\operatorname{bmf}(\overline H)$ consists of the single
factor $\Delta_{\mu+1}^2$, and the group is $G(\Delta_{\mu+1}^2)\cong
F_\mu\times\Z$.  But $H_2(F_\mu\times\Z)\cong\Z^\mu$, whereas
$H_2(\Gamma\times\Z)\cong H_2(\Gamma)\oplus 
H_1(\Gamma)\cong\Z^{\mu-1}\oplus\Z^\mu$,
because $S^3\setminus L$ is aspherical ($L$ being non-split) and
$H_2(S^3\setminus L)\cong\Z^{\mu-1}$.  \ref{d2} By
Proposition~\ref{prop:linking}, $l_m=gz^e$ with $g$ in the center of
$\Gamma$, which lies in $\Gamma'$; so $D=e[z]=e\sum_j\epsilon_j$, that is,
$d_j=e$ for all $j$.  Then $e\ge2$ by \ref{d1}.
\end{proof}

For the unlink, which is split, $\Gamma\times\Z=F_\mu\times\Z$ is realized
in degree $\mu+1$, by $\Delta_{\mu+1}^2$.

\begin{proposition}
\label{prop:doubling-links}
Let $L$ and $\Gamma$ be as in Corollary~\ref{cor:link-degree}, and let
$\overline H$ be a Hurwitz curve with $\pi_1(\C^2\setminus H)\cong\Gamma\times\Z$
whose braid monodromy factorization is a doubling $d(s_1,s_2,s_3,s_4)$, with
$s_i\in\cS_m$.
\begin{enumerate}[label=(\roman*),itemsep=1pt]
\item \label{f1}
The components of $\overline H$ have even degrees, and $m\ge\mu+2$.
\item \label{f2}
If the meridians of $\overline H$ are as in
Corollary~\ref{cor:link-degree}\ref{d2}, then $\overline H$ has degree
$(\mu+1)e$ with $e\ge4$ even.
\end{enumerate}
In particular, if $\beta\in B_\mu$, $\mu\ge2$, is a pure braid whose closure $L$
is non-split with vanishing linking numbers, then $U_\beta$ has the smallest
degree, $4(\mu+1)$, among the doublings with group
$\pi_1(S^3\setminus L)\times\Z$ and the C-structure of
Theorem~\ref{thm:twisted}; for $L_n$ this degree is $16$.
\end{proposition}

\begin{proof}
\ref{f1} By \eqref{eq:double} and \eqref{eq:cmi}, the permutations of the
factors of $d(s_1,\dots,s_4)$ are those of the factors of the $s_i$, acting on
the first $m$ strings, their shifts, and the transpositions of $i$ and $m+i$
given by the bands $c_{m,i}$.  So, by the proof of Proposition~\ref{prop:H1},
each component of $\overline H$ consists of the strings $i$ and $m+i$ for $i$
in an orbit $O$ of the group generated by the permutations of the factors of
$s_1,\dots,s_4$, and has degree $2\lvert O\rvert$.  As there are $\mu+1$
components, $m\ge\mu+1$.  Suppose that $m=\mu+1$.  Then all orbits are points,
so the $s_i$ are pure.  A factor of $s_i$ is then the braid monodromy of the
germ of $\overline H$ at a point $P$ whose $k\ge2$ branches have degree $1$
over the base, so that they are smooth graphs
$v=f_1(u),\dots,v=f_k(u)$ lying on $k$ distinct components, those of the
strings of the germ.  By Theorem~\ref{thm:concurrency}, $k=\mu+1=m$, and any
two of the branches meet with the same multiplicity $r=4\kappa_P\ge1$.  So the
$f_a$ agree modulo $u^r$ and have distinct coefficients $c_a$ of $u^r$, and,
after subtracting their common part, the braid monodromy is that of the
points $c_au^r$, namely $\Delta_m^{2r}$, which is central.  Hence the factors of
$s_i$ are powers $\Delta_m^{2r}$ whose product is $\Delta_m^2$, so
$s_i=\Delta_m^2$, and $\pi_1(\C^2\setminus H)\cong G(\Delta_m^2)\cong
F_\mu\times\Z$ by \eqref{eq:Gdouble}.  This is not $\Gamma\times\Z$, by the
proof of Corollary~\ref{cor:link-degree}\ref{d1}, a contradiction; hence
$m\ge\mu+2$.

\ref{f2} By Corollary~\ref{cor:link-degree}\ref{d2}, all components have the
same degree $e$, which is even by~\ref{f1}; and $e=2$ would give $m=\mu+1$.
For the last assertion, $U_\beta$ has degree $4(\mu+1)$ by
Theorem~\ref{thm:twisted}\ref{b1}, and its meridians are as in
Corollary~\ref{cor:link-degree}\ref{d2} by Theorem~\ref{thm:twisted}\ref{b2},
since $\beta$ is pure.
\end{proof}

\begin{remark}
\label{rem:Ln-pencil}
Let $\Gamma_n=\pi_1(S^3\setminus L_n)$, with the C-structure of
Corollary~\ref{cor:link-degree}\ref{d2}, so that the degree is $4e$ with
$e\ge2$; Theorem~\ref{thm:twisted} gives $e=4$, the smallest value for
doublings by Proposition~\ref{prop:doubling-links}.  If $e=2$, the four
components are conics, that is, have degree $2$, and
Theorem~\ref{thm:concurrency} says that they meet only at points common to
all four, where any two of them have the same local intersection number
$4\kappa_P$; these numbers form a partition of $4$.  This is the intersection
pattern of four members of a pencil of conics, with base locus of type
$(1,1,1,1)$, $(2,1,1)$, $(2,2)$, $(3,1)$, or $(4)$.  By
Proposition~\ref{prop:multivar-fails}, such a curve could not be
algebraic.  By Theorem~\ref{thm:model}, the curves
$\overline H_{\gamma_n}$, which realize $e=4$, are obtained by a twist from four
members of a pencil of quartics, much as the case $e=2$ would require four
members of a pencil of conics.
\end{remark}

\section{Local versus global: the links \texorpdfstring{$L_n$}{L\_n} and horizontal arrangements}
\label{sec:local}

This section compares the local and global pictures.  The semiholomorphic
germs of \cite{Suciu-borro} realizing the links $L_n$ are braided surfaces,
but not Hurwitz curve germs (\S\ref{subsec:braided}).  Every C-group is the
group of a non-conical algebraic curve in a bidisk, whereas globally the
centrality obstruction applies (\S\ref{subsec:localreal}).  The horizontal
$2$-arrangements of \cite{Matei-Suciu-top00,Suciu-planes-r4} are another
family of braided-surface germs (\S\ref{subsec:arrangements}).

\subsection{Braided surfaces and Hurwitz germs}
\label{subsec:braided}
Let $f(u,v,\bar v)$ be a polynomial, monic of degree $d$ in $u$ and
holomorphic in $u$; this is the setting of semiholomorphic polynomials in
the sense of Bode \cite{Bode-jktr19}; here $v$ is the coordinate on the base
and $u$ that on the fiber.  Over each $v$ the zero set $V(f)\subset\C^2$
consists of the $d$ roots of $f(\cdot,v,\bar v)$, which collide only over the zero 
set of the $u$-discriminant, a pair of real equations in $v\in\R^2$.  When this 
zero set is discrete, as for the germs considered below, $V(f)\to\C_v$ is a 
\emph{braided surface}, in a sense slightly more general than that of 
Rudolph \cite{Rudolph-cmh83}, who requires all branch points to be
simple: a branched covering whose local monodromies are arbitrary 
braids, and whose simple branch points may be positive or negative.
A Hurwitz curve is the special case in which every singular fiber is modeled
on a complex-analytic germ with its complex orientation.

\begin{lemma}
\label{lem:expsum}
Let $\beta\in B_d$ be the braid of a Hurwitz curve over a disk $W\subset\C_v$
along $\partial W$.  Then the exponent sum $e(\beta)$ is non-negative, and
$e(\beta)=0$ only if $\beta=1$.
\end{lemma}

\begin{proof}
$\beta$ is the product of the local braid monodromies of the singular
fibers over $W$.  Each of these is quasipositive, i.e., a product of
conjugates of the standard generators $\sigma_j$: a complex-analytic germ can be
deformed, inside a small ball and without changing the boundary braid, to a
smooth complex curve whose projection has only simple, hence positive, branch
points \cite{Rudolph-top83}.  So $\beta$ is a product of $e(\beta)$
conjugates of the $\sigma_j$, and $e(\beta)=0$ forces $\beta=1$.
\end{proof}

For the semiholomorphic cubics $f_n$ of \cite[Thm.~A]{Suciu-borro}, the braid of
$V(f_n)$ around the singular value is conjugate to the braid $\gamma_n$ of
\eqref{eq:gamma}, whose closure is $L_n$.  Since $e(\gamma_n)=0$ and $\gamma_n\neq1$,
Lemma~\ref{lem:expsum} shows that $V(f_n)$ is not a Hurwitz curve over any
disk containing the singular value, even allowing several singular fibers.
More simply, the germ of $V(f_n)$ at the singular point is not complex
analytic: links of plane curve singularities with at least two components
have positive pairwise linking numbers, while those of $L_n$ vanish.

\subsection{Local and global realization}
\label{subsec:localreal}
By \cite[Thm.~7.2]{Kulikov-usp07}, every C-group is $\pi_1(D\setminus C)$ for a
bidisk $D$ and a non-singular algebraic curve $C\subset\C^2$; this applies to
$\Gamma_n=\pi_1(S^3\setminus L_n)$ with its Wirtinger C-structure, and to
$\BS(k,k+1)$ \cite[\S7.2]{Kulikov-usp07}, \cite[Cor.~0.7]{Kulikova-izv05}.
So there is no local obstruction at all; but these realizations are not
conical, i.e., they are not complements of germs at a point, and they are
unrelated to the semiholomorphic germs.  Globally, neither $\Gamma_n$ nor
$\BS(k,k+1)$ is the group of an affine Hurwitz curve, by
Proposition~\ref{prop:linking} and Corollary~\ref{cor:trivcenter}, whereas
$\BS(k,k+1)\times\Z$ and $\Gamma_n\times\Z$ are, in degrees $12$ and $16$, by
Theorems~\ref{thm:explicit} and~\ref{thm:twisted}.

\subsection{Horizontal arrangements}
\label{subsec:arrangements}
For $\sigma\in S_n$, the zero set of
\[
f_\sigma(z,w)=\prod_{j=1}^n\Bigl(z-\tfrac{j+\sigma(j)}{2}\,w
-\tfrac{j-\sigma(j)}{2}\,\bar w\Bigr),
\]
monic and holomorphic in $z$, is the horizontal $2$-arrangement $\A(\sigma)$
of $n$ transverse real planes in $\R^4$ \cite{Matei-Suciu-top00,Suciu-planes-r4}.
It is a braided-surface germ over $\C_w$, unbranched off the origin, with
pure braid monodromy; in suitable linear coordinates, every $2$-arrangement,
horizontal or not, is a germ of this kind \cite[Lemma~3.3]{Suciu-lambda}.  
On the circle $w=e^{\sqrt{-1}\,\theta}$ the roots are
$z_j=j\cos\theta+\sqrt{-1}\,\sigma(j)\sin\theta$, so the linking numbers of
its link are $\ell_{ij}=\operatorname{sgn}\bigl((i-j)(\sigma(i)-\sigma(j))\bigr)$.
As distinct branches of a complex-analytic germ link positively, $\A(\sigma)$
is a Hurwitz germ only for $\sigma=\id$, the central arrangement of the
complex lines $z=jw$, whose complement is formal \cite{Brieskorn-71}.  If
$\sigma\ne\id$ and $\operatorname{inv}(\sigma)\ge n(n-1)/4$, the boundary
braid, which is pure and non-trivial, has exponent sum
$2\sum_{i<j}\ell_{ij}=n(n-1)-4\operatorname{inv}(\sigma)\le0$, so by
Lemma~\ref{lem:expsum} $\A(\sigma)$ is not a Hurwitz curve over any disk
containing the origin.  For $\sigma\ne\id$ formality can fail: for
$\sigma=31425$ the first resonance variety $\RR_1$ is an irreducible cubic
that coincides with the tangent cone $\operatorname{TC}_1\V_1$ at $1$
\cite[Thm.~D, Prop.~9.7]{Suciu-planes-r4}, so $\C^2\setminus\A(31425)$ is not 
$1$-formal \cite[Thm.~B(1)]{Dimca-Papadima-Suciu-duke09},
and its group is not quasi-projective, as for quasi-projective groups the
components of $\V_1$ through $1$ are subtori \cite[Thm.~V.1.6]{Arapura-jag97},
\cite[\S\S4--5]{Artal-Cogolludo-Matei-gt13}; see \cite[Cor.~9.5]{Suciu-planes-r4}.
Here the departure from the algebraic constraints occurs in the tangent cone 
at the identity, whereas for the groups $\widetilde G_k$ of \S\ref{sec:alexander} 
it occurs at an isolated point away from it.

\begin{remark}
\label{rem:exotic24}
Theorem~\ref{thm:twisted} carries this phenomenon over to Hurwitz curves.  
Let $\beta\in B_5$ be the pure braid of \cite[Lemma~2.1]{Suciu-planes-r4} 
whose closure is the link of $\A(31425)$, and let
$\Gamma=\pi_1(\C^2\setminus\A(31425))$.  By Theorem~\ref{thm:twisted},
$G(U_\beta)\cong\Gamma\times\Z$ is the group of a Hurwitz curve
$\overline H_\beta$ of degree $24$, a union of six quartics, and
$\pi_1(\CP^2\setminus\overline H_\beta)\cong\Gamma\times\Z/4$.  By the
K\"unneth argument of Proposition~\ref{prop:cvar},
$\V_1(\Gamma\times\Z)=\iota(\V_1(\Gamma))$ in the coordinates of
Proposition~\ref{prop:V1-central}, with $\iota$ as in
Proposition~\ref{prop:V1-universal}.  The tangent cone of $\V_1(\Gamma)$ at $1$
is the irreducible cubic hypersurface recalled above, which is not a finite
union of linear subspaces; as a torsion-translated subtorus through $1$ is a
subtorus, some irreducible component of $\V_1(\Gamma)$ through $1$ is not a
torsion-translated subtorus.  So $\V_1(G(U_\beta))$ has an exotic component of
dimension $4$, through the trivial character.  The same argument, with $\Z/4$ in
place of $\Z$, shows that neither $\Gamma\times\Z$ nor $\Gamma\times\Z/4$ is
the fundamental group of a smooth quasi-projective variety, so that
$\overline H_\beta$ is not isotopic to an algebraic curve.
\end{remark}

\section{Questions}
\label{sec:questions}

We close with five questions: on the minimal degree of a Hurwitz curve with
group $\BS(2,3)\times\Z$ (Question~\ref{q:mindegree}) or
$\pi_1(S^3\setminus L_n)\times\Z$ (Question~\ref{q:Ln-degree}), on which
Hurwitz C-groups are quasi-projective (Question~\ref{q:qp}), on which
characteristic varieties occur (Question~\ref{q:realize}), and on the
homotopy types of complements of curves with isomorphic groups
(Question~\ref{q:homotopy}).

\begin{question}
\label{q:mindegree}
What is the minimal degree $m_{\min}$ of a Hurwitz curve $\overline H\subset\CP^2$
with $\pi_1(\C^2\setminus H)\cong\BS(2,3)\times\Z$?  Theorem~\ref{thm:explicit}
and Corollary~\ref{cor:deg3} give $4\le m_{\min}\le12$.
Proposition~\ref{prop:min12} shows that no twisted doubling of an iterated
doubling of a full twist does better than $12$, Corollary~\ref{cor:deg3}
excludes doublings of degree $6$, and Remark~\ref{rem:deg10} excludes, in the
most natural degree-$10$ construction, a single band over a degree-$5$
factorization with group $F_2\times\Z$, all bands that join the two conics
and satisfy the hypothesis of Lemma~\ref{lem:fixband}.  Is $m_{\min}<12$?  By
Proposition~\ref{prop:degconstraints}, a realization of degree $m$ has invariants
$(n,h)$ with $nh=m$ and $h\ge2$, detected by the torsion $\Z/(3^h-2^h)$ of
$A_0$.  Already degree $4$ is open: there a realization would consist of two
conics, with $(n,h)=(2,2)$, or of a cubic and a line, with $(n,h)=(1,4)$.
\end{question}

\begin{question}
\label{q:Ln-degree}
Theorem~\ref{thm:twisted} realizes $\pi_1(S^3\setminus L_n)\times\Z$ in
degree $16$.  With the product C-structure the degree is $4e$ with $e\ge2$,
and degree $16$ is the smallest one reached by Kulikov doublings
(Proposition~\ref{prop:doubling-links}).  Is there a
Hurwitz curve of degree $12$ with this group and the product C-structure, or
one of degree $8$, that is, four conics with the intersection pattern of four
members of a pencil (Remark~\ref{rem:Ln-pencil})?  More
generally, what is the minimal degree of a Hurwitz curve with group
$\pi_1(S^3\setminus\hat\beta)\times\Z$ for $\beta\in B_m$, compared with the
upper bound $4(m+1)$ and, when $\hat\beta$ is non-split with $\mu$ components
and vanishing linking numbers, the lower bound $\mu+2$?
\end{question}

\begin{question}
\label{q:qp}
Which Hurwitz C-groups are quasi-projective?  Torsion-translated
characteristic varieties are necessary \cite{Arapura-jag97,
Artal-Cogolludo-Matei-gt13} but not sufficient (Proposition~\ref{prop:multivar-fails}); 
the relative-position and single-essential-variable conditions of 
\cite[Thms.~A, B]{Dimca-Papadima-Suciu-imrn08} are further necessary conditions.  
These depend only on $\V_1$ and the multivariable Alexander polynomial, and they 
do not suffice either: by Remark~\ref{rem:raag}, the right-angled Artin groups
$(\Z*\Z^2)\times\Z^2$ and $F_3\times\Z^2$ are Hurwitz C-groups, with the 
same $\V_1$, a coordinate subtorus of dimension~$3$
\cite[Prop.~11.5]{Dimca-Papadima-Suciu-duke09}, and the same Alexander
polynomial, namely~$1$; but only the second is quasi-projective.  Indeed, 
the tangent space at~$1$ of $\V_1\bigl((\Z*\Z^2)\times\Z^2\bigr)$ is a
$3$-dimensional subspace of $H^1$ on which the cup product does not vanish,
whereas for a quasi-projective group such a subspace is $0$-isotropic, by the
isotropy condition of \cite[Thm.~C(1)]{Dimca-Papadima-Suciu-duke09}.  
Is there a Hurwitz C-group satisfying all of these conditions that is not
quasi-projective?
\end{question}

\begin{question}
\label{q:realize}
Call an irreducible component of $\V_1$ \emph{exotic}\/ if it is not a
torsion-translated subtorus; quasi-projective groups have none
\cite{Arapura-jag97,Artal-Cogolludo-Matei-gt13}.  For a Hurwitz curve 
with components of degrees $d_1,\dots,d_r$, the exotic components lie in 
the hypersurface $\prod_c\lambda_c^{d_c}=1$ (Proposition~\ref{prop:V1-central}), 
and, $\V_1$ being defined over $\Z$, their union is stable under
$\operatorname{Gal}(\bar\Q/\Q)$.  Without a bound on the degree, every
subvariety of $(\C^*)^n$ defined over $\Z$, $n\ge2$, occurs, embedded in
$(\C^*)^{n+1}$ (Proposition~\ref{prop:V1-universal}).  Which exotic
components occur for curves of bounded degree, or with prescribed
singularities?  The curves $\overline H_k$ of degree $12$ have the isolated
exotic points $(\zeta_k,\zeta_k^{-2})$, on $\lambda_1^8\lambda_2^4=1$.  What is
the smallest degree of a Hurwitz curve whose characteristic variety has an
exotic component, or a positive-dimensional one?  The curve of
Remark~\ref{rem:exotic24} has degree $24$ and an exotic component 
of dimension $4$.
\end{question}

\begin{question}
\label{q:homotopy}
Are the affine complements $\C^2\setminus Q_k$ of the curves
$\overline Q_k$ of Proposition~\ref{prop:samegroup} homotopy equivalent?
More generally, are there Hurwitz curves of the same degree and with the
same singularities, with C-isomorphic groups, whose affine complements are
not homotopy equivalent? The complement of an affine Hurwitz curve is 
homotopy equivalent to the $2$-complex of its braid monodromy presentation
\cite[Thm.~4.3]{Kulikov-usp07}, as Libgober proved for algebraic curves
(see \cite[Thm.~4.2]{Cohen-Suciu-cmh97}), so its Euler characteristic is 
determined by the degree and the singularities; when the group is $\pi\times\Z$, 
whose classifying space is a $2$-complex times a circle, Schanuel's lemma shows
that $\pi_2$ of the complement is stably free.  For the trefoil group itself
there are stably free, non-free relation modules \cite{Dunwoody-blms72}, and
infinitely many homotopy types of $2$-complexes with Euler characteristic
one \cite[Thm.~4.5]{Harlander-Jensen-agt06}.
\end{question}

\subsection*{Acknowledgments}
Parts of this work were carried out with the assistance of Claude Opus~5.5
(Anthropic), an AI system, which was used for exact computer checks of
braid identities and group presentations, for checking citations against
the original sources, for suggestions on the exposition, and for producing
the code of the figures.  All mathematical arguments and computations have
been independently verified by the author, who takes full responsibility
for the contents of this paper.

\raggedbottom


\providecommand{\MR}{}
\renewcommand{\MR}[1]{\href{https://mathscinet.ams.org/mathscinet-getitem?mr=#1}{MR~#1}}

\end{document}